\documentclass[10pt]{amsart}
\usepackage{amsmath,amssymb,amsthm}
\newtheorem{thm}{Theorem}[section]
\newtheorem{lem}[thm]{Lemma}
\newtheorem{cor}[thm]{Corollary}
\newtheorem{prop}[thm]{Proposition}
\newtheorem*{Prop}{Proposition}		
\newtheorem*{sublemma}{Sublemma} 
\newtheorem*{bprime}{Theorem~\ref{thm:3.6}(b)$'$}
\newtheorem*{tenprime}{Lemma~\ref{lem:3.10}$'$}
\newtheorem{DEF}[thm]{Definition}
\theoremstyle{remark}                  
\newtheorem*{rem}{Remark}
\newtheorem*{rems}{Remarks}
\newtheorem*{Def}{Definition}
\newtheorem*{Com}{Comment}
\newtheorem*{Problem}{Problem}
\newtheorem*{Problemprime}{Problem$'$}
\def\EE{{\mathbb E}}
\def\EEo{{\EE_\omega}}
\def\complex{{\mathbb C}}
\def\nat{{\mathbb N}}
\def\que{{\mathbb Q}}
\def\real{{\mathbb R}}
\def\zed{{\mathbb Z}}
\def\A{{\mathcal A}}
\def\B{{\mathcal B}}
\def\J{{\mathcal J}}
\def\K{{\mathcal K}}
\def\M{{\mathcal M}}
\def\N{{\mathcal N}}
\def\P{{\mathcal P}}
\def\R{{\mathcal R}}
\def\cS{{\mathcal S}}
\def\U{{\mathcal U}}
\def\bj{{\mathbf j}} 
\def\T{{\tau}}
\def\ep{\varepsilon}
\def\lam{\lambda}
\def\om{\omega}
\def\bom{\underline{\omega}}
\def\wtilde{\widetilde}
\def\Moo{{M_{00}}}
\def\Aut{\operatorname{Aut}}
\def\cb{\operatorname{cb}}
\def\Col{\operatorname{Col}}
\def\Ext{\operatorname{Ext}}
\def\id{\operatorname{id}}
\def\op{\operatorname{op}}
\def\Rad{\operatorname{Rad}}
\def\RadCp{{\Rad C_p}}
\def\Row{\operatorname{Row}}
\def\Ba{\mathop{\B_a}\nolimits}
\def\csim{\mathop{\buildrel C\over\sim}}
\def\defeq{\mathop{\buildrel \text{def}\over =}}
\def\chix{{\raise.5ex\hbox{$\chi$}}}
\def\oco{\overline{\text{co}}}
\def\ellpp{{\ell^{1 +}}}
\def\bone{\text{\bf 1}}
\def\ccarrow{\overset{cc}{\hookrightarrow}}
\title{Banach embedding properties\\ of non-commutative $L^p$-spaces}
\author{U. Haagerup, H.P. Rosenthal and F.A. Sukochev}
\address{U.H.: Department of Mathematics and Computer Science,
Odense University, DK-5230 Odense M, Denmark}
\email{haagerup@imada.sdu.dk}
\address{H.R.: Department of Mathematics, The University of Texas at Austin,
Austin, TX 78712 USA}
\email{rosenthl@math.utexas.edu}
\address{F.S.: Department of Mathematics and Statistics, School of
Informatics and Engineering, The Flinders University of South Australia,
Bedford Park, 5042, SA, Australia}
\email{sukochev@ist.flinders.edu.au}

\begin{document}
\maketitle
\pagestyle{headings}		
\markboth{}{}	
\begin{abstract}
Let $\N$ and $\M$ be von Neumann algebras. 
It is proved that $L^p(\N)$ does not Banach embed in $L^p(\M)$ for $\N$ 
infinite, $\M$ finite, $1\le p<2$. 
The following considerably stronger result is obtained (which implies this, 
since the Schatten $p$-class $C_p$ embeds in $L^p(\N)$ for $\N$ infinite). 

{\bf Theorem.} 
{\it Let $1\le p<2$ and let $X$ be a Banach space with a spanning set 
$(x_{ij})$ so that for some $C\ge 1$, 
\begin{itemize}
\item[(i)] any row or column is $C$-equivalent to the usual $\ell^2$-basis, 
\item[(ii)] $(x_{i_k,j_k})$ is $C$-equivalent to the usual $\ell^p$-basis, 
for any $i_1<i_2<\cdots$ and $j_1<j_2<\cdots$. 
\end{itemize}
Then $X$ is not isomorphic to a subspace of $L^p(\M)$, for $\M$ finite.} 
Complements on the Banach space structure of non-commutative $L^p$-spaces 
are obtained, such as the $p$-Banach-Saks property and characterizations 
of subspaces of $L^p(\M)$ containing $\ell^p$ isomorphically. 
The spaces $L^p(\N)$ are classified up to Banach isomorphism, for $\N$ 
infinite-dimensional,  hyperfinite and semifinite, $1\le p<\infty$, $p\ne 2$. 
It is proved that there are exactly thirteen isomorphism types; the 
corresponding embedding properties are determined for $p<2$ via an eight 
level Hasse diagram.
It is also proved for all $1\le p<\infty$ that $L^p(\N)$ is completely 
isomorphic to $L^p(\M)$ if $\N$ and $\M$ are the algebras associated 
to free groups, or if $\N$ and $\M$ are injective factors of type III$_\lambda$ 
and III$_{\lambda'}$ for $0<\lambda$, $\lambda'\le 1$. 
\end{abstract}
\bigskip

{\parskip=0pt
\centerline{\sc Contents}

\begin{itemize}
\item[\S1.] Introduction.
\item[\S2.] The modulus of uniform integrability 
and weak compactness in $L^1(\N)$.
\item[\S3.] Proof of the Main Theorem.
\item[\S4.] Improvements to the Main Theorem.
\item[\S5.] Complements on the Banach/operator space 
structure of $L^p(\N)$-spaces.
\item[\S6.] The Banach isomorphic classification of the spaces $L^p(\N)$
for $\N$ hyperfinite semi-finite.
\item[\S7.]  $L^p(\N)$-isomorphism results for $\N$ type III hyperfinite or 
a free group von Neumann algebra.
\item[] References
\end{itemize}
}

\section{Introduction} 
\setcounter{equation}{0}

Let $\N$ be a finite von Neumann algebra and $1\le p<2$. 
Our main theorem yields that $C_p$ is not linearly isomorphic to a subspace 
of $L^p(\N)$ (where $C_p$ denotes the Schatten $p$-class). 
It follows immediately that for any infinite von~Neumann algebra $\M$, 
$L^p(\M)$ is not isomorphic to a subspace of $L^p(\N)$, since $C_p$ is then 
isomorphic to a subspace of $L^p(\M)$. 
(It is proved in \cite{S1} that also $C_p$ does not embed in $L^p(\N)$ 
for any $2<p<\infty$.) 

For $\N$ a semi-finite von-Neumann algebra and $\tau$ a faithful normal 
semi-finite trace on $\N$, $L^p(\tau)$ denotes the non-commutative $L^p$ 
space associated with $(\N,\tau)$ (see e.g., \cite{FK}). 
The particular choice of trace $\tau$ is unimportant, for if $\beta$ is 
another such trace, $L^p(\beta)$ is isometric to $L^p(\tau)$. 
We also denote this (isometrically unique) Banach space by $L^p(\N)$. 

Given $C\ge1$ and non-negative reals $a$ and $b$, let $a \csim b$ denote 
the equivalence relation $\frac1C a \le b\le Ca$. 
Sequences $(x_j)$ and $(y_j)$ in Banach spaces $X$ and $Y$ respectively all 
called $C$-equivalent if 
\begin{equation}\label{eq:1.1}
\Big\| \sum_{i=1}^n \alpha_i x_i\Big\| \csim \Big\|\sum_{i=1}^n 
\alpha_i y_i\Big\|\ \text{ for all $n$ and scalars } \alpha_1,\ldots,\alpha_n
\ .
\end{equation}
(Equivalently, there exists an invertible linear map 
$T:[x_i]\to [y_i]$ with $\|T\|,\|T^{-1}\| \le C$, where $[x_i]$ denotes 
the closed linear span of $(x_i)$.)  
$(x_j)$ is called {\it unconditional\/} if there is a constant $u$ so that 
for any $n$ and scalars $c_1,\ldots,c_n$ and $\ep_1,\ldots,\ep_n$ with 
$|\ep_i| =1$ for all $i$, 
$\| \sum_{i=1}^n \ep_i c_i x_i\|\le u\|\sum c_i x_i\|$ 
(then one says $(x_j)$ is $u$-unconditional). 
The usual $\ell^p$-basis refers to the unit vector basis $(e_j)$ of $\ell^p$, 
where $e_j(i) = \delta_{ji}$ for all $i$ and $j$.

Our main result goes as follows. 

\begin{thm}\label{thm:1.1} 
Let $\N$ be a finite von~Neumann algebra, $1\le p<2$, and let $(x_{ij})$ 
be an infinite matrix in $L^p(\T)$ where $\T$ is a fixed faithful normal 
tracial state on $\N$. 
Assume for some $C\ge1$ that every row and column of $(x_{ij})$ is 
$C$-equivalent to the usual $\ell^2$-basis 
and that $(x_{i_k,j_k})_{k=1}^\infty$ is unconditional, whenever $i_1<i_2<
\cdots$ and $j_1<j_2<\cdots$. 
Then there exist $i_1<i_2<\cdots$ and $j_1<j_2<\cdots$ so that setting 
$y_k = x_{i_k,j_k}$ for all $k$, then 
\begin{equation}\label{eq:1.2}
\lim_{n\to\infty}  n^{-1/p} \Big\|\sum_{i=1}^n  y'_i\Big\|_{L^p(\T)} =0
\end{equation}
for all subsequences $(y'_k)$ of $(y_k)$.
\end{thm}

\begin{cor}\label{cor:1.2} 
Let $p$ and $\N$ be as in \ref{thm:1.1}. 
Let $X$ be a Banach space spanned by an infinite matrix of elements 
$(x_{ij})$ so that for some $\lam \ge1$, 
\begin{itemize}
\item[(i)] every row and column of $(x_{ij})$ is $\lam$-equivalent to 
the usual $\ell^2$ basis
\item[(ii)] $(x_{i_n,j_n})_{n=1}^\infty$ is $\lam$-equivalent to the usual 
$\ell^p$-basis, for all $i_1<i_2<\cdots$ and $j_1<j_2<\cdots$. 
\end{itemize}
Then $X$ is not Banach isomorphic to a subspace of $L^p(\T)$.
In particular, $C_p$ does not embed in $L^p(\T)$. 
\end{cor}

The Corollary yields its final statement since the 
standard matrix units $(x_{ij})$ for $C_p$ satisfy (i) and (ii) with 
$\lam=1$. 

To see why \ref{thm:1.1} $\implies$ \ref{cor:1.2}, 
suppose to the contrary that $T:X\to X'
\subset L^p(\T)$ were an isomorphic embedding, where $X$ is as in 1.2. 
Then $(Tx_{ij})$ satisfies the hypotheses of 1.1 with $C=\lam\|T\|\, 
\|T^{-1}\|$. 
However if $(i_k),(j_k)$ satisfies the conclusion of Theorem~1.1, 
$(Tx_{i_k,j_k})$ and hence $(x_{i_k,j_k})$ cannot be equivalent to the 
usual $\ell^p$-basis, a contradiction. 

Let $\RadCp$ denote the ``Rademacher unconditionalized version'' of 
$C_p$ $(1\le p<\infty)$. 
That is, letting $(r_{ij})$ be an independent matrix of $\{1,-1\}$-valued 
random variables with $P(r_{ij}=1) = P(r_{ij}=-1) =\frac12$ for all $i,j$, 
and letting $(c_{ij})$ be a matrix of scalars with only finitely  
many non-zero terms, then 
\begin{equation}\label{eq:1.3}
\|(c_{ij})\|_{\Rad_{C_p}} = \EEo\|(r_{ij}(\om)c_{ij})\|_{C_p}\ .
\end{equation}

\begin{cor}\label{cor:1.3}
Let $p$ and $\N$ be as in \ref{thm:1.1}. 
Then $\RadCp$ is not isomorphic to a subspace of $L^p(\T)$.
\end{cor}

\begin{proof} 
The standard matrix units basis $(x_{ij})$ of $\RadCp$ also satisfies 
the hypotheses of Corollary~\ref{cor:1.2} with $\lam=1$.
\end{proof}

Corollary~\ref{cor:1.3} yields new information in the classical, 
commutative case of $L^p$. 
(Throughout, $L^p$ refers to $L^p$ on the unit interval, endowed with 
Lebesgue measure; i.e., $L^p = L^p(\N)$ where $\N = L^\infty$ acting 
on $L^2$ via multiplication.) 
This also reveals a remarkable difference in the structure of 
$L^p$-spaces, $p<2$ or $p>2$, 
for $\RadCp$ is {\em isometric\/} to a subspace of $L^p$ 
for $2<p<\infty$ (cf.\ Theorem~5 of \cite{L-P}). 
Also, let us note that $\RadCp$ is isometric to a subspace of $L^p$ $(C_p)$
for $1\le p<2$, so we obtain an unconditionalized version of $C_p$ in 
$L^p(\M)$ which also does not embed in $L^p(\N)$, 
for $\N$ finite, where $\M = L^\infty \otimes B(H)$. 
(Throughout, $L^p(X)$ refers to the Bochner-Lebesgue space $L^p(X,m)$, 
where $m$ is Lebesgue measure.) 

It is a classical result of C.A.~McCarthy that $C_p$ does not ``locally'' 
embed in $L^p$, for $1\le p<\infty$ \cite{McC}. 
Corollary~\ref{cor:1.2} yields an ``infinite'' dimensional proof of this 
result for $1\le p<2$, as well as the apparently new discovery that also 
$\RadCp$ does not locally embed in $L_p$ for these $p$. 
To see this, we give the following.

\begin{Def}
Let $1\le p<\infty$, $n\in\nat$, and $\lam\ge1$. 
A finite-dimensional Banach space $X$ is called a $\lam$-$GC_p^n$-space 
provided there is an $(n\times n)$-matrix $(x_{ij})$ spanning $X$ so that 
\begin{itemize}
\item[(i)] any row and column of $(x_{ij})$ is $\lam$-equivalent to the 
usual $\ell_n^2$-basis
\item[(ii)] $(x_{i_k,j_k})_{k=1}^m$ is $\lam$-equivalent to the usual 
$\ell_m^p$ basis for any $m$, 
$$1\le i_1< \cdots < i_m\le n\ \text{ and }\ 1\le j_1<j_2<\cdots < j_m\le n\ .$$
\end{itemize}
An infinite-dimensional space $X$ is called a $\lam$-$GC_p$-space 
provided it admits a spanning matrix $(x_{ij})$ satisfying (i) and (ii) 
of Corollary~\ref{cor:1.2}; finally $X$ is called a $GC_p$-space if it 
is a $\lam$-$GC_p$-space for some $\lam\ge1$.
\end{Def}

$C_p^n$ refers to the $n^2$-dimensional Schatten $p$-class consisting of 
$n\times n$ matrices in the $C_p$ norm; ``$G$'' stands for 
``Generalized''. 
For example, $\Rad C_p^n$ is a 1-$GC_p^n$ space. 
The next result yields that $\lam$-$GC_p^n$-spaces cannot be uniformly 
embedded in $L^p$, hence in particular, we recapture the classical 
fact mentioned above that $L^p$ does not contain $C_p^n$'s uniformly. 
(For isomorphic Banach spaces $X$ and $Y$, 
$d(X,Y) = \inf \{\|T\|\,\|T^{-1}\| : T$ 
is  a surjective isomorphism from $X$ to $Y\}$).  

\begin{cor}\label{cor:1.4}
Let $1\le p<2$ and $\lambda\ge1$. 
Define: 
\begin{equation*}
\beta_{n,\lam} = \inf \{d(X,Y) :X\text{ is a  $\lam$-$GC_p^n$-space and } 
Y\subset L^p\}\ .
\end{equation*}
Then $\lim_{n\to\infty} \beta_{n,\lam} =\infty$. 
\end{cor}

\begin{proof}
Suppose this were false. 
Then we could choose $\lam\ge1$ and $X_1,X_2,\ldots$ subspaces of $L^p$ 
so that $X_n$ is a $\lam$-$GC_p^n$-space for all $n$. 
Choose then $(x_{ij}^n)$ an $n\times n$ matrix of elements of $X_n$, 
satisfying (i) and (ii) of the definition, for all $n$. 
Let $\Moo$ denote the linear space of all infinite matrices of scalars 
with only finitely many non-zero entries. 
Let $U$ be a free ultrafilter on $\nat$. 
Define a semi-norm $\|\cdot\|$ on $\Moo$ by 
\begin{equation}\label{eq:1.4} 
\|(c_{ij})\| = \lim_{n\in U} \| \sum c_{ij} x_{ij}^n\|\ .
\end{equation}

It is easily checked that $\|\cdot\|$ is indeed a semi-norm; let $W$ be 
its null space; $W= \{(c_{ij}) \in \Moo :\|(c_{ij})\| =0\}$, and 
let $X$ denote the completion of $(\Moo,\|\cdot\|)/W$. 
It follows easily that $X$ is a $\lam$-$GC_p$-space. 
By standard ultraproduct techniques, it follows that $X$ is finitely 
representable in $L^p$. 
But then (since ultraproducts of (commutative) 
$L^p(\mu)$ spaces are (commutative) $L^p(\nu)$ spaces 
and any separable subspace of an $L^p(\nu)$ space is isometric to a 
subspace of $L^p$), 
$X$ isometrically embeds in $L^p$. 
This contradicts Corollary~\ref{cor:1.2}. 
\end{proof}

\begin{rem}
Theorem~1.1 may easily be extended to the case of general finite 
von Neumann algebras $\N$, and not just the finite, $\sigma$-finite 
ones covered by its statement. 
Corollaries~1.2 and 1.3 also hold in this setting, as well as the general 
formulations of Theorems~4.1 and 4.2. 
Indeed, in general, one has that $L^p(\N)$ is isometrically isomorphic to 
$L^p(\T)$ for some semi-finite faithful normal trace $\T$ on $\N$. 
Let $(x_{ij})$ be a matrix of elements of $L^p(\T)$ satisfying the 
assumptions of Theorem~1.1, and let $P$ be the supremum of all the 
support projections of $x_{ij}$ and $x_{ij}^*$, $i,j=1,2,\ldots$.  
Then $P$ is a $\sigma$-finite projection in $\N$, and thus $P\N P$ is 
both finite and $\sigma$-finite. 
Moreover all the $x_{ij}$'s belong to $L^p(P\N P,\T') = PL^p(\N,\T)P$, 
where $\T' = \T|P\N P$. 
But in turn, $L^p(P\N P,\T')$ is isometrically isomorphic to 
$L^p(P\N P,\T'')$ for some faithful finite normal trace $\T''$ on $P\N P$. 
This reduces the proof of Theorem~1.1 in the case of general finite 
von Neumann algebras, to those with a finite trace. 
\end{rem}

We now give a description of the results and proof-order of the paper. 

If a matrix satisfies the hypotheses of Theorem~\ref{thm:1.1}, 
then every row and column has the property that the $p^{th}$ powers 
of absolute values of the terms form a uniformly integrable sequence. 
We develop the basic technical tools to explain and exploit this, in 
Section~2, through the device of the $p$-modulus of an element of $L^p(\N)$ 
with respect to a normal tracial state $\T$ on $\N$. 
We give several useful inequalities for this modulus in Lemma~\ref{lem:2.3}. 
Although many of these can be obtained from the literature (e.g., 
\cite{FK}), we give full proofs for the sake of completeness. 
We also obtain equivalences for relative weak compactness in 
$L^1(\N)$ in terms of uniform integrability in Proposition~\ref{prop:2.2}, 
and a useful non-commutative truncation equivalence for general $p$, 
in Corollary~\ref{cor:2.7}.

We give technical information concerning general unconditional sequences 
in $L^p(\N)$ for $p<2$ in Lemmas~\ref{lem:3.1}--\ref{lem:3.3}, 
yielding in particular the following definitive equivalences 
obtained in Corollaries~\ref{cor:3.4} and \ref{cor:3.5}. 
{\it Let $(f_n)$ be a bounded unconditional sequence in $L^p(\N)$.  
Then the following are equivalent.} 
\begin{itemize}
\item[1.] {\it $(f_n)$ has no subsequence equivalent to the 
usual $\ell^p$ basis.}
\item[2.] {\it $(|f_n|^p)$ is uniformly integrable.}
\item[3.] {\it $\lim_{n\to\infty} n^{-1/p} \|\sum_{i=1}^n f'_i\|_{L^p(\T)}=0$ 
for all subsequences $(f'_n)$ of $(f_n)$.}
\end{itemize}

The proof of Theorem~\ref{thm:1.1} is then completed, using the standard 
ultraproduct construction of the finite ultrapower of a finite 
von Neumann algebra $\N$, and a result giving the connection between its 
associated $L^p$ space and the Banach ultrapower of $L^p(\N)$ 
(Lemma~\ref{lem:4.1}). 

Section 4 yields results considerably stronger than Theorem~\ref{thm:1.1}. 
The arguments here do not use the ultraproduct construction in Section~3, and 
are thus more elementary (but also more delicate). 
Theorem~\ref{thm:3.5}  
gives the following result (which immediately implies Theorem~\ref{thm:1.1}). 

{\it If a semi-normalized matrix in $L^p (\N)$ is such that all columns and 
``generalized'' diagonals are unconditional while all rows are 
$u$-unconditional for some fixed $u$, then three alternatives occur:}
{\it Either some column has an $\ell^p$-subsequence, or $\ell_n^p$'s are 
finitely represented in the terms of the rows, or the matrix has a 
``generalized diagonal'' $(y_k)$ satisfying \eqref{eq:1.2} 
of Theorem~\ref{thm:1.1}.}

Using results from Banach space theory, we obtain in Theorem~\ref{thm:3.6} 
that if $p=1$ or if $p>1$ and $\N$ is hyperfinite, the unconditionality 
assumption in \ref{thm:3.5} may be dropped. 
The case $p>1$ also uses recent non-commutative martingale 
inequalities (see \cite{SF}, \cite{PX1}). 
The case $p=1$ uses techniques from \cite{R1}, which yield results 
for sequences in the preduals of arbitrary von~Neumann algebras 
which may be independent interest (see Lemmas~\ref{lem:3.12} 
and \ref{lem:3.13}). 
The proof in this case  also requires an apparently new elementary 
finite disjointness result (Lemma~\ref{lem:3.14}B). 

Section 5 contains rather quick applications of our main results and 
the techniques of their proofs. 
For example, Proposition~\ref{prop:5.1} asserts that neither the Row 
nor Column operator spaces completely embed in the predual of a finite 
von~Neumann algebra; this is a quick consequence of our main result. 
Theorem~\ref{thm:5.4} shows that for $1\le p<2$ and $\N$ finite, a 
subspace of $L^p(\N)$ contains $\ell_n^p$'s uniformly iff it contains 
an almost disjointly supported sequence (which of course is then almost 
isometric to $\ell^p$), extending the previously known commutative 
case \cite{R2}). 
We give the concepts of the $p$-Banach-Saks and strong $p$-Banach-Saks 
properties in Definition~\ref{def:5.5}, and extend the classical 
results of Banach-Saks \cite{BS} and Szlenk \cite{Sz} in 
Proposition~\ref{prop:5.6}. 
This result also yields that for $p$ and $\N$ as above, a weakly 
null sequence in $L^p(\N)$ has the property that every subsequence 
has a strong $p$-Banach-Saks subsequence if and only if the $p^{th}$ 
powers of absolute values of its terms are uniformly integrable.

The main result of Section 6 shows that there are precisely thirteen 
Banach isomorphism types among the spaces $L^p(\N)$ for  $\N$ hyperfinite 
semi-finite, $1\le p<\infty$, $p\ne 2$. 
The embedding properties of the various types for $p<2$ 
are given in an eight-level 
Hasse diagram, in Theorem~\ref{thm:6.2}. 
This work completes the classification and embedding properties of the 
type~I case given in \cite{S2}. 
The main work in establishing this Theorem is found in the non-embedding 
results given in Theorems~\ref{thm:6.3} and \ref{thm:6.9}; we also give a 
new proof of a non-embedding result in the type~I case, established in 
\cite{S2}, in our Proposition~\ref{prop:6.5}.
The most delicate of these is Theorem~\ref{thm:6.9}, which yields that 
if $\M$ is a type~II$_\infty$ von-Neumann algebra, and $L^p(\M)$ embeds 
in $L^p(\N)$, then also $\N$ must have a type~II$_\infty$ or 
type~III summand ($1\le p<2$). 
Of course this reduces directly to the case where $\M$ is the hyperfinite 
type~II$_\infty$ factor; the proof requires our Theorem~\ref{thm:3.5}, 
and also rests upon recent discoveries of M.~Junge \cite{J} and 
Pisier-Xu \cite{PX2}.

Our methods do not cover the following case, which remains a 
fascinating open problem: 
Is it so that the predual of a type~III von-Neumann algebra does not 
Banach embed in the predual of one of type~II$_\infty$?
In fact, we do not know if the predual of the injective 
type~II$_\infty$~factor can be Banach isomorphic to the predual of 
an injective type~III-factor. 
We show in Theorem~7.2 that  such factors cannot in general be distinguished 
by the Banach space isomorphism class (or even operator space isomorphism 
class) of their preduals. 
Letting $R_\lambda$ denote the Powers injective factor of type 
III$_\lambda$ and $R_\infty$ denote the Araki-Woods injective factor 
of type III$_1$, we show that $(R_\lambda)_*$ is completely isomorphic 
to $(R_\infty)_*$ for all $0<\lambda <1$. 
(For a von~Neumann algebra $\N$, $\N_*$ denotes its predual, also denoted 
here by $L^1(\N)$.) 
Thus there are uncountably many isomorphically distinct injective factors, 
all of whose preduals are completely isomorphic. 
We also show in Theorem~7.2 that there are uncountably many isomorphically 
distinct injective type~III$_0$-factors, all of whose preduals are 
completely isomorphic to $(R_\infty)_*$. 

We show in Theorem 7.3 that the famous open isomorphism problem for free 
group von~Neumann algebras cannot be resolved by the Banach (or even 
operator) space structure of the predual. 
Namely, we prove that the preduals of the $L(F_n)$'s are all completely 
isomorphic, for $2\le n\le\infty$, where $F_n$ is the free group on 
$n$ generators and $L(F_n)$ its associated von~Neumann algebra. 
This extends the result of A.~Arias \cite{Ar}, showing that the $L(F_n)$'s 
themselves are completely isomorphic as operator spaces. 
The proof of Theorem~7.3 relies basically on the deep result of 
D.~Voiculescu that $L(F_\infty) \cong M_k(L(F_\infty))$ as 
von~Neumann algebras, for $k= 2,3,\ldots$ (cf.\ \cite{Vo} 
or \cite{VDN}). 

The results in Section 7 also extend to the case of the non-commutative 
spaces $L^p(\N)$, for $1<p<\infty$ (see Theorem~7.5). 
These isomorphism results (as well as the ``positive'' isomorphism 
results in Section~6) rely on the operator space version of the so-called 
Pe{\l}czy\'nski decomposition method (see Lemma~6.13). 
Thus, one actually shows for von~Neumann algebras $\N$ and $\M$, that each 
of the spaces $L^p(\N)$ and $L^p(\M)$ is completely isometric to a 
completely contractively complemented subspace of the other, and also 
(e.g., in the free group case $\M=L(F_\infty))$, that say $L^p(\M)$ also 
has the property that 
$(L^p(\M) \oplus \cdots \oplus L^p(\M)\oplus\cdots)_{\ell^p}$ 
completely contractively factors through $L^p(\M)$, which then implies 
the operator space isomorphism of these two spaces. 
Thus the proofs of these operator space isomorphism results are actually 
based on natural isometric embedding properties of the 
$L^p(\N)$ spaces themselves.

\section{The modulus of uniform integrability and weak compactness 
in $L^1(\N)$}
\setcounter{equation}{0}

Let $\N$ be a finite von Neumann algebra, acting on a Hilbert space $H$. 
Let $\P = \P(\N)$ denote the set of all (self-adjoint) projections in $\N$. 
We shall assume that $\N$ is endowed with a faithful normal tracial state 
$\T$, which is {\it atomless\/}. 
That is, for all $P\in\P$ with $P\ne 0$, there is a $Q\le P$, $Q\in\P$, 
with $0<\T(Q) <\T(P)$. 
(Equivalently, $0\ne Q\ne P$, since $\T$ is faithful.) 

These assumptions cause no loss in generality. 
Indeed, if $\N$ has a faithful normal trace $\gamma$, then simply replace 
$\N$ by $\tilde\N = \N\bar\otimes L^\infty$, where $\tilde\N$ is equipped 
with the atomless trace $\gamma = \T\otimes m$, with $m$ the trace on 
$L^\infty$ given by integration with respect to Lebesgue measure on $[0,1]$. 
$\N$ is ($*$-isomorphic to) a subalgebra of $\tilde\N$, and hence $L^p(\N)$ 
is isometric to a subspace of $L^p(\tilde\N)$, so we may as well assume 
our space $X$ in Theorem~\ref{thm:1.1} is already contained in 
$L^p(\tilde \N)$. 

Now if $\M\subset \N$ is a MASA, it follows easily that also $\T|\M$ is 
atomless. 
Indeed, were this false, we could choose $P\ne0$, $P\in\M$ so that 
$0\le Q\le P$, $Q\in\M$ implies $Q=0$ or $Q=P$. 
But then choosing $Q\in \P(\N)$, $0\le Q\le P$ with 
$0<\T(Q) <\T(P)$, we obtain that if $\wtilde \M$ is the von~Neumann algebra 
generated by $\M$ and $Q$, $\wtilde \M$ is also commutative and 
$\wtilde \M\ne\M$, a contradiction. 


\begin{DEF}\label{def:2.1}
Given $f\in \N_* = L^1(\T)$, we define the modulus of uniform integrability 
of $f$ as the function on $\real^+$, $\ep \to \om (f,\ep)$ given by  
\begin{equation}\label{eq:2.1}
w(f,\ep) = \sup \{\T(|fP|),\ P\in\P,\ \T(P)\le\ep\}\ .
\end{equation}
We also define  the lower modulus of $f$, $\ep\to\bom (f,\ep)$, as 
\begin{equation}\label{eq:2.2} 
\bom (f,\ep ) = \sup \{|\T(fP)| : P\in \P ,\T(P) \le \ep\}\ .
\end{equation}
\end{DEF}

To handle the case $p\ne1$ in our Main Theorem, we also use the following 
$p$-moduli. 
(When $\T$ is fixed, we set $\|f\|_p = \|f\|_{L^p(\T)} = (\T(|f|^p))^{1/p}$. 
Also, for $f\in \N$, we set $\|f\|_\infty = \|f\|_{\N}$.) 

\begin{DEF}\label{def:2.2}
Let $0<p<\infty$ and $f\in L^p(\T)$. 
The $p$-modulus of $f$, $\om_p(f,\cdot)$, the symmetric $p$-modulus of $f$, 
$\om_p^s (f,\cdot)$, and the spectral $p$-modulus of $f$, $\tilde\om_p(f,
\cdot)$ are given, for $0\le \ep \le1$, by 
\begin{align}
\om_p(f,\ep) &= \sup \{\|fP\|_p :P\in \P,\ \T(P)\le\ep\}\ ,\label{eq:2.3}\\
\om_p^s(f,\ep) & = \sup\{ \|PfP\|_p : P\in\P,\ \T(P) \le\ep\}\ ,\label{eq:2.4}\\
\tilde\om_p(f,\ep)&=\sup
\biggl\{\Big(\int_{(r,\infty)} t^pd(\T\circ E_{|f|}(t))  \Big)^{1/p} : 
\T\circ E_{|f|} ((r,\infty))\le\ep\biggr\}\label{eq:2.5}
\end{align}
where for $g$ self-adjoint, $E_g$ denotes the spectral measure for $g$.
\end{DEF}

It is trivial that all these moduli are increasing (i.e., non-decreasing) 
functions on $\real^+$, which are continuous at $0$, thanks to the 
assumption that $f\in L^p(\T)$. 
Actually, the assumption that $\T$ is atomless yields that $\om_p(f,\cdot)$, 
$\bom(f,\cdot)$ and $\om_p^s(f,\cdot)$ are {\it absolutely continuous\/} 
on $[0,1]$. 

We now give some basic properties of these moduli. 
The most important of these is that several of them reduce to the uniform 
integrability modulus given in Definition~\ref{def:2.1}. 
In particular, we obtain for $f\in L^p(\T)$ and $\ep>0$ that 
\begin{equation*}
\om_p^s (f,\ep) \le \om_p (f^*,\ep) 
= \om_p (f,\ep) = (\om(|f|^p,\ep) )^{1/p} \le 2\om_p^s(|f|,\ep)\ .
\end{equation*}

For any $f$ affiliated with $\N$, we let $t\to\mu(f,t)$ denote the 
decreasing rearrangement of $|f|$ on $[0,1]$; $\mu (f,t) = \inf \{r\ge0: 
\T\circ E_{|f|} ((r,\infty)) \le t\}$. 

\begin{lem}\label{lem:2.3} 
Let $1\le p<\infty$, $f,g\in L^p(\T)$, and $\ep>0$.
\begin{align}
\om_p (f+g,\ep) & \le \om_p (f,\ep) + \om_p (g,\ep) \label{eq:2.6}\\
\noalign{\hbox{and}}
\om_p^s (f+g,\ep) & \le \om_p^s (f,\ep) +\om_p^s (g,\ep)\ .\notag
\end{align}
If $f$ is self-adjoint, then 
\begin{equation}\label{eq:2.7i}
\begin{split}
\om_p (f,\ep) = \om_p^s(f,\ep) 
& = (\bom (|f|^p,\ep ))^{1/p}\\
& = \max \{\|fP\|_p : Pf = fP,\ P\in\P,\ \text{ and } \T(P)=\ep \}\\
& = \biggl( \int_0^\ep \mu^p (f,t)\,dt\biggr)^{1/p} 
\end{split}
\end{equation}
and
\begin{equation}\label{eq:2.7ii}
\om (f,\ep) \le 2\bom (f,\ep)\ \text{ when }\ p=1\ .
\end{equation}

In general,
\begin{equation}\label{eq:2.8}
\begin{split}
\om_p^s (f,\ep ) \le \om_p (f,\ep ) & = \om_p(f^*,\ep )\\
&= \om_p (|f|,\ep ) = (\bom(|f|^p,\ep ))^{1/p}
\le 2\om_p^s (f,\ep)
\end{split}
\end{equation}
and in case $p=1$,
\begin{equation}\label{eq:2.9}
\bom (f,\ep ) \le \om(f,\ep ) \le 4\bom (f,\ep)\ .
\end{equation}

Finally, let $r= \ep^{-1/p} \|f\|_p$. 
There exists a spectral projection $P$ for $|f|$ so that $fP\in \N$ with 
\begin{equation}\label{eq:2.10}
\|fP\|_\infty \le r\text{ and } \|f(I-P)\|_p \le \tilde \om_p (f,\ep) 
\le \om_p(f,\ep)\ .
\end{equation}
\end{lem}

The case $p>1$ uses the following classical submajorization inequality, 
due to H.~Weyl \cite{W}. 

\begin{sublemma}
Let $f$ and $g$  be decreasing non-negative functions on $(0,1]$ so that 
\begin{equation*}
\int_0^x f(t)\,dt \le \int_0^x g(t)\,dt \ \text{ for all }\ 0<x\le1\ .
\end{equation*}
Then also 
\begin{equation*}
\int_0^x f^p(t) \,dt \le \int_0^x g^p(t)\,dt\ \text{ for all }\ 1<p<\infty\ ,
\end{equation*}
all $0<x\le1$.
\end{sublemma}

\begin{rems}
1. This follows easily from the corresponding ``discrete'' formulation, cf. 
\cite{GK}. 
Also, the result holds in greater generality; one does not need the 
functions to be non-negative, and moreover the conclusion generalizes 
to assert that 
\begin{equation*}
\int_0^x \Phi \circ f(t)\,d 
\le \int_0^x \Phi \circ g(t)\,dt\ \text{ for all }\ 0<x\le1
\end{equation*}
all continuous convex functions $\Phi$.

2. All the assertions of Lemma~\ref{lem:2.3} hold for {\em semi-finite\/} 
von~Neumann algebras $\N$ that are {\em atomless\/} (i.e., have no 
minimal projections), endowed with a faithful normal trace $\T$. 
Several of its assertions can also be deduced from results in \cite{FK} 
and \cite{CS}. 
For example, once one {\em proves\/} the equality of the first and last 
terms in \eqref{eq:2.7i}, one may apply Lemma~4.1 of \cite{FK} to obtain 
several of the other equalities in \eqref{eq:2.7i}, for $p=1$; 
one then has that $\om (T,\ep) = \Phi_\ep (T)$ in the notation of 
\cite{FK}, and some other results in Lemma~\ref{lem:2.3} follow from 
Theorem~4.4 of \cite{FK}. 
However we prefer to give a ``self-contained'' treatment, in part because 
we take the modulus $\om (f,\ep)$ as the primary concept in our 
development. 
\end{rems}

\begin{proof}[Proof of Lemma~\ref{lem:2.3}] 
Let $p,f,g$ and $\ep$ be as in the statement. 
\eqref{eq:2.6} is a trivial consequence of the fact that $\|\cdot\|_p$ 
is a norm (i.e., the triangle inequality). 
Also, we easily obtain that 
\begin{align}
&\om_p^s (f,\ep) \le \om_p (f,\ep) = \om_p (|f|,\ep)\label{eq:2.11i}\\
&\tilde\om_p (f,\ep) \le \om_p (f,\ep) \label{eq:2.11ii}
\end{align}
and in case $p=1$, 
\begin{equation}\label{eq:2.11iii} 
\bom (f,\ep ) \le \om (f,\ep )\ .
\end{equation}

Indeed, if $P\in\P$, then 
\begin{equation}\label{eq:2.12} 
|fP| = (Pf^* fP)^{1/2} = (P|f|^2P)^{1/2} = \big|\,|f|P\,\big|
\end{equation}
which immediately yields the equality in \eqref{eq:2.11i}. 
Since compression reduces the $L^p(\T)$ norm, we have 
\begin{equation}\label{eq:2.13} 
\|PfP\|_p  = \|P(fP)P\|_p \le \|fP\|_p 
\end{equation}
which gives the inequality in \eqref{eq:2.11i}. 
If $0\le r$ and $\T \circ E_{|f|} ((r,\infty))\le\ep$, then setting 
$P = E_{|f|} ((r,\infty))$, 
\begin{equation}\label{eq:2.14} 
\biggl( \int_{(r,\infty)} t^p \,d\tau \circ E_{|f|}(t)\biggr)^{1/p} 
= \big\|\, |f|P\big\|_p \le \om_p (f,t)\ ,
\end{equation}
yielding the inequality in \eqref{eq:2.11ii}. 
\eqref{eq:2.11iii} is trivial, since for any $P\in \P$, 
\begin{equation}\label{eq:2.15} 
|\T(fP)| \le \T(|fP|) = \|fP\|_1\ .
\end{equation}

For the non-trivial assertions of the Lemma, we need the following 
basic identities (cf. \cite{FK}, \cite{CS}). 
\begin{equation}\label{eq:2.16} 
\|f\|_p^p = \int_0^\infty t^p \,d\tau \circ E_{|f|} (t) 
\le \int_0^1 \mu^p (f,t)\,dt\ .
\end{equation} 
(The final inequality is also an equality, but this follows from 
the conclusion of our Lemma.) 

Now let $f$ be self-adjoint. 
Let $\N(f)$ denote the von Neumann algebra generated by $f$, and let $\M$ 
be a MASA contained in $\N$ with $\N(f)\subset\M$. 
Then by our initial remarks, $\tau|\M$ is atomless. 
Let us identify (as we may), $\M$ and $\tau|\M$ with an atomless 
probability space $(\Omega, \cS,\nu)$. 
It follows that we may choose a countably generated $\sigma$-subalgebra 
$\cS_0$ of $\cS$ so that $f$ is $\cS_0$-measurable and also $\nu|\cS_0$ is 
atomless. 
Denote the corresponding von-Neumann algebra by: $L^\infty (\nu|\cS_0)=\M_0$. 

It then follows that $(\Omega,\cS_0,\nu)$ is measure-isomorphic to 
$([0,1],\B,m)$ (where $\B$ denotes the Borel subsets of $[0,1]$ and $m$ 
denotes Lebesgue measure on $\B$), and moreover the measure-isomorphism 
may be so chosen that the ``random-variable'' $f$ is carried over to 
the decreasing function $t\to \mu(f,t)$ (cf. Lemma~4.1 of \cite{CS}). 
It now follows that 
\begin{equation}\label{eq:2.17} 
\int_0^x \mu^p (f,t)\,dt \le \om_p^p (f,x)\ .
\end{equation}
Indeed, it follows that there exists a set $S\in\cS_0$ with $\nu (S)=x$ and 
$\int_S|f|^p\,d\nu = \tau (|\chix_Sf|^p) = \int_0^x \mu^p (f,t)\,dt$ 
(where $\chix_S$ may be interpreted as the projection in $\M_0$ obtained via 
multiplication). 
Now we define a quantity $\beta$ (depending on $x$) by 
\begin{equation}\label{eq:2.18} 
\beta = \sup \{\|f\psi\|_1 : \psi \in \N,\ \|\psi\|_\infty \le 1,\ 
|\T (\psi)|\le x\}\ .
\end{equation}

We are going to prove that there exists a $G\in \P(\M_0)$ with 
$\T(G)=x$ and 
\begin{equation}\label{eq:2.19} 
\T(|fG|) = \T (|f|G)  = \beta\ .
\end{equation}
Note that the first equality in \eqref{eq:2.19} is trivial, since 
$G\leftrightarrow f$. 
But then all the equalities in \eqref{eq:2.7i} for the case $p=1$, 
follow immediately, for we have also that then  
$|f|G = G|f|G = |GfG|$ and so trivially $\T(|f|G) \le \bom (|f|,x)\le\beta$ 
and $\T(|f|G) \le \om_1^s (f,x) \le\beta$; 
of course also $\om (f,x)\le\beta$, hence by \eqref{eq:2.19}, 
$\beta =\om(f,x)$. 
Moreover by the argument for \eqref{eq:2.17} and \eqref{eq:2.19} we have 
that $\beta = \T(|f|G) = \int_0^x \mu (f,t)\,dt$. 

Before proving this basic claim, let us see why it also yields 
\eqref{eq:2.7i} for $p>1$ (via the Sublemma). 
Still keeping $x$ fixed, assume $0<x\le \ep\le1$, and suppose $P\in\P$ with 
$\T(P)\le\ep$. 
Now setting $g = |fP|$, $g$ is self-adjoint and ``supported'' on $P$, whence 
it easily follows that $\mu (g,t) =0$ for $t>\ep$. 

But now we obtain that 
\begin{equation}\label{eq:2.20} 
\int_0^x \mu (g,t)\,dt \le \int_0^x \mu (f,t)\,dt\ .
\end{equation}
Indeed, 
\begin{equation}\label{eq:2.21}
\begin{split}
\int_0^x \mu (g,t)\,dt 
& \le \om (g,x) = \om (fP,x)\\
& = \sup \{\|fPQ\|_1 : \T(Q) \le x\}\\
& = \sup \{|\T(fPQ\varphi)| : \varphi \in\N,\ \|\varphi\|_\infty \le1\}
\ \text{ (by duality)}\\
&\le \beta
\end{split}
\end{equation}
(since $PQ\in\N$, $\|PQ\|_\infty \le 1$, and $|\T(PQ)|\le \T(Q)\le x$). 

Now (temporarily) unfixing $x$, we also have that \eqref{eq:2.20} holds 
for $x> \ep$, since $\mu (g,t) = 0$ for all $t>\ep$. 
Thus the Sublemma yields that 
\begin{equation}\label{eq:2.22} 
\int_0^\ep \mu^p (g,t)\,dt \le \int_0^\ep \mu^p (f,t)\,dt\ .
\end{equation}
Hence in view of \eqref{eq:2.16}, 
\begin{equation}\label{eq:2.23} 
\|fP\|_p^p \le \int_0^\ep \mu^p (f,t)\,dt\ ,
\end{equation}
and so at last 
\begin{equation}\label{eq:2.24} 
\om_p (f,\ep) \le \biggl( \int_0^\ep \mu^p (f,t)\,dt\biggr)^{1/p}\ .
\end{equation}

Of course \eqref{eq:2.17} combined with \eqref{eq:2.24} now yields that 
\begin{equation}\label{eq:2.25} 
\om_p (f,\ep ) = \biggl( \int_0^\ep \mu^p (f,t)\,dt\biggr)^{1/p}\ ,
\end{equation}
and now all the equalities in \eqref{eq:2.7i} follow for $p>1$ as well. 

We now establish \eqref{eq:2.19}. 
Using the polar decomposition of $f$ and duality, we have that 
\begin{equation}\label{eq:2.26} 
\begin{split}
\beta & = \sup\{ |\T(f\psi\varphi)| :\psi,\, \varphi\in\N,\, \|\psi\|_\infty,\, 
\|\varphi\|_\infty \le 1\text{ and } |\T(\psi)|\le x\}\\
& = \sup \{\T(|f|\psi) :\psi\in\N,\, 0\le \psi\le1,\ \T(\psi)\le x\}\\
& = \sup \{\T(|f|\psi) :\psi\in\M,\, 0\le \psi\le1,\ \T(\psi)\le x\}\ .
\end{split}
\end{equation}
The last equality follows by a conditional expectation argument from 
classical probability theory. 

Indeed, given $0\le \psi\le1$ in $\N$ with $\T(\psi)\le x$, there exists a 
unique $\tilde\psi \in\M_0$ such that 
\begin{equation}\label{eq:2.27} 
\T(g\psi) = \T(g\tilde\psi)\ \text{ for all }\ g\in L^1 (\M_0)\ .
\end{equation}
It follows that then $0\le \tilde\psi\le1$ and $\T(\tilde\psi)\le x$; 
this yields the desired equality.

Now let $K$ be defined: 
\begin{equation}\label{eq:2.28} 
K= \{\psi \in \M_0 :0\le \psi \le1 \text{ and }  \T(\psi) \le x\}\ .
\end{equation}
Then $K$ is a weak* compact convex set, thus 
\begin{equation}\label{eq:2.29} 
K = \om^* - \oco \{\varphi :\varphi\in \Ext K\}
\end{equation}
and moreover 
\begin{equation}\label{eq:2.30} 
\beta = \sup \{\T (|f|\varphi) :\varphi \in \Ext K\}\ .
\end{equation}

Now we claim that if $\varphi \in \Ext K$, $\varphi$ is a {\it projection\/}. 
To see this, again identifying $\M_0$ with $L^\infty (\Omega,\cS_0,\nu|\cS_0)$, 
we regard $\varphi$ as an $\cS_0$-measurable function on $\Omega$. 
Were $\varphi$ not a projection, we could choose $0<\delta<\frac12$ so 
that setting $F= \{\om\in\Omega :\delta\le \varphi (\om) \le 1-\delta\}$, 
then $\mu(F)>0$. 
Since $\mu$ is atomless, choose a measurable $E\subset F$ with $\mu(F)=
\frac12\mu (E)$. 
Now define $g$ by 
\begin{equation}\label{eq:2.31} 
g = \frac{\delta} 2 \chix_E - \frac{\delta}2 \chix_{F\sim E}\ .
\end{equation}
Then $g \ne 0$, $\T(g)=0$, and $0\le \varphi \pm g\le 1$. 
But then $\T(\varphi\pm g)\le\ep$, hence $\varphi\pm g\in K$ and 
$\varphi = \frac{(\varphi+g)+ (\varphi-g)}2$, contradicting 
the fact that $\varphi \in \Ext K$. 
(For a proof of this claim in a more general setting, see \cite{CKS}.) 

We finally observe that the supremum in \eqref{eq:2.26} is actually 
attained, thanks to the $\om^*$-compactness of $K$. 
But it then follows that this is attained at an extreme point of $K$, 
i.e., there indeed exists a $G\in \P(\M_0)$ with $\T(G)=x$, 
satisfying \eqref{eq:2.19}. 

We may now also easily obtain \eqref{eq:2.7ii}. 
Letting $f= f^+ - f^-$ where $f^+ \cdot f^- = 0$ and $f^+ ,f^- \ge0$, 
we have (by the proof of \eqref{eq:2.7i}) 
\begin{equation}\label{eq:2.32} 
\begin{split} 
\om (f,\ep) & = \sup \{\T(|f|P) :P\in \P (\M_0),\ \T(P)\le \ep\}\\
& = \sup \{\T(f^+ P) + \T(f^- P) :P\in \P (\M_0),\ \T(P)\le \ep\}\\
& \le 2\sup \{|\T(fP)| :P\in \P(\M_0),\ \T(P)\le\ep\}\\
& \le 2\bom (f,\ep)
\end{split}
\end{equation}

The first equality in \eqref{eq:2.8} follows from the fact that for a 
general $f$ affiliated with $\N$, there exists a unitary $U$ in $\N$ with 
$f= U|f|$ (thanks to the finiteness of $\N$). 
But then $|f|$ and $|f^*|$ are unitarily equivalent, which yields that 
$\mu (f,t) = \mu (f^*,t)$ for all $t$, and hence the desired equality 
follows by the final equality in \eqref{eq:2.7i}. 

It remains to prove the last inequalities in \eqref{eq:2.8} and 
\eqref{eq:2.9}, and the final statement of the lemma. 
Let $f= g+ih$ with $g$ and $h$ self-adjoint (and so in $L^p(\T)$). 
Then 
\begin{equation}\label{eq:2.33}
\begin{split}
\om_p (f,\ep) & \le \om_p (g,\ep) + \om_p (h,\ep)\ \text{ by \eqref{eq:2.6}}\\
& = \om_p^s (g,\ep) + \om_p^s (h,\ep)\ \text{ by \eqref{eq:2.7i}}\ .
\end{split}
\end{equation}

But if $\varphi=g$ or $h$, then 
\begin{equation}\label{eq:2.34}
\om_p^s (\varphi,\ep) \le\om_p^s (f,\ep)\ .
\end{equation}
Indeed, if $P\in \P$, $\T(P)\le\ep$, then $PfP = PgP+iPhP$. 
But $PgP$ and $PhP$ are both self adjoint, hence $\|P\varphi P\|_p 
\le \|PfP\|_p$, yielding \eqref{eq:2.34}. 
Of course \eqref{eq:2.33} and \eqref{eq:2.34} yield the final inequality 
in \eqref{eq:2.8}. 
Similarly, in case $p=1$, 
\begin{equation}\label{eq:2.35} 
\begin{split}
\om (f,\ep) & \le \om (g,\ep) + \om (h,\ep)\ \text{ by \eqref{eq:2.6}}\\
&\le 2\bom (g,\ep) + 2\bom (h,\ep)\ \text{ by \eqref{eq:2.7ii}}\\
&\le 4\bom (f,\ep)
\end{split}
\end{equation}
since we also have for $\varphi = g$ or $h$, that 
$\bom (\varphi,\ep) \le\bom (f,\ep)$ (by an argument similar to that 
for \eqref{eq:2.34}). 

To obtain the final assertion of the lemma, let $r= \mu (f,\ep)$, 
and let $E= E_{|f|}$. 
Now if $\bar{\ep} = \T(E[r,\infty))$ 
then since 
\begin{equation}\label{eq:2.36} 
E([r,\infty)) = \bigwedge \{E([s,\infty)): s<r\}\ ,
\end{equation}
we have $\ep \le \bar\ep$.
Thus 
\begin{equation}\label{eq:2.37}
r^p \ep \le r^p \bar \ep \le \int_{[r,\infty)} t^p \,d\tau \circ E (t)
\le \int_{[0,\infty)} t^p\,d\tau \circ E (t) 
= \|f\|_p^p\ .
\end{equation}
Hence 
\begin{equation}\label{eq:2.38}
r\le \ep^{-1/p} \|f\|_p \ .
\end{equation}
Now also by the definition of $r$, $\T (E(r,\infty))\le\ep$, and so 
\begin{equation}\label{eq:2.39} 
\T (|f|^p E_{(r,\infty)}) =  \int_{(r,\infty)} t^p d \T\circ E(t) 
\le \tilde\om_p(f,\ep)^p\ .
\end{equation}
Finally, let $f= U|f|$ be the polar decomposition of $f$.
In particular, $U$ is a partial isometry belonging to $\N$.
Then $P= E([0,r])$ satisfies \eqref{eq:2.10}.
Indeed, $fP = U|f|P$ and $\|\, |f|P\|_\infty \le r$, so also
$\|U|f|P\|_\infty \le r$, and 
\begin{equation*}
\begin{split}
\|U|f|(I-P)\|_p &\le \|\, |f|(I-P)\|_p = (\T(|f|^p E_{(r,\infty)})^{1/p}\\
& \le \tilde\om_p (f,\ep)\ \text{ by \eqref{eq:2.39}. }\quad \qed
\end{split}
\end{equation*}
\renewcommand{\qed}{}
\end{proof}

\begin{rems}
1. We have given a self-contained proof of the basic inequality 
\eqref{eq:2.24} for the sake of completeness. 
An alternate deduction may be obtained as follows. 
The remarks preceding \eqref{eq:2.17} actually yield that for any 
$g\in L^p(\T)$, $\|g\|_p = \|\mu (g,\cdot)\|_p$. 
Let $f$ be as in the proof of \eqref{eq:2.24} and fix a $P\in \P$ with 
$\T (P)=\ep$. 
We apply this observation to $g=fP$. 
First, Proposition~1.1 of \cite{CS} yields that for any $0<x\le1$, 
\begin{equation*}
\int_0^x \mu (fP,t)\,dt \le \int_0^x \mu (f,t) \mu (P,t)\,dt\ .
\end{equation*}
Hence applying the Sublemma and the observation, 
\begin{equation*}
\begin{split}
\|fP\|_p^p = \int_0^1 \mu (fP,t)^p\,dt 
& \le \int_0^1 (\mu (f,t)\mu(P,t))^p\,dt\\
& = \int_0^\ep \mu^p (f,t)\,dt
\end{split}
\end{equation*}
which of course yields \eqref{eq:2.23} and hence \eqref{eq:2.24}. 

2. Rather than making use of the measure isomorphism of $(\Omega,\cS_0,\nu|
\cS_0)$ with $([0,1],\B,m)$, one can use the following more elementary
procedure, in demonstrating \eqref{eq:2.17}.
Let $r=\mu (f,x)$.
Then it follows that setting $P= E_{|f|} ((r,\infty))$, $\T(P)\le x$
and $\T(E_{|f|} ([r,\infty))) \ge x$.
Using that $\T|\M$ is atomless, choose $Q\in \P(\M)$ with $Q\le E_{|f|}
(\{r\})$ so that $\T(Q) +\T(P)=x$.
Then
\begin{equation*}
\begin{split}
\T(|f(P+Q)|^p)
& = \T(|f|^p (P+Q))\\
& = r\T(Q) + \int_{(r,\infty)} t^p\,d\tau \circ E_{|f|}(t)\\
& = \int_0^x \mu^p (f,t) \,dt\ .
\end{split}
\end{equation*}
Here, the first two equalities are trivial; however the third one follows
by a direct elementary (but somewhat involved) argument.
(We are indebted to Ken Davidson for this Remark.)
\end{rems}

We next use the modulus of uniform integrability to establish a
criterion for relative weak compactness.

\begin{DEF}\label{def:2.3}
A subset $W$ of $L^1(\T)$ is called uniformly integrable if
\begin{equation*}
\lim_{\ep\to0} \sup_{f\in W} \om (f,\ep ) =0\ .
\end{equation*}
\end{DEF}

\begin{Com}
The assumption that $\T$ is atomless implies uniformly integrable subsets
are bounded in $L^1(\T)$.
In fact, it then follows that if $W$ satisfies that
$\sup_{f\in W} \om(f,\ep_0) <\infty$ for some $\ep_0>0$, $W$ is bounded.
\end{Com}

\begin{prop}\label{prop:2.2}
Let $(f_n)$ be a given sequence in $L^1(\T)$.
The following are equivalent
\begin{itemize}
\item[(i)] $(f_n)$ is relatively weakly compact in $L^1(\T)$.
\item[(ii)] $(f_n)$ is uniformly integrable.
\item[(iii)] $(|f_n|)$ is relatively weakly compact.
\item[(iv)] $(f_n)$ is bounded in $L^1(\T)$ and $\lim_{\ep\to0} \sup_n
\tilde\om_1(f_n ,\ep)=0$.
\item[(v)] For all $\ep>0$, there exists an $r<\infty$ so that for all $n$,
\begin{equation*}
d_{L^1(\T)} (f_n,r \Ba(\N))<\ep\ .
\end{equation*}
\end{itemize}

Moreover if $(f_n)$ is bounded in $L^1(\T)$ and
\begin{equation}\label{eq:2.41}
\eta = \lim_{\ep\to0} \sup_n \om(f_n,\ep)>0\ ,
\end{equation}
there exists a sequence $P_1,P_2,\ldots$ of pairwise
orthogonal projections in $\P$ and
$n_1<n_2<\cdots $ so that
\begin{equation}\label{eq:2.42}
|\T(f_{n_k} P_k)| > \frac{\eta}5 \ \text{ for all } \ k\ .
\end{equation}
\end{prop}

\begin{rem}
$\Ba (\N)$ denotes the closed unit ball of $\N$; thus $r\cdot\Ba(\N) =
\{f\in\N :\|f\|_\infty \le r\}$.
For $W\subset L^1(\T)$ and $f\in L^1(\T)$, $d_{L^1(\T)}(f,W) =
\inf \{\|f-w\|_1 :w\in W\}$ by definition.
Our proof of (iv) $\implies$ (v) reduces, via the proof of Lemma~\ref{lem:2.3},
to a standard truncation argument in the case of commutative $\N$.
\end{rem}

\begin{proof}
Once (i) $\Leftrightarrow$ (ii) is established, the other equivalences in
this Proposition follow easily from 2.3.
Indeed, we have by the equalities in \eqref{eq:2.8} that
\begin{equation*}
\lim_{\ep\to0} \sup_n \om (f_n,t) = \lim_{\ep\to0} \sup_n \om (|f_n|,\ep)\ ,
\end{equation*}
whence we have the equivalence of (i)--(iii).
Now trivially (ii) $\implies$ (iv) since $\tilde\om_1(f,\ep)\le\om(f,\ep)$
for any $f\in L^1(\T)$ and $\ep>0$ (see \eqref{eq:2.10}).
Suppose first that $(f_n)$ satisfies (v).
Then given $\ep>0$,  for each $n$ we may choose $\psi_n \in\N$,
$\|\psi_n\|_\infty \le r$, with
\begin{equation}\label{eq:2.43}
\|f_n-\psi_n\|_{L^1(\T)} <\ep\ .
\end{equation}
But then for any $\delta <\ep$,
\begin{equation}\label{eq:2.44}
\om(f_n,\delta) \le \om (f_n-\psi_n,\delta)
+ \om(\psi_n,\delta) <\ep + r\delta\ .
\end{equation}
Hence $\varlimsup_{\delta\to0} \sup_n \om (f_n,\delta) \le \ep$,  proving (ii).
On the other hand, suppose (iv) holds.
Let $\ep>0$, and choose $\delta >0$ so that
\begin{equation}\label{eq:2.45}
\tilde\om_1(f_n,\delta) <\ep\ \text{ for all } \ n\ .
\end{equation}
Also, let $M=\sup \|f_n\|_{L^1(\T)}$.
Then setting $r= \delta^{-1}M$, it follows by the
final statement of Lemma~\ref{lem:2.3} that for each $n$, we may choose
$\psi_n \in r\Ba \N$ with
\begin{equation*}
\|\psi_n -f_n\|_{L^1(\T)} \le \tilde\om_1(f,\delta) <\ep\ ,
\end{equation*}
proving (iv) $\implies$ (v).

To prove the equivalences of (i) and (ii), we use the following
classical criterion due to C.~Akemann \cite{A}:
{\em A bounded set $W$ in the predual of a von-Neumann algebra $\M$ is
relatively compact if and only if for any sequence $P_1,P_2,\ldots$ of
disjoint projections in $\M$,}
\begin{equation}\label{eq:2.46}
\lim_{j\to\infty} \sup_{w\in W} |P_j(w)| =0\ .
\end{equation}

Now suppose first that $(f_n)$ is {\em not\/} relatively weakly compact;
then choosing disjoint $P_j$'s as in the above criteria, we obtain that
\begin{equation}\label{eq:2.47}
\mathop{\varlimsup}\limits_{j\to\infty} \sup_n |\T(P_jf_n)| = \delta >0\ .
\end{equation}
But $\lim \T(P_j) =0$, since the $P_j$'s  are disjoint.
It follows immediately that
\begin{equation}\label{eq:2.48}
\lim_{\ep\to0}\sup_n \bom (f_n,\ep) \ge\delta\ ,
\end{equation}
which together with \eqref{eq:2.9}, proves that (ii) $\implies$ (i).

Finally, to show that (i) $\implies$ (ii), assume instead that $\eta >0$,
where $\eta $ is given in \eqref{eq:2.41}.
It now suffices to demonstrate the final  assertion of 2.5, for then $(f_n)$
is not relatively weakly compact by Akemann's criterion.
Let $0<\ep<\eta $ with $\frac{\eta }4 -\ep >\frac{\eta}5$.
By \eqref{eq:2.41}, choose $n_1$ with
\begin{equation}\label{eq:2.49}
\om \left(f_{n_1},\frac12\right) > \eta -\ep\ .
\end{equation}
Then choose (by \eqref{eq:2.9} of Lemma~\ref{lem:2.3}), $Q_1\in \P$ with
$\T(Q_1) \le 1/2$ and
\begin{equation}\label{eq:2.50}
|\T(f_{n_1}Q_1)| > \frac{\eta-\ep}4\ .
\end{equation}
Since $f_{n_1}$ is integrable, $\{f_{n_1}\}$ is uniformly integrable, so
we may choose $0<\ep_2<1$ so that
\begin{equation}\label{eq:2.51}
\om (f_{n_1},\ep_2) < \frac{\ep}{2}\ .
\end{equation}
Next by \eqref{eq:2.41}, choose $n_2>n_1$ with
\begin{equation}\label{eq:2.52}
\om (f_{n_2},\ep_2) > \eta -\ep\ .
\end{equation}
(It is easily seen, thanks to the uniform integrability of finite
sets in $L^1(\T)$, that in fact $\eta =\break \lim_{\ep\to0}
\varlimsup_{n\to\infty} \om (f_n,\ep)$; thus we may insure that $n_2$ may be
chosen larger than $n_1$.)
Again using  \eqref{eq:2.52}  and \eqref{eq:2.9},
choose $Q_2\in\P$ with $\T(Q_2)
\le \frac{\ep_2}{2^2}$ and
\begin{equation}\label{eq:2.53}
|\T(f_{n_2} Q_2)| > \frac {\eta -\ep}4\ .
\end{equation}
Then choose $\ep_3 < \ep_2$ so that
\begin{equation}\label{eq:2.54}
 \om (f_{n_2}, \ep_3) < \frac{\ep}{2} \ .
 \end{equation}

Continuing by induction, we obtain $n_1<n_2<\cdots$,
$1 = \ep_1>\ep_2 > \cdots$, and projections $Q_1,Q_2,\ldots$ in
$\P$ so that for all $k$,
\begin{gather}
\T(Q_k)\le  \frac{\ep_k}{2^k}
\label{eq:2.55i}\\
\om (f_{n_k},\ep_{k+1}) < \frac{\ep}{2}\label{eq:2.55ii}
\end{gather}
and
\begin{equation}\label{eq:2.55iii}
|\T (f_{n_k} Q_k) | > \frac{\eta-\ep}4\ .
\end{equation}
Now set $P_k = Q_k \wedge (\wedge_{j>k} (1-Q_j))$, for $k=1,2,\ldots$.
Evidently the $P_k$'s are pairwise orthogonal. 
For each $i$, let $\tilde Q_i = Q_i-P_i$. 
Now by subadditivity of $\T$, 
\begin{equation*}
\begin{split}
\T(P_i) & \ge \T(Q_i) - \biggl( 1-\T \bigwedge_{j>i} (1-Q_j)\biggr)\\
& \ge \T(Q_i) - \sum_{j>i} \T(Q_j)\ .
\end{split}
\end{equation*}
But 
\begin{equation*}
\begin{split}
\sum_{j>i} \T(Q_j) \le \sum_{j>i} \frac{\ep_j}{2^j}
& < \ep_{i+1} \sum_{j>i} \frac1{2^j} \text{ by (2.57)}\\
& < \ep_{i+1}\ .
\end{split}
\end{equation*}
Hence we have 
\begin{equation}\label{eq:2.60} 
\T(\tilde Q_i) \le \sum_{j>i} \T(Q_j) <\ep_{i+1}\ .
\end{equation}
Thus by (2.58), 
\begin{equation}\label{eq:2.61} 
\|f_{n_i} \tilde Q_i\|_1 \le  \om (f_{n_i},\ep_{i+1}) <\frac{\ep}2\ .
\end{equation}
Hence 
\begin{equation*}
\begin{split}
|\T(f_{n_i} P_i)| & = |\T(f_{n_i} Q_i - f_{n_i} \tilde Q_i)|\\
& \ge \frac{\eta-\ep}4 - \frac{\ep}2\ \text{ by \eqref{eq:2.61}}\\
& \ge \frac{\eta}5\ .
\end{split}
\end{equation*}
\end{proof}

\begin{rem}
The proof of the implication (i) $\implies$ (ii) itself, may quickly 
be achieved, using instead Theorem~3.5 of \cite{DSS}. 
\end{rem}

The following result is an immediate consequence of \ref{prop:2.2}.

\begin{cor}\label{cor:2.3}
A subset of $L^1(\T)$ is relatively weakly compact if and only if it
is uniformly integrable.
\end{cor}

\begin{proof}
Let $W$ be the subset, and suppose first $W$ is relatively weakly
compact, yet $\lim_{\ep\to0} \sup_{f\in W} \om (f,\ep)\break \defeq \eta >0$.
Then for each $n$, choose $f_n\in W$ with $\om (f_n,\frac1{2^n})
> \eta -\frac1{2^n}$.
It follows immediately that also $\lim_{\ep\to0} \sup_n \om (f_n,\ep)=\eta$,
hence $(f_n)$ is not relatively weakly compact by Proposition~\ref{prop:2.2}.
If $W$ is uniformly integrable, then $W$ is bounded, and then $W$ is
relatively weakly compact by Akemann's criterion, (stated
preceding \eqref{eq:2.46}).
\end{proof}

\begin{rem}
Suppose $\|f_i\|_1 \le 1$ for all $i$, and $(f_i)$ satisfies
\eqref{eq:2.41}.
Letting the $n_1<n_2 <\cdots$ be as in the proof
of \ref{prop:2.2}, we show in Section~3, using
arguments in \cite{R1}, that there exists a subsequence $(f'_i)$ of
$(f_{n_i})$ so that $(f'_i)$ is $\frac5{\eta}$-equivalent to the
usual $\ell^1$-basis, with also $[f'_i]$ $\frac5{\eta}$-complemented
in $L^1(\T)$.
Hence $(f_i)$ has a subsequence equivalent to the $\ell^1$-basis, so
of course $(f_i)$ is not relatively weakly compact.
\end{rem}

We note finally a consequence of the proof of \ref{prop:2.2}, valid
for all $1\le p<\infty$ and arbitrary (not necessarily atomic) finite
von~Neumann algebras.

\begin{cor}\label{cor:2.7}
Let $1\le p<\infty$, let $\M$ be a finite von~Neumann algebra endowed with
a faithful normal tracial state $\T$, and let $W$ be a bounded subset of
$L^p(\T)$.
Then the following are equivalent.
\begin{itemize}
\item[(i)] $\{|w|^p :w\in W\}$ is uniformly integrable.
\item[(ii)] $\lim_{\ep\to0} \sup_{f\in W} \tilde\om_p (f,\ep)=0$.
\item[(iii)] $\lim_{r\to\infty} g_W (r)=0$,
\end{itemize}
where the function $g_W$ is defined by
\begin{equation}\label{eq:2.62}
g_W (r) = \sup_{w\in W} d_{L^p(\T)} (w,r\Ba (\M))\ \text{ for }\ r>0\ .
\end{equation}
\end{cor}

\begin{proof}
(i) $\implies$ (ii) follows immediately from the (obvious) 
inequality  $\tilde\om_p(f,\ep)\le \om_p (f,\ep)$ (stated as part of 
\eqref{eq:2.10} in Lemma~\ref{lem:2.3}).

(ii) $\implies$ (iii). Assume that $\|w\|_p \le M$ for all $w\in W$.
For $r$ sufficiently large, define $\ep (r) = \ep >0$ by
\begin{equation}\label{eq:2.63}
r= \ep^{-1/p} M\ .
\end{equation}
Let $f\in W$.
Since $\ep^{-1/p} \|f\|_p \le r$, by the final assertion  of
Lemma~\ref{lem:2.3}, we may choose $P$ a spectral projection for $|f|$
so that
\begin{equation} \label{eq:2.64}
fP\in r\Ba (\M)\ \text{ and }\ \|f(I-P)\|_p \le \tilde\om_p (f,\ep)\ .
\end{equation}
It follows immediately that
\begin{equation}\label{eq:2.65}
g_W(r) \le \sup_{f\in W} \tilde \om_p (f,\ep)\ .
\end{equation}
Thus (iii) holds by (ii), since $\ep(r) \to 0$ as $r\to\infty$.
(Note also that the final assertion of \ref{lem:2.3} does not involve
the ``atomless'' hypothesis, since $\tilde\om_p(f,\ep)$ is defined
in terms of the spectral measure for $|f|$.)

(iii) $\implies$ (i).
Given $f\in W$ and $\ep>0$, choose $\psi \in r\cdot\Ba (\M)$ with
\begin{equation}\label{eq:2.66}
\|f-\psi\|_{L^p(\T)} <\ep\ .
\end{equation}
Then for any $\delta <\ep$,
\begin{equation}\label{eq:2.67}
\om_p (f,\delta) \le \om_p (f-\psi,\delta) + \om_p (\psi,\delta)
< \ep + r\delta\ .
\end{equation}
Hence $\varlimsup_{\delta\to0} \sup_{f\in W} \om_p (f,\delta)\le\ep$,
proving that (i) holds, since $\ep>0$ is arbitrary and $\om_p (f,t)
= (\om (|f|^p,t))^{1/p}$ for any $f$ and $t$, by \eqref{eq:2.8}
of Lemma~\ref{lem:2.3}.
\end{proof}

\section{Proof of the Main Theorem}
\setcounter{equation}{0}

We first assemble some preliminary lemmas, perhaps useful in a wider context.
$\N$ and $\T$ are assumed to be as in Section~2.
Let $r_1,r_2,\ldots$ denote the Rademacher functions on $[0,1]$; 
equivalently, an independent sequence of $\{1,-1\}$-valued random 
variables $(r_j)$ with $P(r_j=1) = P(r_j = -1) = \frac12$ for all $j$.

\begin{lem}\label{lem:3.1}
Let $1\le p<2$ and $(f_n)$ be a bounded unconditional basic sequence
in $L^p(\T)$, so that $(|f_i|^p)_{i=1}^\infty$ is uniformly integrable.
Then $\lim_{n\to\infty} n^{-1/p} \|f_1+\cdots + f_n\|_{L^p(\T)}=0$.
\end{lem}

\begin{rem}
Recall from the introduction that a sequence $(x_n)$ in 
a Banach space is called
{\it unconditional\/} if there is a constant $u$ so that
\begin{equation}\label{eq:3.1}
\begin{split}
&\biggl\{\Big\|\sum_{i=1}^n \alpha_i c_i x_i\Big\|
\le u\Big\|\sum_{i=1}^n c_i x_i\Big\|\biggr\}\text{ for all $n$ and scalars}\\
&c_1,\ldots, c_n\text{ and } \alpha_1,\ldots,\alpha_n
\text{ with } |\alpha_i| =1\text{ for all } i\ .
\end{split}
\end{equation}
$(x_n)$ is called $u$-unconditional if \eqref{eq:3.1} holds.
\end{rem}

\begin{proof}[Proof of \ref{lem:3.1}]
Suppose $(f_n)$ is $u$-unconditional.
Then  $(f_n)$ is $u$-equivalent to $(f_n\otimes r_n)$ in
$L^p (\N \bar\otimes L^\infty)$,  so it
suffices to prove the same conclusion for $(f_n\otimes r_n)$ instead.
Let $\beta = \T\otimes m$, where $m$ is Lebesgue measure on $[0,1)$.
We may also assume without loss of generality that $\|f_n\|_{L^p(\T)} \le 1$
for all $n$.
Now let $\ep>0$, and choose $\delta >0$ so that
\begin{equation}\label{eq:3.2}
\om (|f_n|^p,\delta) \le \ep\text{ for all } n
\end{equation}
(using that $(|f_n|^p)$ is uniformly integrable).
By the final statement of Lemma~\ref{lem:2.3}, we may by \eqref{eq:3.2}
choose for each $j$ a $P_j \in \P = \P(\N)$ so that
$f_j P_j \in\N$ with
\begin{equation}\label{eq:3.3}
\|f_jP_j\|_\infty \le \frac1{\delta}\ \text{ and }\
\|f_j  (I-P_j)\|_p^p \le \ep\ .
\end{equation}
Then fixing $n$,
\begin{equation}\label{eq:3.4}
\Big\| \sum_{i=1}^n f_i\otimes r_i\Big\|_{L^p(\beta)}
\le \Big\| \sum_{i=1}^n f_i P_i \otimes r_i\Big\|_{L^p(\beta)}
+ \Big\| \sum_{i=1}^n f_i (I-P_i)\otimes r_i\Big\|_{L^p(\beta)}\ .
\end{equation}
But
\begin{equation}\label{eq:3.5}
\Big\|\sum_{i=1}^n f_iP_i\otimes r_i\Big\|_{L^p(\beta)}
\le \Big\| \sum_{i=1}^n f_i P_i\otimes r_i\Big\|_{L^2(\beta)}
\le \frac{\sqrt n}{\delta}
\end{equation}
since $\|f_iP_i\|_\infty \le \frac1{\delta}$ for all $i$.

On the other hand, since $L^p(\M)$ is type $p$ with type $p$ constant 1 for
any von-Neumann algebra $\M$,
\begin{equation}\label{eq:3.6}
\begin{split}
\Big\| \sum_{i=1}^n f_i (I-P_i)\otimes r_i\Big\|_{L^p(\beta)}
& \le \biggl( \sum_{i=1}^n \|f_i(I-P_i)\|_{L^p(\T)}^p \biggr)^{1/p}\\
&\le \ep n^{1/p}\ \text{ by \eqref{eq:3.3}.}
\end{split}
\end{equation}
(This fact follows by Clarkson's inequalities --- see the discussion in
the proof of the next lemma.)
We thus have that
\begin{equation}\label{eq:3.7}
\mathop{\varlimsup}\limits_{n\to\infty}
n^{-1/p} \Big\| \sum_{i=1}^n f_i \otimes r_i
\Big\|_{L^p(\beta)}
\le \lim_{n\to\infty} \frac{n^{1/2}}{\delta n^{1/p}} + \ep =\ep
\end{equation}
by \eqref{eq:3.5} and \eqref{eq:3.6}.
Since $\ep>0$ is arbitrary, the
conclusion of the lemma follows. 
\end{proof}

\begin{rems}
1. It follows easily from the above proof that in fact if $(f_n)$ satisfies
the hypothesis of \ref{lem:3.1}, then
$\lim_{n\to\infty} n^{-1/p} \|f'_1 +\cdots + f'_n\|_p =0$
uniformly over all subsequences $(f'_n)$ of $f_n$.

2. The proof of Lemma~\ref{lem:3.1} yields the following quantitative result.
{\it Fix $\ep>0$, and let $(f_j)$ be a bounded sequence in $L^p(\T)$ so that
there exists an $r<\infty$ with $d_{L^p(\T)} (f_j,r\Ba\N) <\ep$ for all $j$.
Then $\varlimsup_{n\to\infty} \EEo n^{-1/p} \|\sum_{j=1}^n  r_j
(w)f_j\|_{L^p(\T)} \le \ep$.}
Indeed, for each $j$, choose $\varphi_j\in r\Ba \N$ with $\|f_j-\varphi_j\|
_{L^p(\T)} \le\ep$.
Then fixing $n$, \eqref{eq:3.4}--\eqref{eq:3.6} yield
\begin{equation*}
\begin{split}
\Big\|\sum_{i=1}^n f_i \otimes r_i\Big\|_{L^p(\beta)}
& \le \Big\|\sum_{i=1}^n \varphi_i\otimes r_i\Big\|_{L^p(\beta)}
+ \Big\|\sum_{i=1}^n (f_i -\varphi_i)\otimes r_i\Big\|_{L^p(\beta)}\\
& \le r\sqrt n + \ep n^{1/p}\ .
\end{split}
\end{equation*}
Hence $\varlimsup_{n\to\infty} n^{-1/p} \|\sum_{i=1}^n f_i
\otimes r_i\|_{L^p(\beta)} \le \ep$ as desired.\qed
\end{rems}

We next give a criterion for a finite or infinite sequence in
$L^p(\T)$ to be equivalent to the usual $\ell^p$ basis.

\begin{lem}\label{lem:3.2}
Let $u\ge 1$, $\delta >0$, $1\le p<2$, and $f_1,\ldots, f_n$ elements of
$\Ba (L^p(\N))$ be given so that $(f_i)_{i=1}^n$ is $u$-unconditional.
Assume there exist pairwise 
orthogonal projections $P_1,\ldots,P_n$ in $\P$ so that
\begin{equation}\label{eq:3.8}
\T(|P_j f_j P_j|^p) \ge \delta^p\ \text{ for all }\ 1\le j\le n\ .
\end{equation}
Then $(f_i)_{i=1}^n$ is $C$-equivalent to the usual $\ell_n^p$ basis,
where $C= u\sqrt3\, \delta^{-1}$.
\end{lem}

\begin{proof}
We first note that (using interpolation),
$L^p(\T)$ satisfies Clarkson's inequalities:
for all $x,y\in L^p(\T)$,
\begin{equation}\label{eq:3.9}
\|x+y\|_p^p + \|x-y\|_p^p
\le 2(\|x\|_p^p + \|y\|_p^p)\ .
\end{equation}
It follows immediately by induction on $n$ that $L^p(\T)$ is type $p$
with constant one; that is, for any $x_1,\ldots,x_n$ in $L^p(\T)$,
\begin{equation}\label{eq:3.10}
\begin{split}
\sum_{A\vee\pm} \|\pm x_1 \pm\cdots \pm x_n\|_p^p
& = \int_0^1 \Big\|\sum_{i=1}^n r_i(\om)x_i\Big\|_p^p\,d\om\\
&\le \biggl( \sum_{i=1}^n \|x_i\|_p^p \biggr)\ .
\end{split}
\end{equation}
Now let scalars $a_1,\ldots,a_n$ be given, and let $f= \sum_{i=1}^n a_if_i$.
We obtain from \eqref{eq:3.10} that since $(f_i)$ is $u$-unconditional,
\begin{equation}\label{eq:3.11}
\|f\|_p \le u\biggl( \sum_{i=1}^n |a_i|^p\biggr)^{1/p} \ .
\end{equation}

Now fix $\om$ and set $f_\om = \sum_{i=1}^n a_i  r_i (\om ) f_i$.
Then
\begin{equation}\label{eq:3.12}
\|f_\om \|_p^p \ge \sum_{j=1}^n \|P_j f_\om P_j\|_p^p\ .
\end{equation}
Thus integrating over $\om$ and again using unconditionality,
\begin{equation}\label{eq:3.13}
\begin{split}
\|f\|_p^p & \ge \frac1{u^p} \int_0^1 \|f_\om\|_p^p\,d\om\\
&\ge \frac1{u^p} \sum_{j=1}^n \int_0^1 \|P_j f_\om P_j\|_p^p \,d\om
\ \text{ by \eqref{eq:3.12}.}
\end{split}
\end{equation}
But fixing $j$, since $L^p (\T)$ is cotype 2 with constant at most $3^{1/2}$,
\begin{equation}\label{eq:3.14}
\begin{split}
\int_0^1 \|P_j f_\om P_j \|_p^p\,d\om
&\ge \frac1{3^{p/2}} \biggl( \sum_i \|P_j a_i f_i P_j\|_p^2\biggr)^{p/2}\\
&\ge \frac1{3^{p/2}} \|P_j a_j f_j P_j\|_p^p\\
&\ge \frac1{3^{p/2}} |a_j|^p \delta^p\ \text{ by \eqref{eq:3.8}.}
\end{split}
\end{equation}
Thus in view of \eqref{eq:3.13},
\begin{equation}\label{eq:3.15}
\|f\|_p^p \ge \frac{\delta^p}{u^p 3^{p/2}}
\biggl( \sum_{j=1}^n |a_j|^p\biggr)\ ,
\end{equation}
so \eqref{eq:3.11} and \eqref{eq:3.15} now imply the conclusion of
Lemma~\ref{lem:3.2}.
\end{proof}

Our last preliminary result yields an estimate for equivalence to the
$\ell_n^p$ basis in terms of $p$-moduli.

\begin{lem}\label{lem:3.3}
Let $0<\ep<\eta/2$, $n\ge1$, and $f_1,\ldots,f_n\in\Ba L^p(\T)$ be such that
$(f_1,\ldots,f_n)$ is $u$-unconditional and there are $\delta_1\ge\delta_2
\ge \cdots\ge \delta_n >0$ so that for all $1\le j\le n$ and all $k$ with
$j<k$ (if $j<n$)
\begin{equation}\label{eq:3.16}
\om_p (f_j,\delta_j) >\eta \ \text{ and }\ \om_p (f_j,\delta_k +
\delta_{k+1} +\cdots + \delta_n) < \frac{\ep}{2}\ .
\end{equation}
Then $(f_1,\ldots,f_n)$ is $C$-equivalent to the $\ell_n^p$ basis where
$$C\le u\sqrt3 \Big( \frac{\eta}2 -\ep\Big)^{-1}\ .$$
\end{lem}

\begin{proof}
By Lemma~\ref{lem:2.3}, (see \eqref{eq:2.8}), 
we have, fixing $1\le j\le n$, that 
\begin{equation}\label{eq:3.17} 
\om_p^s (f_j,\delta_j) > \frac{\eta}2\ .
\end{equation}
Hence we may choose $Q_j \in \P$ with 
\begin{equation}\label{eq:3.18}
\|Q_jf_jQ_j\|_p > \frac{\eta}2\ \text{ and }\ 
\T(Q_j) \le \delta_j\ .
\end{equation}
Define projections $P_j$ and $\tilde Q_j$ by 
\begin{equation}\label{eq:3.19} 
P_j = Q_j \wedge  \bigwedge\limits_{k>j} (1-Q_k)\ \text{ and }\ 
\tilde Q_j = Q_j - P_j\ .
\end{equation}
Then 
\begin{equation}\label{eq:3.20}
Q_j f_j Q_j = P_j f_j P_j + \tilde Q_j f_j P_j + Q_j f_j \tilde Q_j\ .
\end{equation}
Now we have by subadditivity of $\T$ that 
$\T(\bigwedge_{k>j} (1-Q_k)) \ge 1-\sum_{k> j}\delta_k$,
and so again by subadditivity, 
\begin{equation*}
\begin{split}
\T(P_j) 
& \ge \T(Q_j) - \biggl( 1-\T\biggl( \bigwedge_{k>j} 1-Q_k\biggr) \biggr)\\
& \ge \T(Q_j) - \sum_{k>j} \delta_k\ .
\end{split}
\end{equation*}
Thus $\T(\tilde Q_j) < \sum_{k>j} \delta_k$.
Hence we have 
\begin{equation}\label{eq:3.21}
\begin{split}
\|\tilde Q_j f_j P_j\|_p \le \|\tilde Q_j f_j\|_p 
&\le \om_p \biggl( f_j^*,\sum_{k>j} \delta_k\biggr)\\
&= \om_p \biggl( f_j,\sum_{k>j} \delta_k\biggr) 
\le \frac{\ep}2\ \text{ by \eqref{eq:3.16}).}
\end{split}
\end{equation}
By the same argument,
\begin{equation}\label{eq:3.21ii}
\|Q_j f_j \tilde Q_j \|_p \le \frac{\ep}2\ .
\end{equation}
Thus from \eqref{eq:3.18}, \eqref{eq:3.20}, \eqref{eq:3.21} and 
\eqref{eq:3.21ii}, we obtain 
\begin{equation}\label{eq:3.22} 
\|P_j f_j P_j \|_p \ge \frac{\eta}2 -\ep\ .
\end{equation}
Of course $P_1,\ldots,P_n$ are {\it pairwise orthogonal\/}; 
hence Lemma~\ref{lem:3.2}
now immediately yields the conclusion of \ref{lem:3.3}. 
\end{proof}

Lemma \ref{lem:3.3} immediately yields an infinite dimensional 
conclusion as well. 
Combining this and Lemma~\ref{lem:3.1} we obtain the following definitive 
result. 

\begin{cor}\label{cor:3.4}
Let $(f_n)$ be a bounded unconditional sequence in $L^p(\T)$, $1\le p<2$. 
The following are equivalent:  
\begin{itemize}
\item[(a)] $(f_n)$ has a subsequence equivalent to the usual $\ell^p$ basis.
\item[(b)] $(|f_n|^p)$ is not uniformly integrable.
\end{itemize}
\end{cor}

\begin{proof}
(a) $\implies$ (b) follows immediately from Lemma~\ref{lem:3.1}. 
Assume that (b) holds and also assume without loss of generality that 
$\|f_n\|_p \le 1$ for all $n$. 
Then by Lemma~\ref{lem:3.1}, 
\begin{equation}\label{eq:3.23}
\eta \defeq \lim_{\ep\to0} \sup_n \om_p (f_n,\ep) >0\ .
\end{equation}
Now Lemma~\ref{lem:3.3}	
yields that there is a subsequence $(f'_n)$ of 
$(f_n)$ so that 
\begin{equation}\label{eq:3.24}
(f'_n) \text{ is $\frac{cu}{\eta}$-equivalent to the $\ell^p$ basis,}
\end{equation}
where $c$ is an absolute constant.

Indeed, fix $0<\ep<\frac{\eta}2$. 
Choose $\delta_1 \le1$ and $n_1$ so that 
\begin{equation}\label{eq:3.25} 
\om_p (f_{n_1},\delta_1) >\eta-\ep\ .
\end{equation}
Suppose $n_1<\cdots < n_j$ and $\delta_1>\delta_2 \cdots > \delta_j$ chosen
so that 
\begin{equation*}
\om_p (f_{n_i},\delta_{i+1} +\cdots + \delta_j) <\frac{\ep}2
\ \text{ for all }\ 1\le i<j\ .
\end{equation*}
By continuity of the functions $t\to\om_p (f_{n_i},t)$ for 
$i<j$ and the fact that $f_{n_j}\in L^p(\T)$, 
choose $\bar\delta_{j+1} <\delta_j$ so that 
\begin{equation}\label{eq:3.26}
\om_p (f_{n_i},\delta_{i+1} +\cdots + \delta_j +\bar\delta_{j+1}) 
<\frac{\ep}{2}\ \text{ for all }\ 
1\le i\le j\ .
\end{equation}
Then choose $\delta_{j+1} \le\bar\delta_{j+1}$ and $n_{j+1}>n_j$ so that 
\begin{equation}\label{eq:3.27}
\om_p (f_{n_{j+1}},\delta_{j+1}) >\eta-\ep\ .
\end{equation}
This completes the inductive choice of $n_1<n_2<\cdots$.

Setting $f'_k = f_{n_k}$, then $(f'_1,\ldots,f'_n)$ satisfies the 
hypotheses of Lemma~\ref{lem:3.3} for all $n$, and hence $(f'_n)$ is 
$u\sqrt3 (\frac{\eta}2-\ep)^{-1}$-equivalent to the $\ell^p$ basis by 
\ref{lem:3.3}. 
By taking $\ep$ small enough, we obtain $c\le 7$ in \eqref{eq:3.24}.
\end{proof}

\begin{rem}
The hypothesis that $(f_n)$ is unconditional may be omitted when $p=1$, 
as pointed out in the remark following the proof of 
Proposition~\ref{prop:2.2}. 
Also, it's not hard to show that the sequence $(f'_n)$ constructed above 
has its closed linear span complemented in $L^p(\T)$. 
Finally, it follows from known (rather non-trivial) results that 
if $1<p<\infty$ 
and $\N$ is  hyperfinite, then {\em every\/} semi-normalized 
{\em weakly null sequence in $L^p(\N)$ has an unconditional subsequence\/}. 
Indeed, assuming (as we may) that $\N$ acts on a separable Hilbert space, 
$L^p(\N)$ has an unconditional finite dimensional decomposition 
(see \cite{SF}, \cite{PX1}), which yields the above statement. 
Thus also in the hyperfinite case, the hypothesis that $(f_n)$ is 
unconditional may be omitted. 
We do not know, however, if this is so for general $\N$. 
\end{rem}

\begin{cor}\label{cor:3.5}
Let $(f_n)$ be a bounded unconditional sequence in $L^p(\T)$, $1\le p<2$. 
The following are equivalent.
\begin{itemize}
\item[(a)] For every subsequence $(f'_n)$ of $(f_n)$ 
\begin{equation*}
\lim_{n\to\infty} n^{-1/p} \Big\| \sum_{i=1}^n f'_i\Big\|_{L^p(\T)} =0\ .
\end{equation*}
\item[(b)] $(|f_n|^p)$ is uniformly integrable.
\end{itemize}
\end{cor}

\begin{proof} 
(a) $\implies$ (b): 
Assume (b) is false. 
Then by Corollary~3.4 there exists a subsequence $(f'_n)$ equivalent 
to the usual $\ell^p$-basis. 
In particular 
\begin{equation*}
\liminf_{n\to\infty} n^{-1/p} \Big\|\sum_{i=1}^n f'_i\Big\|_{L^p(\T)}>0\ .
\end{equation*}
which contradicts (a). 

(b) $\implies$ (a). 
This follows from Lemma~3.1, since condition (b) implies that 
$(|f'_n|)^p$ is uniformly integrable for any subsequence $(f'_n)$ of $(f_n)$.
\end{proof}

We now turn to the proof of the Main Theorem. 
First we 
give some preliminary results concerning ultrapowers of Banach spaces 
and the standard construction of the ultrapower of a finite von~Neumann 
algebra (cf. \cite{McD}, \cite{V}). 

Fix $U$ a free ultrafilter on $\nat$. 
For a given Banach space $X$, let $\ell^\infty (X)$ denote the set of 
bounded sequences in $X$, under the norm $\|(x_n)\| = \sup_n \|x_n\|$, and 
set 
\begin{equation}\label{eq:4.1} 
E_U = \{(x_n)\in\ell^\infty (X) : \lim_{n\in U} \|x_n\| =0\}\ .
\end{equation}
Then $X_U$, the ultrapower of $X$ with respect to $U$, is given by 
\begin{equation}\label{eq:4.2} 
X_U = \ell^\infty (X) /E_U\ .
\end{equation}

Now fix $\N$ a finite von Neumann algebra with a normal faithful tracial 
state $\T$, and define $I_U$ by 
\begin{equation}\label{eq:4.3} 
I_U = \{(x_n) \in \ell^\infty (\N): \lim_{n\in U} \T (x_n^* x_n) =0\}\ .
\end{equation}
Then $I_U$ is a norm-closed two-sided ideal in $\ell^\infty (X)$; we define 
$\N^U$ (a different object than $\N_U$!) by 
\begin{equation}\label{eq:4.4} 
\N^U = \ell^\infty (\N)/I_U\ .
\end{equation}
Then by the references cited above,  
$\N^U$ is a $W^*$-algebra (i.e., an abstract von~Neumann algebra) 
with a normal faithful tracial state $\T_U$ given by 
\begin{equation}\label{eq:4.5} 
\T_U (\pi (x_n)) = \lim_{n\in U} \T (x_n)
\end{equation}
where $\pi :\ell^\infty (\N) \to \N^U$ is the quotient map. 

The next result yields that $L^p(\N^U)$ may be regarded as a subspace of 
the Banach space ultrapower $L^p(\N)^U$. 

\begin{lem}\label{lem:4.1} 
Let $1\le p<\infty$ and let $Y_p$ denote the closure of $\ell^\infty (\N)$ 
in the Banach space $\ell^\infty (L^p(\N))$. 
Then $\pi$ has a unique extension to a bounded linear map $\tilde\pi :Y_p 
\to L^p (\N^U)$. 
Moreover, for $(x_n)\in Y_p$, 
\begin{equation}\label{eq:4.6} 
\|\tilde\pi ((x_n))\|_{L^p(\T_U)} 
= \lim_{n\in U} \|x_n\|_{L^p(\T)}\ .
\end{equation} 
\end{lem}

Fixing $p$ as in \ref{lem:4.1} and letting $\rho :\ell^\infty (L^p(\N))
\to L^p (\N)^U$ be the quotient map, Lemma~\ref{lem:4.1} yields there is 
a unique isometric embedding $i:L^p(\N^U) \to L^p(\N)^U$ so that 
the following diagram commutes: 
\begin{equation}\label{eq:4.7}
\begin{picture}(100,70)
\put(54,55){$L^p(\N^U)$}
\put(18,32){$\scriptstyle \tilde\pi$}
\put(70,32){$\scriptstyle i$}
\put(34,12){$\scriptstyle \rho$}
\put(20,20){\vector(3,4){25}}		
\put(66,50){\vector(0,-4){30}}		
\put(22,8){\vector(1,0){25}}		
\put(54,4){$L^p(\N)^U$}
\put(4,4){$Y_p$}
\end{picture}
\qquad\qquad \raise3ex\hbox{.}
\end{equation}

\begin{proof}
Since $\pi$ is a $*$-homomorphism of $\ell^\infty (\N)$ onto $\N^U$, we 
have for any continuous function $f:[0,\infty)\to \complex$ and any 
$x= (x_n)\in\ell^\infty (\N)$, 
\begin{equation}\label{eq:4.8} 
\pi \left( (f(x_n^* x_n))_{n=1}^\infty \right) 
= f(\pi (x^*) \pi (x))\ .
\end{equation}
Applying this to $f(t) = |t|^{p/2}$, we get by the trace formula 
\eqref{eq:4.5} that 
\begin{equation}\label{eq:4.9} 
\|\pi (x)\|_{L^p(\T_U)} = \lim_{n\in U} \|x_n\|_{L^p(\T)}\ .
\end{equation}
In particular, 
\begin{equation}\label{eq:4.10}
\begin{split}
\|\pi (x)\|_{L^p (\T_U)} 
& \le \sup_n \|x_n\|_{L^p(\T)}\\ 
& = \|x\|_{\ell^\infty (L^p (\N))}\ .
\end{split}
\end{equation}
Thus $\pi$ extends by continuity to a contraction $\tilde\pi :Y_p \to 
L^p(\N^U)$. 
Now let $x= (x_n)$ belong to $Y_p$, and let $\ep>0$. 
Then choose $y = (y_n)$ in $\ell^\infty (\N)$ so that 
\begin{equation}\label{eq:4.11} 
\|x-y\|_{\ell^\infty (L^p(\N))} <\ep\ .
\end{equation}

It follows from \eqref{eq:4.11} that 
\begin{equation}\label{eq:4.12} 
\Big|\ \|\pi (x)\|_{L^p(\T_U)} - \|\pi (y) \|_{L^p(\T_U)}\Big|  <\ep 
\end{equation}
and 
\begin{equation}\label{eq:4.13} 
\Big|\lim_{n\in U}  \|x_n\|_{L^p(\T)} 
- \lim_{n\in U} \|y_n\|_{L^p(\T)}\Big| <\ep\ .
\end{equation} 
Since \eqref{eq:4.9} holds, replacing ``$x$'' by ``$y$'' 
in its statement, we have 
from \eqref{eq:4.12} and \eqref{eq:4.13} that 
\begin{equation}\label{eq:4.14} 
\Big|\, \| \pi (x)\|_{L^p (\T_U)} - \lim_{n\in U} \|x_n\|_{L^p(\T)}\Big| 
< 2\ep\ .
\end{equation}
Since $\ep>0$ is arbitrary, \eqref{eq:4.6} holds for all $x= (x_n)$ 
in $Y_p$.
\end{proof}

\begin{lem}\label{lem:3.7new}
Let $1\le p<2$, and let $(x_{ij})$ be an infinite matrix in $L^p(\N)$ so that 
for some $C\ge1$, each row and each column of $(x_{ij})$ is $C$-equivalent 
to the usual $\ell^2$-basis. 
Then for every free ultrafilter $U$ on $\nat$ 
\begin{equation}\label{freeultrafilter}
\sup_{j\in\nat} \lim_{i\in U} d_{L^p(\T)} (x_{ij},r\Ba(\N)) \to0\text{ as } 
r\to\infty 
\end{equation}
\end{lem}

\begin{proof} 
Define for each $j\in \nat$ a function $g_j: \real^+\to\real^+$ by 
\begin{equation*}
g_j(r)  = \sup_i d_{L^p(\T)} (x_{ij} ,r\Ba(\N))\ .
\end{equation*}
For fixed $j$, $(x_{ij})_{i=1}^\infty$ is $C$-equivalent to the usual 
$\ell^2$-basis, so by Corollary~3.4 and Corollary~2.7, 
$(|x_{ij}|^p)_{i=1}^\infty$ is uniformly integrable and 
\begin{equation}\label{eq:3.44new}
\lim_{r\to\infty} g_j(r)=0\ .
\end{equation}
Now \eqref{eq:3.44new} implies that $(x_{ij})_{i=1}^\infty$ belongs 
to $Y_p$. 
Let $\tilde\pi$ be as in the statement of Lemma~3.6 and define $x_j$ by 
\begin{equation*}		
x_j = \tilde\pi \Big( (x_{ij})_{i=1}^\infty\Big) \in L^p (\N^U)\ .
\end{equation*} 
Now we claim that 
\begin{equation}\label{eq:4.19} 
(x_j)\text{ is $C$-equivalent to the $\ell^2$-basis.}
\end{equation}
Indeed, using the hypotheses of Theorem~\ref{thm:1.1} and 
Lemma~\ref{lem:4.1}, we have for any $n$ and scalars $c_1,\ldots,c_n$, that 
\begin{equation*}		
\begin{split} 
\Big\| \sum_{j=1}^n c_j x_j\Big\|_{L^p(\T_U)} 
& = \Big\| \tilde\pi \left( \biggl( \sum_{j=1}^n c_j x_{ij}\biggr)_{i=1}^\infty
\right)\Big\|_{L^p(\T_U)}\\
& = \lim_{i\in U} \Big\| \sum_{j=1}^n c_j x_{ij}\Big\|_{L^p(\T)}
\ \text{ by \eqref{eq:4.6}}\\
& \buildrel C\over \sim \Big(\sum |c_j|^2\Big)^{1/2}\ .
\end{split}
\end{equation*}

Now define $g:\real^+\to\real^+$ by 
\begin{equation*}			
g(r) = \sup_j d_{L^p(\T_U)} (x_j,r\Ba (\N^U))\ .
\end{equation*}
Again by \eqref{eq:4.19} and Corollary~\ref{cor:3.4}, 
$(|x_j|^p)_{j=1}^\infty$ is uniformly integrable in $L^p(\T_U)$, so by 
Corollary~\ref{cor:2.7}  we have that 
\begin{equation}\label{eq:4.22} 
\lim_{r\to\infty} g(r)=0\ .
\end{equation}

Now let $\ep>0$. 
Since $\pi$ is a quotient map of $\ell^\infty(\N)$ onto $\N^U$, 
it follows that fixing $j$, there exists for every $r>0$, 
$(y_{ij})_{i=1}^\infty \in r\Ba (\N)$ so that 
\begin{equation*}			
\|x_j -\pi ((y_{ij})_{i=1}^\infty)\|_{L^p(\T_U)} < g(r)+\ep \ .
\end{equation*}
Hence by Lemma~\ref{lem:4.1}, 
\begin{equation*}			
\lim_{i\in U} \|x_{ij} -y_{ij}\|_{L^p(\T)} < g(r)+\ep\ ,
\end{equation*}
which implies that 
\begin{equation*}
\lim_{i\in U} d_{L^p(\T)} (x_{ij},r\Ba (\N)) < g(r) +\ep\ .
\end{equation*}
Hence by \eqref{eq:4.22} 
\begin{equation*}
\limsup_{r\to\infty} \biggl( \sup_{j\in\nat} \lim_{i\in U} d_{L^p(\T)} 
(x_{ij},r\Ba(\N))\biggr) \le \ep\ .
\end{equation*}
Since $\ep>0$ was arbitrary, we get \eqref{freeultrafilter}.
\end{proof} 

\begin{proof}[Proof of Theorem~\ref{thm:1.1}]
Let $1\le p<2$, and let $(x_{ij})$ be as in Theorem~\ref{thm:1.1}, and let 
$U$ be a free ultrafilter on $\nat$. 
Put 
\begin{equation}\label{eq:3.47new}
h(r) = \sup_j\lim_{i\in U} d_{L^p(\T)} (x_{ij},r \Ba (\N)),\qquad r\in\real_+\ .
\end{equation}
Then $h:\real^+\to\real^+$ is a decreasing function and by 
\eqref{freeultrafilter} 
\begin{equation}\label{eq:3.48new}
\lim_{r\to\infty} h(r)=0\ .
\end{equation}
We claim that \eqref{eq:3.47new} and \eqref{eq:3.48new} imply that for a 
suitable choice of natural numbers $i_1<i_2<\cdots$ one has 
\begin{equation}\label{eq:3.49new}
(|x_{i_j,j}|^p)_{j=1}^\infty\text{ is uniformly integrable.}
\end{equation}
To prove \eqref{eq:3.49new} put for $j\in\nat$ 
\begin{equation}\label{eq:3.50new}
G_j = \bigcap_{r=1}^j G_{j,r}
\end{equation}
where for $j,r\in\nat$, 
\begin{equation}\label{eq:3.51new}
G_{j,r} = \left\{i\in\nat \mid d_{L^p(\T)} (x_{ij},r\Ba (\N)) <h(r) 
+ \frac1r\right\}\ .
\end{equation}
By \eqref{eq:3.47new} each $G_{j,r}\in U$, and hence also $G_j\in U$ for 
all $j\in\nat$. 
Since $U$ is a free ultrafilter, each $G_j$ is infinite, so we can choose 
successively $i_1<i_2<\cdots$ such that $i_j\in G_j$ for all $j$. 
Put $y_j = x_{i_j,j}$, $j\in\nat$ and $W= \{y_j,j\in\nat\}$, and put as 
in Corollary~2.7 
\begin{equation}\label{eq:3.52new} 
g_W (r) = \sup_{j\in\nat} d_{L^p(\T)} (y_j,r\Ba (\N))\ ,\qquad r\in\real^+\ .
\end{equation}
To prove \eqref{eq:3.49new} we just have to show that $g_W (r)\to0$ 
when $r\to\infty$ (cf.\ Corollary~2.7). 
Let $\ep>0$. 
By \eqref{eq:3.48new} we can choose $r_0\in\nat$ such that 
\begin{equation}\label{eq:3.53new} 
h(r_0) +\frac1{r_0} <\ep\ .
\end{equation}
When $j\ge r_0$, $i_j \in G_j \subseteq G_{j,r_0}$. 
Hence by \eqref{eq:3.51new} and \eqref{eq:3.53new} 
\begin{equation}\label{eq:3.54new} 
d_{L^p(\T)} (y_j,r_0\Ba (\N)) <\ep\ ,\qquad j\ge r_0\ .
\end{equation}
Since $\N=\bigcup_{r>0} r\Ba (\N)$ is dense in $L^p(\T)$ we have for every 
$j\in\nat$, 
\begin{equation*}
\lim_{r\to\infty} d_{L^p(\T)} (y_j,r\Ba (\N)) =0\ .
\end{equation*}
Hence, we may choose $r_1\ge r_0$, such that 
\begin{equation}\label{eq:3.55new} 
d_{L^p(\T)} (y_j,r_1 \Ba (\N)) < \ep\ ,\qquad j=1,\ldots,r_0-1\ .
\end{equation}
By \eqref{eq:3.54new} and \eqref{eq:3.55new}, $g_W (r)<\ep$ for all $r\ge r_1$.
This shows that $\lim_{r\to\infty} g_W (r)=0$ and hence by Corollary~2.7, 
$(|y_j|^p)_{j=1}^\infty$ is uniformly integrable, i.e., (3.49) holds.  
Thus by the assumption that $(y_j)$ is unconditional, 
Corollary~3.5 yields that 
for any subsequence $(y'_j)$ of $(y_j)$,
\begin{equation}\label{eq:3.56new}
\lim_{n\to\infty} n^{-1/p} \Big\|\sum_{j=1}^n y'_j\Big\|_{L^p(\T)} =0\ .
\end{equation}
Putting now $j_k=k$, we have $y_k = x_{i_k,j_k}$ and Theorem~\ref{thm:1.1} 
follows.
\end{proof}

\section{Improvements to the Main Theorem}	
\setcounter{equation}{0}

We obtain here results that are stronger than the Main Theorem. 
In particular, Theorem~4.1 is also needed in Section~6 (specifically, for 
the proof of Theorem~6.9). 
The arguments in this section do not use the ultraproduct construction and 
technique of Section~3. 
They are in a sense more elementary, and also more delicate, than 
those of Section~3. 

We use the following terminology: given a matrix $(x_{ij})$, a 
sequence $(x_{i_k,j_k})$ of elements of the matrix is called a 
{\em generalized diagonal\/} if $i_1<i_2<\cdots$ and $j_1<j_2<\cdots$. 
A set $W$ (or matrix $(x_{ij})$) in a Banach space is called 
{\em semi-normalized\/} if there are $0<\delta\le K<\infty$ with 
$\delta\le \|w \| \le K$ for all $w\in W$. 
The main theorem follows also immediately from the following result. 

\begin{thm}\label{thm:3.5} 
Let $\N$ be a finite von-Neumann algebra, $1\le p<2$, and $(x_{ij})$ be 
an infinite semi-normalized matrix in $L^p(\N)$. 
Assume that every column and generalized diagonal is unconditional, and 
there is a $u\ge1$ so that every row is $u$-unconditional. 
Then one of the following three alternatives holds. 
\begin{itemize}
\item[I.] Some column has a subsequence equivalent to the usual $\ell^p$ 
basis.
\item[II.] There is a $C\ge1$ so that for all $n$, there exists a row 
which contains $n$ elements $C$-equivalent to the usual $\ell_n^p$ basis.
\item[III.] There is a generalized diagonal $(y_k)$ so that 
$$n^{-1/p} \Big\| \sum_{i=1}^n y'_i\Big\|_p \to 0\ \text{ as }\ 
n\to\infty$$
for all subsequences $(y'_i)$ of $(y_i)$.
\end{itemize}
\end{thm}

To recover the Main Theorem from \ref{thm:3.5}, 
let $(x_{ij})$ be as in the hypotheses of the Main Theorem, 
and simply note that 
Cases~I and II of \ref{thm:3.5} are impossible, since otherwise one would 
obtain a constant $\lambda$ so that the $\ell_n^p$ and $\ell_n^2$ bases 
are $\lambda$-equivalent for all $n$. 
Case~III now yields the conclusion of the Main Theorem. 

\begin{rem}
Let us say that {\it the rows of $(x_{ij})$ contain $\ell_n^p$-sequences\/} 
if condition~II of \ref{thm:3.5} holds, with a similar definition for the 
columns. 
Since obviously we can interchange rows and columns in the statement of 
\ref{thm:3.5}, we then obtain the following immediate consequence: 
{\it Let $\N$, $p$ and $(x_{ij})$ be as in the first sentence of 
Theorem~\ref{thm:3.5}. 
Assume that every generalized diagonal 
is unconditional and there is a $u\ge1$ so 
that every row and column are $u$-unconditional. 
Then one of the following~holds.}
\begin{itemize}
\item[I.] {\it Some column or some row has a subsequence equivalent to 
the usual $\ell^p$ basis.}
\item[II.] {\it Both the rows and the columns contain $\ell_n^p$-sequences.} 
\item[III.] {\it Condition III of \ref{thm:3.5} holds.}
\end{itemize}
\end{rem}

\section*{Proof of Theorem~\ref{thm:3.5}} 

We may assume without loss of generality that $\|x_{ij}\|_p \le1$ for all 
$i$ and $j$. 
We introduce the following notation, for all $\ep>0$ and all $i,j=1,2,\ldots$.
\begin{eqnarray}
\om_{ij}(\ep)&=& \om_p (x_{ij},\ep) \label{eq:3.28i}\\
\om_j (\ep)  &=& \sup_i \om_{ij} (\ep)\ .\label{eq:3.28ii}
\end{eqnarray}

Now assume that Case I of Theorem~\ref{thm:3.5} does not occur. 
We then have by Corollary~\ref{cor:3.4} (and Lemma~\ref{lem:2.3}) that 
$(|x_{ij}|^p)_{i=1}^\infty$ is uniformly integrable for all $j$, and hence 
\begin{equation}\label{eq:3.29} 
\lim_{\ep\to0} \om_j (\ep) =0\ \text{ for all }\ j\ .
\end{equation}

We now use the following (hopefully intuitive) convention. 
A set of rows $\R$ of $(x_{ij})$ is identified with a set 
$\J \subset \{1,2,\ldots\}$ via $\R = \{R_i:i\in\J\}$ where 
$R_i = \{x_{ij} :j=1,2,\ldots\}$ for all $i\in \J$. 
Columns are just identified with $j\in \nat $; i.e., 
$j\sim C_j = \{ x_{ij} :i=1,2,\ldots\}$. 

\subsection*{Case II}
There is an $\eta>0$ and an infinite set of rows $\J$ so that for all 
further infinite sets of rows $\J'\subset\J$, all $\delta>0$, and all 
columns $j_0$, there is a column $j>j_0$ so that 
\begin{equation}\label{eq:3.30} 
\{i\in\J' :\om_{i,j} (\delta) >\eta\}\ \text{ is infinite.}
\end{equation}

Intuitively, the final statement means that looking down the $j^{th}$ 
column of the submatrix with rows $\J'$, then infinitely many of the moduli 
$\om_{i,j}(\delta)$ are bigger than $\eta$. 

We shall show that Case II yields II of Theorem~\ref{thm:3.5}. 
In fact, we shall show that then, via Lemma~\ref{lem:3.3}, 
\begin{equation}\label{eq:3.31}
\left\{
\begin{split}
&\text{for every $n$, there exists a row $R_i$ and elements $x_{ij_1}, \ldots,
x_{ij_n}$ in}\\
&\text{$R_i$, $j_1<\cdots <j_n$, with $(x_{ij_k})_{k=1}^n$ 
$\frac{7u}{\eta}$-equivalent to the $\ell_n^p$ basis.}
\end{split}
\right.
\end{equation}

Let $\J_0$ be the initial set of rows satisfying Case~II. 
Let $\delta_1 =1/2$, and choose $j_1$ so that 
\begin{equation}\label{eq:3.32} 
\J_1 \defeq \{i\in\J_0 :\om_{ij_1} (\delta_1) >\eta\}\ \text{ if infinite.}
\end{equation}

Next, using \eqref{eq:3.29}, choose $\bar\delta_2<\delta_1$ so that 
\begin{equation}\label{eq:3.33}
\om_{j_1} (\bar\delta_2) <\frac{\ep}{2}\ ,
\end{equation}
and choose $\delta_2 <\bar\delta_2$.
Now using the assumptions of Case~II, choose $j_2>j_1$ so that 
\begin{equation}\label{eq:3.34}
\J_2 \defeq \{i\in \J_1 :\om_{ij_2} (\delta_2) >\eta\}
\ \text{ is infinite.}
\end{equation}

For the general inductive step, suppose $n>1$, infinite $\J_1\supset\cdots
\supset \J_{n-1}$ and $j_1<\cdots < j_{n-1}$, 
$\delta_1 >\bar\delta_2 >\delta_2 >\cdots > \bar\delta_{n-1}>\delta_{n-1}>0$ 
have been chosen so that for all $1\le \ell<n-1$, 
$\om_{j_\ell}(\bar\delta_{\ell+1} ) <\frac{\ep}2$ and 
$\delta_{\ell+1} +\cdots + \delta_{n-1} <\bar\delta_{\ell+1}$. 
Using \eqref{eq:3.29}, choose $0<\bar\delta_n<\delta_{n-1}$ so that 
$\om_{j_{n-1}}(\bar\delta_n) <\frac{\ep}2$; then choose 
$0<\delta_n<\bar\delta_n$ so that also $\delta_{\ell+1}+\cdots + \delta_n 
<\bar\delta_{\ell+1}$ for all $1\le \ell < n-1$. 
We thus have that 
\begin{equation}\label{eq:3.35}
\om_{j_\ell} (\delta_{\ell+1}+\cdots + \delta_n) <\frac{\ep}2
\ \text{ for all }\ 1\le \ell\le n-1\ .
\end{equation}
Then choose $j_n>j_{n-1}$ so that 
\begin{equation}\label{eq:3.36} 
\J_n \defeq \{i\in\J_{n-1} :\om_{ij_n} (\delta_n) >\eta\}
\ \text{ is infinite.}
\end{equation}

This completes the inductive construction. 
Now fix $n$, let $i\in \J_n$, and let $f_k = x_{ij_k}$ for $1\le k\le n$.
Then $(f_1,\ldots,f_n)$ satisfies the assumption of Lemma~\ref{lem:3.3}. 
Indeed, the $f_i$'s are $u$-unconditional by hypothesis, and for each 
$k$, $1\le k\le n$ 
\begin{equation}\label{eq:3.37i}
\om_{ij_k} (\delta_k) = \om_p (f_k,\delta_k) >\eta 
\end{equation}
and 
\begin{equation}\label{eq:3.37ii}
\om_p (f_k,\delta_m + \delta_{m+1}+\cdots + \delta_n) 
\le \om_{j_k} (\delta_m + \delta_{m+1}+\cdots + \delta_n) <\frac{\ep}{2}
\ \text{ for }\ k<m\le n\ .
\end{equation}
Thus $(x_{ij_k})_{k=1}^n$ satisfies the conclusion of \eqref{eq:3.31} in 
view of Lemma~\ref{lem:3.3}, proving Case~II of \ref{thm:3.5} holds. 

We now suppose that Case II does not hold, i.e., we have 

\subsection*{Case III} 
For all $\eta>0$ and infinite sets of rows $\J$, there exists an infinite 
set of rows $\J'\subset \J$, a $\delta>0$, and a column $\bj$ so that 
for all columns $j\ge \bj$, 
\begin{equation}\label{eq:3.38} 
\om_{i'j}(\delta)\le \eta \text{ for all but finitely many } i'\in\J'\ .
\end{equation}
(Note that we get $j\ge \bj$ instead of $j>\bj$ by just replacing $\bj$ 
by $\bj+1$). 

Intuitively, the final statement means that now, looking down the $j^{th}$ 
column of the submatrix with rows $\J'$, then all but finitely many of the 
moduli $\om_{i',j}(\delta)$ are no bigger than $\eta$. 

We shall now construct $i_1<i_2<\cdots$ and $j_1<j_2<\cdots$ so that 
\begin{equation}\label{eq:3.39} 
\lim_{\ep\to0} \sup_k \om_{i_k,j_k} (\ep) =0\ .
\end{equation} 
Thus we obtain that $(|x_{i_kj_k}|^p)_{k=1}^\infty$ is uniformly 
integrable, and hence Case~III of Theorem~\ref{thm:3.5} holds by 
Corollary~3.5. 

We first claim that we may choose infinite sets of rows $\J_1\supset 
\J_2\supset\cdots$, columns $j_1<j_2<\cdots$, and numbers 
$1\ge \delta_1$, $\frac12 \ge \delta_2$, $\frac13\ge \delta_3\cdots$ 
so that for all $k$, 
\begin{equation}\label{eq:3.40} 
\text{ for all $j\ge j_k$, $\om_{ij}(\delta_k)\le \frac1{2^k}$ for all 
but finitely many } i\in \J_k\ .
\end{equation}

Indeed, first choose $\J_1$ an infinite set of rows, $j_1\in \nat$ and 
$\delta_1>0$ so that for all $j\ge j_1$, \eqref{eq:3.38}  holds, where
$\J'=\J$, $\eta=1/2$, and $\delta_1=\delta$. 

Now suppose $\J_k$, $j_k$, and $\delta_k$ have been chosen. 
Setting $\eta = 1/2^{k+1}$, choose an infinite $\J_{k+1}\subset \J_k$, 
$\bj >j_k$ and a $\delta>0$ so that for all $j\ge \bj$, 
\eqref{eq:3.38} holds for $\J' = \J_{k+1}$. 
Now simply let $\delta_{k+1} = \min \{\delta,2^{-1}\delta_k ,\frac1{k+1}\}$. 
Since the functions $\om_{i\ell}$ are non-decreasing, we have that also for 
all $j>\bj$, $\om_{ij}(\delta_{k+1}) \le 1/2^{k+1}$ for all but 
finitely many  $i\in \J_{k+1}$. 
This completes the inductive construction, with \eqref{eq:3.40} holding 
for all $k$. 

Now choose $i_1\in \J_1$ with $\om_{i_1,j_1}(\delta_1) \le 1/2$. 
Then also for all but finitely many $i\in \J_2$, $\om_{i,j_2}(\delta_1) 
\le 1/2$ and $\om_{i,j_2}(\delta_2)\le 1/4$. 
Hence we can choose $i_2 >i_1$ $(i_2 \in \J_2)$, with 
\begin{equation}\label{eq:3.41} 
\om_{i_2,j_2} (\delta_1) \le \frac12\ \text{ and }\ \om_{i_2,j_2} 
(\delta_2 ) \le \frac14\ .
\end{equation}
But we can also choose $0<\ep_2\le\delta_2$ so that 
\begin{equation}\label{eq:3.42} 
\om_{i_1,j_1}(\ep_2) \le \frac14\ .
\end{equation}
Thus also 
\begin{equation}\label{eq:3.43} 
\om_{i_2,j_2} (\ep_2) \le \frac14\ .
\end{equation}

Now suppose $i_1<\cdots <i_n$ and $\delta_1 =\ep_1,\ldots,\ep_n$ have been 
chosen so that $\ep_j\le\delta_j$ for all $j\le n$ and 
\begin{equation}\label{eq:3.44} 
\om_{i_k,j_k} (\ep_i) \le \frac1{2^i} \ \text{ for all }\ 
1\le k\le n,\ 1\le i\le n\ .
\end{equation}
Now by \eqref{eq:3.40}, choose $i_{n+1} >i_n$ ($i_{n+1}\in\J_{n+1}$) so that 
\begin{equation}\label{eq:3.45} 
\om_{i_{n+1},j_{n+1}} (\delta_\ell) \le \frac1{2^\ell} \ \text{ for all }\ 
1\le \ell \le n+1\ .
\end{equation}
This is possible, since for each $\ell$, $\om_{i,j_{n+1}} (\delta_\ell) 
\le 1/2^\ell$ for all but finitely many $i\in \J_{n+1}$.

Again, since the $\ep_\ell$'s are smaller than the $\delta_\ell$'s, 
\begin{equation}\label{eq:3.46} 
\om_{i_{n+1}, j_{n+1}} (\ep_\ell) \le \frac1{2^\ell}\ \text{ for all } \ 
1\le \ell\le n\ .
\end{equation} 
Finally, choose $\ep_{n+1} \le \delta_{n+1}$ so that 
\begin{equation}\label{eq:3.47} 
\om_{i_\ell,j_\ell} (\ep_{n+1}) \le \frac1{2^{n+1}} \ \text{ for all }\ 
1\le \ell\le n\ .
\end{equation}
Again, we also have 
\begin{equation}\label{eq:3.48} 
\om_{i_{n+1},j_{n+1}} (\ep_{n+1}) \le \frac1{2^{n+1}}\ .
\end{equation}

This completes the inductive construction of $i_1<i_2<\cdots$ and 
$\ep_1,\ep_2,\ldots$. 
Then for each $i$, we have 
\begin{equation}\label{eq:3.49} 
\sup_k \om_{i_k,j_k} (\ep_i) \le \frac1{2^i}\ .
\end{equation} 

It then follows immediately that \eqref{eq:3.39} holds, since if 
$\ep\le \ep_i$, then also 
\begin{equation}\label{eq:3.50} 
\sup_k \om_{i_k,j_k} (\ep) \le \frac1{2^i}\ .
\end{equation}
This completes the proof of Theorem~\ref{thm:3.5}, in view of the 
comment after \eqref{eq:3.39}.\qed
\medskip

Using theorems from Banach space theory, we next obtain a stronger 
version of \ref{thm:3.5}. 

\begin{thm}\label{thm:3.6}
Let $\N$, $p$ and  $(x_{ij})$ be as in \ref{thm:3.5}. 
The conclusion of \ref{thm:3.5} holds under the following assumptions: 
\begin{itemize}
\item[(a)] $1<p$, and every column is an unconditional basic sequence, 
every generalized diagonal is a basic sequence, and there is a $\lambda\ge1$ 
so that every row is a $\lambda$-basic sequence 
\end{itemize}
or
\begin{itemize}
\item[(b)] $p=1$ and every generalized diagonal is a basic sequence.
\end{itemize}

Moreover the unconditional assumption in {\rm (a)} may be dropped if $\N$ 
is hyperfinite. 
\end{thm}

\begin{rem}
We do not know if the unconditional assumption in (a) may be dropped in 
general. 
However our proof of \ref{thm:3.6} yields the following result, for 
arbitrary finite $\N$ and $1<p<2$. 
{\it Assume {\rm (a)} with ``unconditional'' deleted. 
Then the following three alternatives hold:  {\rm II} or {\rm III} of 
Theorem~\ref{thm:3.5}, or \newline 
{\rm I}$'$. 
There is a $C\ge1$ and a column so that for all $n$, the column contains 
$n$ elements $C$-equivalent to the usual $\ell_n^p$ basis\/}. 
\end{rem}

To obtain the case $p>1$, we require the following remarkable result, due 
to Brunel and Sucheston (\cite{BrS1}, \cite{BrS2}; see also \cite{G}). 
(A sequence $(x_j)$ of non-zero elements in a Banach space is called 
{\it $\beta$-suppression unconditional\/} 
if for all $n$, scalars $c_1,\ldots,c_n$, 
and $F\subset \{1,\ldots,n\}$, $\|\sum_{j=F} c_j x_j\| 
\le \beta \|\sum_{j=1}^n c_j x_j\|$. 
It is easily seen that if $(x_j)$ is $\lambda$-suppression unconditional, 
it is $2\lambda$-unconditional over real scalars and $4\lambda$-unconditional 
over complex scalars. 
Actually, a neat result of Kaufman-Rickert yields that such a sequence 
is $\pi \lambda$-unconditional (over complex scalars) \cite{KR}.) 

\begin{lem}\label{lem:3.7} 
Let $(x_n)$ be a semi-normalized weakly null sequence in a Banach space $X$, 
and let $\ep>0$. 
Then there exists a subsequence $(y_j)$ of $(x_j)$ so that for any 
$k\le j_1<j_2<\cdots <j_{2^k}$, $(y_{j_i})_{i=1}^{2^k}$ is 
$(1+\ep)$-suppression unconditional (and hence $\pi   (1+\ep)$-unconditional). 
\end{lem}

\begin{rems} 
1. Actually, the results of Brunel-Sucheston yield much more than this. 
They obtain that under the hypotheses of Lemma~\ref{lem:3.7}, there 
exists a Banach space $E$ with a suppression 1-unconditional semi-normalized 
basis $(e_j)$ and a basic subsequence $(y_j)$ of $(x_j)$ so that: 
\begin{itemize}
\item[(i)] $(e_j)$ is isometrically equivalent to all of its subsequences and
\item[(ii)] for all $\ep>0$ and $k$ large enough, and any $k\le j_1<\cdots 
< j_{2^k}$, $(y_{j_i})_{i=1}^{2^k}$ is $(1+\ep)$-equivalent to 
$(e_1,\ldots,e_{2^k})$. 
\end{itemize}
In the standard Banach space terminology, $(e_j)$ is called a subsymmetric 
basis for $E$, and a {\it spreading model\/} for $(x_j)$.

2. A classical result of Bessaga-Pe{\l}czy\'nski yields that any seminormalized 
weakly null sequence in a Banach space has a basic subsequence (in fact, 
for every $\ep>0$, a subsequence which is $(1+\ep)$-basic). 
However it is obtained in \cite{MR} that there exists a normalized weakly 
null sequence in a certain Banach space with no unconditional subsequence, 
and in \cite{GM} that there exists an (infinite dimensional) reflexive 
Banach space with no (infinite) unconditional basic sequences at all. 
Thus in a sense, Lemma~\ref{lem:3.7} is the best possible positive result 
in this direction. 
\end{rems}

We now give consequences of this lemma that are needed for 
Theorem~\ref{thm:3.6}. 
The first one follows from Lemma~\ref{lem:3.1} and Lemma~4.3.

\begin{cor}\label{cor:3.8}
Let $1\le p<2$ and $(f_n)$ be a weakly null sequence in $L^p(\T)$ so that 
$(|f_i|^p)_{i=1}^\infty$ is uniformly integrable. 
Then there is a subsequence $(f'_i)$ of $(f_i)$ so that 
\begin{equation*}
\lim_{n\to\infty} n^{-1/p} \|\ep_1 y_1 +\cdots + \ep_n y_n\|_{L^p(\T)}=0
\end{equation*}
uniformly over all subsequences $(y_i)$ of $(f'_i)$ and all choices 
$(\ep_j)$ of scalars with $|\ep_j| \le1$ for all $j$. 
\end{cor}

\begin{rem}
The result shows (and also follows from): 
{\it any spreading model for $(f_j)$ is not equivalent to the $\ell^p$-basis.}
\end{rem}

\begin{proof}[Proof of \ref{cor:3.8}] 
We may assume without loss of generality that $\|f_j\|_p \le1$ for all $j$. 
Let $\ep>0$ be such that $\pi (1+\ep)\le4$, and choose $(y_j)$ a 
subsequence of $(f_j)$ satisfying the conclusion of Lemma~\ref{lem:3.7}. 
Let $(r_j)$ denote the Rademacher functions on $[0,1]$ (as defined in 
Section~3), set $\tilde \N = \N\bar\otimes L^\infty$, and let $g_j=
y_j \otimes r_j$ for all $j$. 
Then $(g_j)$ is 2-unconditional (over complex scalars) and of course 
$(|g_j|^p)$ is also uniformly integrable in $L^1(\tilde\N)$, whence 
by Lemma~\ref{lem:3.1},  
\begin{equation}\label{eq:3.44a}
\lim_{n\to\infty} n^{-1/p} \|g_1+\cdots + g_n\|_{L^p(\tilde\N)} =0\ .
\end{equation} 
Let $\ep>0$, and choose $N$ so that if $n\ge N$, then 
\begin{equation}\label{eq:3.45a}
n^{-1/p} \|g_1+\cdots + g_n\|_{L^p(\tilde\N)} < \frac{\ep}{16}
\end{equation}
and
\begin{equation}\label{eq:3.46a}
n^{-1/p} (1+ \log_2 n) <\frac{\ep}2\ .
\end{equation} 
Now fix $n$, and choose $k$ with 
\begin{equation}\label{eq:3.47a}
2^{k-1} \le n<2^k\ .
\end{equation}
Of course then 
\begin{equation}\label{eq:3.48a}
k\le 1+ \log_2 n\ .
\end{equation}
Now if $\ep_1,\ldots,\ep_n$ are given scalars of modulus at most one, then 
\begin{equation}\label{eq:3.49a}
\Big\|\sum_{j=k+1}^n \ep_j y_j\Big\|_{L^p(\N)} 
\le 16 \Big\| \sum_{j=k+1}^n g_j\Big\|_{L^p(\tilde \N)}\ .
\end{equation}
Indeed, $y_{k+1},\ldots,y_n$ is 4-unconditional by the conclusion 
of Lemma~\ref{lem:3.7} (since $n-k<n<2^k$), yielding \eqref{eq:3.49a}. 
On the other hand, 
\begin{equation}\label{eq:3.50a} 
\Big\|\sum_{j=1}^k \ep_j y_j \Big\|_{L^p(\N)} 
\le k\le 1+\log_2 n\ \text{ by \eqref{eq:3.48a}.}
\end{equation}
Thus we have 
\begin{equation}\label{eq:3.51} 
\begin{split}
n^{-1/p} \Big\|\sum_{j=1}^n \ep_j y_j\Big\|_p 
& \le n^{-1/p} \Big\|\sum_{j=1}^k \ep_j y_j\big\|_p 
+ n^{-1/p} \Big\| \sum_{j=k+1}^n \ep_j y_j\Big\|_p\\
\noalign{\vskip6pt}
&\le n^{-1/p} (1+\log_2 n) + 8 n^{-1/p} 
\Big\|\sum_{j=k+1}^n g_j\Big\|_{L^p(\tilde\N)} \\
\noalign{\vskip6pt}
&\le \frac{\ep}2 + 8n^{-1/p} \Big\|\sum_{j=1}^n g_j\Big\|_{L^p(\tilde \N)}\\
\noalign{\vskip6pt}
&\le \frac{\ep}2 + \frac{\ep}2 = \ep\ .
\end{split}
\end{equation}
(The last inequality holds by \eqref{eq:3.45a}; 
the next to the last by \eqref{eq:3.46a} and the fact that 
$(g_j)$ is 1-unconditional over real scalars.) 
The uniformity of the limit over all {\it subsequences\/} of $(y_i)$ 
follows from the fact that the limit in \eqref{eq:3.44a} is uniform 
over all subsequences of $(g_i)$, thanks to the proof of 
Lemma~\ref{lem:3.1}. 
\end{proof}

We next note a general consequence of Lemma~\ref{lem:3.7}, which 
follows from ultraproducts. 

\begin{cor}\label{cor:3.9} 
Let $X$ be a uniformly convex Banach space and let $\lambda\ge1$, $\ep>0$, 
and $k$ be given. 
Then there is an $n\ge k$ so that for any $\lambda$-basic sequence 
$(x_1,\ldots,x_n)$ in $X$, there exist $1\le j_1<j_2<\cdots < j_k$ so that 
$(x_{j_1},\ldots,x_{j_k})$ is suppression $(1+\ep)$-unconditional 
(and hence $\pi (1+\ep)$-unconditional). 
\end{cor}

\begin{proof}
Suppose the conclusion were false. 
Then we could find for every $n\ge k$, an $n$-tuple $(x_1^n,\ldots,x_n^n)$ 
of elements in $X$ so that $(x_1^n,\ldots,x_n^n)$ is $\lambda$-basic, 
yet no $k$ terms are suppression $(1+\ep)$-unconditional. 
By homogeneity, we may assume that $\|x_i^n\|=1$ for all $n$ and $i\le n$. 
Now let $\U$ be a non-trivial ultrafilter on $\nat$ and let $X_{\U}$ 
denote the ultrapower of $X$ with respect to $\U$. 
(That is, we let $E_{\U}$ denote the subspace of $\ell^\infty (X)$ 
consisting of all bounded sequences $(x_j)$ in $X$ with 
$\lim_{j\in\U} \|x_j\|=0$, and then set $X_{\U} = \ell^\infty (X)/E_{\U}$.)
Since $X$ is uniformly convex, so is $X_{\U}$. 
Now define a sequence $(\tilde x_j)$ in $X_{\U}$ by 
$\tilde x_j = \pi   (x_j^n)_{n=1}^\infty$, for all $j$, where 
$\pi   :\ell^\infty (x) \to X_{\U}$ is the quotient map and we set 
$x_j^n=0$ if $n<j$. 
It then follows that $(\tilde x_j)$ is also $\lambda$-basic and 
normalized; since $X_{\U}$ is {\it reflexive\/}, $(\tilde x_j)$ is 
weakly null. 
But then by Lemma~\ref{lem:3.7}, there exist $k$ terms 
$\tilde x_{j_1},\ldots,\tilde x_{j_k}$ of this sequence with 
$(\tilde x_{j_i})_{i=1}^k$ $(1+\frac{\ep}2)$-suppression unconditional. 
Standard ultraproduct techniques yield that $\eta>0$ given, there exists 
an $n>j_k$ so that $(\tilde x_{j_1},\ldots,\tilde x_{j_k})$ is 
$(1+\eta)$-equivalent to $(x_{j_1}^n,\ldots,x_{j_k}^n)$ and hence the 
latter is $(1+\eta)$ $(1+\frac{\ep}2)$-suppression unconditional. 
Of course we have a contradiction if $(1+\eta)(1+\frac{\ep}2) <1+\ep$.
\end{proof}

\begin{proof}[Proof of Theorem~\ref{thm:3.6}(a)]

We use the same notations and assumptions as in the proof of 
Theorem~\ref{thm:3.5} (e.g., we assume that $\|x_{ij}\|_p \le1$ for all 
$i$ and $j$). 
Assume that Case~I of \ref{thm:3.5} does not occur. 
Then again we have that $(|x_{ij}|^p)_{i=1}^\infty$ is uniformly 
integrable for all $j$, and hence \eqref{eq:3.22} holds. 
This is also the case if $\N$ is hyperfinite and the unconditional 
assumption in (a) is dropped. 
For suppose to the contrary that for some $i$, $(f_j) \defeq (x_{ij})$ 
has the property that $(|f_j|^p)$ is not uniformly integrable. 
Then setting $g_j = f_j\otimes r_j$ in $L^p(\tilde \N)$ (as defined in 
the proof of Corollary~\ref{cor:3.8}), $(g_j)$ is unconditional and 
again $(|g_j|^p)$ is not uniformly integrable, hence there exist 
$n_1<n_2<\cdots$ with $(g_{n_j})$ equivalent to the usual $\ell^p$-basis, 
by Corollary~\ref{cor:3.4}). 
But  $(f_{n_j})$ has an unconditional subsequence $(f'_j)$ by 
\cite{SF}, \cite{PX1}. 
Of course then $(f'_j)$ is equivalent to $(g'_j) \defeq (f'_j\otimes r_j)$, 
a subsequence of $(g_{n_j})$, whence $(f'_j)$ is equivalent to the 
$\ell^p$ basis. 

Now replace the entire matrix $(x_{ij})$ by 
$(\tilde x_{ij}) \defeq (x_{ij} \otimes r_{ij})$ in $L^p (\tilde\N)$ 
(where $\tilde\N = \N\bar\otimes L^\infty$), 
where $r_{ij}$ is just a ``renumbering'' of $(r_j)$ via $\nat\times\nat$ 
(precisely, let $\varphi :\nat\times\nat\to\nat$ be a bijection, and 
set $r_{ij} = r_{\varphi(i,j)}$). 
Now $\om_p (x_{ij},\ep) = \om_p (\tilde x_{ij},\ep)$ for all $i,j$, 
and $\ep$; hence assuming Case~II in the proof of Theorem~\ref{thm:3.5} 
occurs, we have that II of \ref{thm:3.5} holds for the matrix 
$(\tilde x_{ij})$. 
But then since $L^p(\N)$ is uniformly convex, II holds for $(x_{ij})$ 
itself, by Corollary~\ref{cor:3.9}. 
Indeed, let $C$ be as in II of \ref{thm:3.5}, let $k$ be given. 
Choose $n\ge k$ satisfying the conclusion of \ref{cor:3.9} for 
$X= L^p(\N)$ (with $\pi (1+\ep)\le 4$, say). 
Choose $i$ and $m_1<\cdots < m_n$ with $(\tilde x_j)_{j=1}^n$ 
$C$-equivalent to the $\ell_n^p$ basis where  we set $x_j = x_{im_j}$ 
and $\tilde x_j = \tilde x_{im_j}$ for all $j$. 
Then choose $j_1<\cdots j_k$ with $(x_{j_i})$ 4-unconditional. 
But then $(x_{j_i})$ is 8-equivalent to $(\tilde x_{j_i})$, and is 
hence $8C$-equivalent to the $\ell_k^p$ basis. 

If Case II in the proof of \ref{thm:3.5} does not occur, we have by 
Case~III that there exists a generalized diagonal 
$(\tilde x_{i_n,j_n})_{n=1}^\infty$ of $(x_{ij})$ so that 
$(|\tilde x_{i_n,j_n}|^p)_{n=1}^\infty$ is uniformly integrable. 
Hence immediately, $(|x_{i_n,j_n}|^p)_{n=1}^\infty$ is uniformly 
integrable, and so by Corollary~\ref{cor:3.8}, $(x_{i_n,j_n})$ 
has a subsequence $(y_k)$ (which is of course also a generalized diagonal) 
satisfying III of \ref{thm:3.5}. 
This completes the proof of Theorem~\ref{thm:3.6}(a).\qed 

To obtain \ref{thm:3.6}(b), we need two further ``Banach'' properties 
of preduals of von~Neumann algebras. 
The first one holds in complete generality. 

\begin{lem}\label{lem:3.10} 
Let $\M$ be a von-Neumann algebra, and let $(f_n)$ be a bounded sequence 
in $\M_*$ such that $(f_n)$ is not relatively weakly compact. 
Then $(f_n)$ has a subsequence equivalent to the $\ell^1$-basis.
\end{lem} 

We give a ``quantitative'' proof of this result at the end of this 
section, using the case for commutative $\N$ established in \cite{R1}. 
In fact, Lemma~\ref{lem:3.10} is due to H.~Pfitzner \cite{Pf}. 
However, the second result we need is a ``localization'' of our proof, 
which does not seem to follow directly from previously known material. 
This result yields that if $n$ elements  of $\Ba (\N_*)$ ($\N$ finite) 
have mass at least $\eta$ on pairwise orthogonal projections, 
then $k$ of these are 
$C$-equivalent to the $\ell_k^1$-basis. 
Here, $C$ depends only on $\eta$, $n$ on $k$ and $\eta$. 
To make this more manageable, let us simply say that 
{\it $n$ elements $f_1,\ldots, f_n$ of the predual of a von-Neumann algebra 
$\M$ are $\eta$-disjoint provided there exist pairwise orthogonal projections 
$P_1,\ldots,P_n$ in $\M$ such that} 
\begin{equation}\label{eq:3.51a}
\|P_i f_iP_i\|_1 \ge \eta\ \text{ for all }\ i\ .
\end{equation}
(Here, if $P\in\M$ and $f\in\M_*$, $PfP$ is defined by: 
$\langle T,PfP \rangle = \langle PTP,f\rangle$ for all $T\in\M$. 
Also, $\|\cdot\|_1$ denotes the predual norm on $\M_*$.) 
(We shall also say $f_1,\ldots,f_n$ {\it are disjoint\/} provided there are 
pairwise orthogonal projections $P_1,\ldots,P_n$ 
in $\M$ with $f_i = P_i f_iP_i$ for all $i$. 
Evidently if the $f_i$'s are normalized, they are disjoint iff they 
are 1-disjoint.) 

\begin{lem}\label{lem:3.11} 
Given $\eta>0$, then if $C>\frac1{\eta}$, then for all $k\ge1$, there is an 
$n\ge k$ so that for any von-Neumann algebra $\N$ and $\eta$-disjoint 
elements $f_1,\ldots,f_n$ in $\Ba (\N_*)$, there exist $j_1<\cdots <j_k$ 
with $(f_{j_i})_{i=1}^k$ $C$-equivalent to the $\ell_k^1$ basis. 
\end{lem}

We delay the proof of this result, and complete the proof of 
Theorem~\ref{thm:3.6}, i.e., the case $p=1$. 
Again we make the same assumptions and use the same notation as in the 
proof of \ref{thm:3.5}(a). 
Now suppose that Case~I of Theorem~\ref{thm:3.5} does not occur. 
We now have, immediately from Proposition~\ref{prop:2.2} and 
Lemma~\ref{lem:3.10}, that $(x_{ij})_{j=1}^\infty$ is uniformly integrable 
for all $j$, and hence again \eqref{eq:3.22} holds. 
Now again assume Case~II of the proof \ref{thm:3.5} holds. 
Then the proof of \ref{thm:3.5}II yields that for all $n$, there exists 
a row $i$ and $j_1<\cdots < j_n$ so that $(f_k)_{k=1}^n$ is 
$\frac{\eta}3$-disjoint, where $f_k=x_{ij_k}$ for all $k$. 

Indeed, we obtain there (following formula \eqref{eq:3.29}), that for all 
$n$, there is a sequence $(f_1,\ldots,f_n)$ satisfying the assumptions 
of Lemma~\ref{lem:3.3} (for $\eta >0$ and $0<\ep<\frac{\eta}2$) 
{\it except\/} for the $u$-unconditionality assumption. 
But the {\it proof\/} of Lemma~\ref{lem:3.3} yields precisely that 
$(f_1,\ldots,f_n)$ is $\frac{\eta}2 - \ep$ disjoint; the unconditionality 
assumption was only used, in invoking Lemma~\ref{lem:3.2}. 
Of course we may choose $\ep = \frac{\eta}6$, and so $(f_1,\ldots,f_n)$ 
is then $\frac{\eta}3$-disjoint.

Then Lemma~\ref{lem:3.11} immediately yields Case~II of Theorem~\ref{thm:3.5}. 
Finally, assuming Case~II of the proof of \ref{thm:3.5} does not occur, 
we obtain again from the proof of Case~III that there exists a 
generalized diagonal $(g_k)$ of $(x_{ij})$ with $(g_k)$ uniformly integrable. 
Hence there exists a weakly convergent subsequence $(f_j)$ of $(g_k)$, 
by Proposition~\ref{prop:2.2}. 
But since we assume the generalized diagonals of $(x_{ij})$ are basic 
sequences, $(f_j)$ must be weakly null. 
Now Corollary~\ref{cor:3.8} immediately yields Case~III of 
Theorem~\ref{thm:3.5}.
\end{proof}

\begin{rem} 
The case $p=1$ of Theorem~\ref{thm:3.6} may be alternatively formulated 
as follows (with essentially no assumptions at all on the matrix 
$(x_{ij})$). 
\end{rem}

\begin{bprime}
Let $\N$ be a finite von-Neumann algebra and let $(x_{ij})$ be an infinite 
semi-normalized matrix in $\N_*$. 
Then one of the following holds.  
\begin{itemize}
\item[I.] Some column has a subsequence equivalent to the usual $\ell^1$ basis.
\item[II.] There is a $C\ge1$ so that for all $n$, there exists a row with 
$n$ elements $C$-equivalent to the usual $\ell_n^1$ basis.
\item[III.] Some generalized diagonal of $(x_{ij})$ is weakly convergent.
\end{itemize}
\end{bprime}

It remains to prove Lemma~\ref{lem:3.11}. 
This is an immediate consequence of the following two results, which 
in turn follow from the techniques in \cite{R1}. 
(We denote the ``predual norm'' of a general von-Neumann algebra 
by $\|\cdot\|_1$.) 

\begin{lem}\label{lem:3.12} 
Let $\N$ be an arbitrary von-Neumann algebra, and $f_1,f_2,\ldots$ be 
a finite or infinite sequence in $\N_*$ with $\|f_i\|_1\le1$ for all $i$. 
Assume there are pairwise orthogonal projections $P_1,P_2,\ldots$ in $\N$ 
and $0<\ep<\delta \le1$ so that for all $i$, 
\begin{equation}\label{eq:3.52} 
\|P_i f_iP_i\|_1 \ge \delta\ \text{ and }\ 
\sum_{j\ne i} \|P_j f_i P_j\|_1 \le \ep\ .
\end{equation}
Then $f_1,f_2,\ldots$ is $\frac1{\delta-\ep}$ equivalent to the usual basis 
of $\ell^1$ (resp. $\ell^1_n$ if the sequence has $n$ terms). 
\end{lem} 

\begin{lem}\label{lem:3.13} 
Let $k\ge1$ and $0<\ep<1$ be given. 
There is an $n\ge k$ so that given any von Neumann algebra $\N$, 
$f_1,\ldots,f_n\in \Ba(\N_*)$, and pairwise orthogonal projections 
$P_1,\ldots,P_n$ in $\N$, there exist $j_1<j_2<\cdots < j_k$ so that 
for all $1\le i\le k$, 
\begin{equation}\label{eq:3.53} 
\sum_{r\ne i} \|P_{j_r} f_{j_i} P_{j_r}\|_1 <\ep\ .
\end{equation}
\end{lem}

\begin{rem}
We obtain that we may choose $n= k^\ell$ where $\ell = [1/\ep] +1$. 
\end{rem}

\begin{proof}[Proof of Lemma~\ref{lem:3.11}] 
Let $C>\frac1{\eta}$ and choose $0<\ep<\eta$ with $\frac1{\eta-\ep} < C$. 
Let $n$ be as in Lemma~\ref{lem:3.13}, $f_1,\ldots,f_n$ as in the 
hypotheses of \ref{lem:3.11}, and choose $j_1,\ldots,j_k$ satisfying 
the conclusion of \ref{lem:3.13}. 
Then $(f_{j_i})_{i=1}^k$ is $C$-equivalent to the $\ell^1_k$ basis 
by Lemma~\ref{lem:3.12}.
\end{proof}

\begin{proof}[Proof of Lemma~\ref{lem:3.12}] 
Let $n<\infty$ be less than or equal to the number of terms in the sequence,  
and let $c_1,\ldots,c_n$ be given scalars with 
\begin{equation}\label{eq:3.54} 
\sum_{i=1}^n |c_i| = 1\ .
\end{equation}
Let $g= \sum_{i=1}^n c_i f_i$. 
Since the $P_j$'s are pairwise orthogonal, we have that 
\begin{equation}\label{eq:3.55} 
\|g\|_1 \ge \sum_{j=1}^n \|P_j g P_j\|_1\ .
\end{equation}
Now fixing $j$, 
\begin{equation}\label{eq:3.56}
\begin{split}
\|P_j gP_j\|_1 
& \ge \|P_j c_j f_j P_j + P_j \sum_{i\ne j} c_i f_i P_j\|_1\\
& \ge |c_j|\delta - \sum_{i\ne j} |c_i|\, \|P_j f_i P_j\|_1
\end{split}
\end{equation}
by \eqref{eq:3.52} and the triangle inequality. 
Hence using \eqref{eq:3.55} and \eqref{eq:3.56}, 
\begin{equation}\label{eq:3.57}
\begin{split}
\|g\|_1 &\ge \sum_{j=1}^n |c_j|\delta - \sum_{j=1}^n \sum_{i\ne j} 
|c_i|\, \|P_j f_i P_j\|_1\\
& = \delta - \sum_{i=1}^n |c_i| \sum_{j\ne i} \|P_j f_iP_j\|_1
\ \text{ by \eqref{eq:3.54}}\\
&\ge \delta-\ep\ \text{ by \eqref{eq:3.54} and \eqref{eq:3.52}.}
\end{split}
\end{equation}
This completes the proof.
\end{proof}

We finally deal with Lemma~\ref{lem:3.13}.  
This result follows from the simplest possible setting: 
$\N = \ell_n^\infty$, the $f_i$'s are in $\ell_n^{1+}$ 
(i.e., the positive part of $\N_* =\ell_n^1$), and the orthogonal 
projections $P_i$ correspond to multiplication by $\chix_{\{i\}}$ for all $i$.
That is, we finally have the following elementary disjointness result. 

\begin{lem}\label{lem:3.14}
{\rm A.} Let $f_1,f_2,\ldots$ be a bounded infinite subset of 
$\ellpp$, and let $\ep>0$. 
There exist $n_1<n_2<\cdots$ so that for all $i$, 
\begin{equation}\label{eq:3.58}
\sum_{j\ne i} f_{n_i} (n_j) <\ep\ .
\end{equation}
{\rm B.} Let $k\in\nat$ and $\ep>0$ be given. 
There exists an $N\ge k$ so that given $f_1,\ldots,f_N\in 
\Ba \ell_N^{1+}$, there exist  $n_1<n_2<\cdots < n_k$ so that 
for all $1\le i\le k$, \eqref{eq:3.58} holds.
\end{lem}

\begin{rem}
Part A is a special case of Lemma~1.1 of  \cite{R1}. 
Part~B appears to be new. 
We obtain in fact that we may let $N=k^\ell$ where $\ell = [1/\ep]+1$. 
\end{rem}

\begin{proof}[Proof of Lemma~\ref{lem:3.13}] 
Let $\ep>0$ and $N$ be as in the conclusion of \ref{lem:3.14}B. 
Let the $f_i$'s and $P_i$'s be as in the statement of \ref{lem:3.13}. 
For each $i$, define $\tilde f_i$ in $\ellpp$ by 
$\tilde f_i (j) = \|P_j f_i P_j\|_1$ for all $1\le j\le N$. 
Then 
\begin{equation}\label{eq:3.59} 
\sum_{j=1}^{N} \|P_j f_i P_j\|_1 = 
\|\tilde f_i\|_1 
\le \|f_i\|_1\le1
\end{equation}
for all $i$. 
Now the conclusion of B yields $j_1<\cdots < j_k$ so that 
\begin{equation}\label{eq:3.60} 
\sum_{r\ne i} \tilde f_{j_i} (j_r) <\ep\ \text{ for all }\ 1\le i\le k\ .
\end{equation}
Then $f_{j_1},\ldots, f_{j_k}$ satisfies the conclusion of 
Lemma~\ref{lem:3.13}.
\end{proof}

At last, we give the proof of Lemma~\ref{lem:3.14}. 

We first prove A, using  an argument due to J.~Kupka \cite{Ku}. 
We then adapt this argument to obtain Part~B. 
We regard elements of $\ellpp$ as finite measures on subsets of $\nat$ 
and use the notation:  $f(E) = \sum_{j\in E} f(j)$ for $f\in \ellpp$ 
and $E\subset\nat$. 
Thus, the conclusion of A may be restated: 
{\it There exists an infinite $M\subset \nat $ so that} 
\begin{equation}\label{eq:3.61} 
f_i (M\sim \{i\}) < \ep\ \text{ for all }\ i\in M\ .
\end{equation}
Let $N_1,N_2,\ldots$ be pairwise disjoint infinite subsets of 
$\nat$ with $\nat = \bigcup_{j=1}^\infty N_j$. 

\subsection*{Case I} 
For each $i$, there exists $n_i\in N_i$ so that 
\begin{equation}\label{eq:3.62} 
f_{n_i} (\nat\sim N_i) <\ep\ .
\end{equation}
It then follows that $M= \{n_1,n_2,\ldots\}$ satisfies \eqref{eq:3.61}. 
Indeed, for all $i$, 
\begin{equation}\label{eq:3.63} 
\{n_1,n_2,\ldots,n_{i-1},n_{i+1},\ldots\}\subset \nat \sim N_i
\end{equation}
since the $N_i$ are disjoint, so \eqref{eq:3.61} follows from 
\eqref{eq:3.62} and \eqref{eq:3.63}. 

\subsection*{Case II} 
Case I fails.
Thus we may choose $i_1$ so that 
\begin{equation}\label{eq:3.64}
f_j(\nat \sim N_{i_1}) \ge \ep\ \text{ for all }\ j\in N_{i_1}\ .
\end{equation}
Now repeat the same procedure; let $M_1 = N_{i_1}$, and choose $M_1^1,
M_1^2,\ldots$ disjoint infinite subsets of $M_1$ with 
$M_1=\bigcup_{j=1}^\infty M_1^j$. 
If Case~I fails for $M_1$, we will obtain $M_2\defeq M_1^j$ (for 
some $j$) so that 
\begin{equation}\label{eq:3.65} 
f_j (M_1\sim M_2) \ge \ep\ \text{ for all }\ j\in M_2\ .
\end{equation}
Again divide up $M_2$. 
This ``failure of Case~I'' must terminate before $\ell$ steps, where 
$\|f_j\|_1 <\ell \ep$ for all $j$. 
Indeed, otherwise, we finally obtain 
$\nat = M_0\supset M_1\supset M_2\supset \cdots M_\ell$ and a $j\in M_\ell$ 
with 
\begin{equation}\label{eq:3.66} 
f_j(M_{i-1} \sim M_i) \ge \ep\ \text{ for all }\ i\ ,
\end{equation}
whence $\|f_j\| \ge \ell\ep$, a contradiction. 

\begin{proof}[Proof of Part B] 
Let  $\ell = [1/\ep] +1$ and let $N= k^\ell$. 
Let then $f_1,\ldots,f_N \in \Ba (\ell_N^{1+})$ be given. 
Of course the conclusion of Part~B may be restated: 
{\it There exists an $M\subset \{1,\ldots,N\}$ with $\# M = k$ so that 
\eqref{eq:3.61} holds.}

Let $N_1,\ldots,N_k$ be disjoint subsets of $\{1,\ldots,N\}$, each of 
cardinality $k^{\ell-1}$, and just repeat the argument for Part~A, Case~I. 
If Case~I fails, we repeat again the rest of the argument: 
that is, we find $i_1$ satisfying \eqref{eq:3.64} and set $M_1=N_{i_1}$. 
Now we just choose $M_1^1,\ldots,M_1^k$ disjoint subsets of $M_1$, each 
of cardinality $k^{\ell-2}$; if Case~I fails for $M_1$, we continue as 
before, with $M_2$ satisfying \eqref{eq:3.65} and $M_2\subset M_1$, 
$\# M_2 = k^{\ell-2}$. 
If Case~I fails for $\ell$ steps, we obtain finally $\{1,\ldots,N\} = 
M_0 \supset M_1\supset \cdots M_\ell$ with $\#M_i = k^{\ell-i}$ for all $i$, 
so $\# M_\ell=1$; and for $j$ the unique number of $M_\ell$, \eqref{eq:3.66} 
holds, whence again $\|f_j\|\ge \ell\ep>1$, a contradiction.
\end{proof}

Let us say that a finite or infinite sequence $(f_i)$ satisfying the 
hypotheses of Lemma~\ref{lem:3.12} is {\it $(\delta,\ep)$-relatively 
disjoint\/}. 
It then follows from arguments in \cite{R1} that 
{\it the closed linear span of such a sequence is $K$-complemented in 
$\N_*$, where $K$ depends only on $\delta$ and $\ep$.}
Indeed, let $W$ denote the closed linear span of the $f_i$'s; let 
$P_1,P_2,\ldots$ be as in the statement of \ref{lem:3.12}, and let 
$g_j = P_j f_j P_j$ for all $j$, then let $Z$ denote the closed linear 
span of the $g_j$'s. 
Of course then $Z$ is {\it isometric\/} to $\ell^1$ (or $\ell_n^1$ if 
the sequence has $n$ terms). 
We may easily define a contractive projection $R:\N_* \to Z$ as follows. 
For each $j$, choose by duality an element $\varphi_j \in \N$ of 
norm one with $\varphi_j = P_j \varphi_j P_j$ and 
\begin{equation}\label{eq:3.67} 
\langle \varphi_j,g_j\rangle = \|g_j\|_1\ .
\end{equation}
(Note that $1\ge \|g_j\|_1 \ge \delta$ for all $j$.) 
Then define 
\begin{equation}\label{eq:3.68} 
R(f) = \sum \langle \varphi_j,f\rangle \|g_j\|_1^{-1} g_j 
\end{equation} 
for $f\in \N_*$. 
Next, define an operator $U:W\to Z$ by 
\begin{equation}\label{eq:3.69} 
U (\sum c_j f_j) = \sum c_j g_j
\end{equation}
for all $c_j$'s with $\sum |c_j| <\infty$. 
Then Lemma~\ref{lem:3.12} yields that $U$ is invertible with 
\begin{equation}\label{eq:3.70} 
\|U^{-1}\| \le (\delta -\ep)^{-1}\ .
\end{equation}
Now a simple computation yields that 
\begin{equation}\label{eq:3.71} 
\|U(w) - R(w)\| \le \frac{\ep}{\delta} \|U(w)\| \ \text{ for all }\ w\in W\ .
\end{equation}
It then follows that $R|W$ is an  isomorphism mapping $W$ onto $Z$, with 
\begin{equation}\label{eq:3.72} 
\|(R|W)^{-1}\| \le \left[ \left(1-\frac{\ep}{\delta}\right) 
(\delta-\ep)\right]^{-1} \defeq K\ .
\end{equation} 
Finally, $Q\defeq (R|W)^{-1} R$ is thus a projection from $\N$ onto $W$, 
with $\|Q\| \le K$. 
It then follows that the elements satisfying the conclusion of 
Lemma~\ref{lem:3.11} have a ``well-complemented'' linear span. 

We also obtain finally, a quantitative proof of Lemma~\ref{lem:3.10}, 
yielding also the result of H.~Pfitzner \cite{Pf} that the preduals of 
von~Neumann algebras have Pe{\l}czy\'nski's property $(V^*)$. 

\begin{tenprime}\label{lem:3.tenprime}
Let $\N$ be an arbitrary von~Neumann algebra, and $W$ be a subset 
of $\Ba \N_*$ so that there exists a sequence $P_1,P_2,\ldots$ of 
orthogonal projections in $\N$ with 
\begin{equation}\label{eq:3.73} 
\varlimsup_j \sup_{w\in W} |\langle P_j,w\rangle| \defeq \eta >0\ .
\end{equation}
Then given $C>\frac1{\eta}$, there exists a sequence 
$w_1,w_2,\ldots$ in $W$ which is $C$-equivalent to the usual $\ell^1$-basis, 
with closed linear span $C$-complemented in $\N_*$. 
\end{tenprime}

\begin{rem}
By Akemann's criterion \cite{A}, it thus follows that any bounded 
non-relatively weakly compact subset of $\N_*$ contains a sequence 
equivalent to the $\ell^1$-basis, with complemented span. 
This is an equivalent formulation of property $(V^*)$.
\end{rem}

\begin{proof} 
It follows easily that we may choose $(f_i)$ a sequence in $W$ and 
$n_1<n_2<\cdots$ so that 
\begin{equation}\label{eq:3.74} 
\varliminf |\langle P_{n_j} ,f_j\rangle| \ge \eta\ .
\end{equation}
Then given $0<\ep<\eta'<\eta$, Lemma~\ref{lem:3.14}A yields 
a subsequence $(f'_j)$ of $(f_j)$ so that $(f'_j)$ is 
$(\eta',\ep)$-relatively disjoint. 
Finally, since $\eta'$ may be arbitrarily close to $\eta$ and $\ep$ 
arbitrarily small, we deduce from Lemma~\ref{lem:3.12} and \eqref{eq:3.72} 
that given $C>\frac1{\eta}$, $(f'_i)$ may be chosen 
$C$-equivalent to the $\ell^1$-basis with span $C$-complemented in $\N_*$.
\end{proof}

\section{Complements on the Banach/operator space structure of 
$L^p(\N)$-spaces} 
\setcounter{equation}{0}

We give here several applications of our main result, and the techniques 
used in its proof. 
For the first one, we let $\Row$ (resp. $\Col$) denote the operator row 
(resp. column) space. 
We also follow the notation in \cite{Pi2}: 
for a given operator space $X$, $X^{\op}$ (the ``opposite'' of $X$) denotes 
the following operator space: if $X\subset B(H)$ and $(x_{ij})$ is an element 
of $\K \otimes_{\op} X$, regarded as a matrix, then $X^{\op} \defeq 
\{(x_{ji}) : (x_{ij}) \in \K\otimes_{\op} X\}$, where $\K$ denotes 
the space of compact opertors on $\ell^2$.
One then has that $\Row^* = \Row^{\op} = \Col$. 

\begin{prop}\label{prop:5.1} 
Let $\N$ be a finite von~Neumann algebra. 
Then neither $\Row$ nor $\Col$ is completely isomorphic to a subspace 
of $L^1(\N)$. 
\end{prop}

\begin{proof}
Suppose to the contrary that there exists an $X\subset L^1(\N)$ with $X$ 
completely isomorphic to $\Row$. 
But then $X^{\op}\subset L^1(\N^{\op})$ is completely isomorphic to $\Col$. 
Let then $\M = \N^{\op}\bar\otimes \N$. 
$\M$ is again a finite von-Neumann algebra, and now $X^{\op} \hat\otimes X$ 
{\it is\/} a subspace of $L^1(\M)$; 
that is, $\Col \hat\otimes \Row$ is completely 
isomorphic to a subspace of $L^1(\M)$. 
But $\Col \hat\otimes\Row$ {\it is\/} (completely isometric to) $C_1$; 
this contradicts our main result.
\end{proof} 

\begin{rem}
An operator space $X$ is called {\it homogeneous\/} if every bounded linear 
operator on $X$ is completely bounded; $X$ is called Hilbertian if it is 
Banach isomorphic to a Hilbert space. 
The above argument then yields the following generalization (since $\Row$ 
is indeed a homogeneous Hilbertian operator space). 
\end{rem}

\begin{Prop}
Let $X$ be an infinite dimensional Hilbertian homogeneous operator space 
so that $X^*$ is completely isomorphic to $X^{\op}$, and let $\N$ be 
a finite von~Neumann algebra. 
Then $X$ is not completely isomorphic to a subspace of $L^1(\N)$. 
\end{Prop}

To obtain this, first observe that the hypotheses yield that 
$X^*\otimes_{\op} X$ is Banach isomorphic to $\K$. 
Hence $X^*\hat\otimes X$ is Banach isomorphic to $C_1$. 
But $X^* \hat\otimes X$ is completely isomorphic to $X^{\op}\hat\otimes X$ 
by hypothesis; as above, if we then assume that $X\subset L^1(\N)$, 
we obtain that $C_1$ Banach embeds in $L^1(\M)$, again contradicting 
our main result.\qed

Our next result yields characterizations of those subspaces of 
$L^p(\N)$ which contain $\ell^p$ isomorphically ($1\le p<2$, $\N$ finite). 
We have need of the following concept. 
(For isomorphic Banach spaces $X$ and $Y$, $d(X,Y) = \inf 
\{\|T\|\, \|T^{-1}\| : T:X\to Y$ is a surjective isomorphism). 

\begin{DEF}\label{def:5.2}
Let $1\le p\le\infty$. 
A Banach space $X$ is said to contain $\ell_n^p$'s if there is a $C\ge1$ 
so that for all $n$, there exists a subspace $X_n$ of $X$ with 
$d(X_n,\ell_n^p)\le C$.
\end{DEF}

A remarkable result of J.L.~Krivine yields that if a Banach space contains 
$\ell_n^p$'s, it contains them almost isometrically (\cite{Kr}; cf.\ 
also \cite{R3}, \cite{L}). 
That is, then for every $\ep$ and $n$, one can choose $X_n\subset X$ 
with $d(X_n,\ell_n^p) <1+\ep$. 
(Of course the famous Dvoretzky theorem yields that every infinite 
dimensional Banach space contains $\ell_n^2$'s almost isometrically; 
also the case $p=1$ or $\infty$ in Krivine's Theorem was established 
previously by Giesy-James \cite{GJ}.) 

We also need the following natural notion. 

\begin{DEF}\label{def:5.3} 
Let $\N$ be a von~Neumann algebra and $1\le p<\infty$. 
A sequence $(g_n)$ in $L^p(\N)$ is called {\it disjointly supported\/} 
provided there exists a sequence $P_1,P_2,\ldots$ of pairwise orthogonal 
projections in $\N$ so that $g_j = P_jg_jP_j$ for all $j$. 
A semi-normalized sequence $(f_n)$ in $L^p(\N)$ is called 
{\it almost disjointly supported\/} if there exists a disjointly 
supported sequence $(g_j)$ in $L^p(\N)$ so that 
$\lim_{n\to \infty} \|f_n-g_n\|_{L^p(\N)} =0$. 
\end{DEF} 

Of course a disjointly supported sequence of non-zero elements spans 
a subspace isometric to $\ell^p$. 
A standard elementary perturbation argument then yields that an almost 
disjointly supported sequence in $L^p(\N)$ has, for every $\ep>0$, a 
subsequence spanning a subspace $(1+\ep)$-isomorphic to $\ell^p$. 
The next result yields in particular that for $\N$ finite, and 
$1\le p<2$, subspaces of $L^p(\N)$ which are isomorphic to $\ell^p$ 
always contain almost disjointly supported sequences. 

\begin{thm}\label{thm:5.4}
Let $1\le p<2$ and $\N$ be a finite von~Neumann algebra; let $\T$ be a 
faithful normal tracial state on $\N$. 
Let $X$ be a closed linear subspace of $L^p(\N)$. 
The following assertions are equivalent. 
\begin{itemize}
\item[1.] $X$ contains a subspace isomorphic to $\ell^p$. 
\item[2.] $X$ contains $\ell_n^p$'s. 
\item[3.] $\{|x|^p :x\in \Ba (X)\}$ is not uniformly integrable. 
\item[4.] $\sup_{f\in\Ba(X)} \om_p (f,\ep) = \sup_{f\in\Ba(X)} \tilde \om_p
(f,\ep) =1$ for all $\ep>0$. 
\item[5.] The $p$ and $1$ norms on $X$ are not equivalent (in case $p>1$).
\item[6.] $X$ contains an almost disjointly supported sequence. 
\item[7.] For all $\ep>0$, $X$ contains a subspace $(1+\ep)$-isomorphic 
to $\ell^p$.
\end{itemize}
\end{thm}

\begin{rems}
1. This result is established for the commutative case in \cite{R2}; the 
case $p>2$ is also valid, and follows (with some extra work for 
assertion 5) from the results in \cite{S1}. 
Again, the commutative case for $p>2$ is immediate from the classical 
work of Kadec-Pe{\l}czy\'nski \cite{KP}. 
Also, condition~5 may be replaced by the following one, valid also 
for $p=1$:
\begin{itemize}
\item[5$'$.] {\it The $p$ and $q$ quasi-norms are not equivalent on $X$ 
for all $0<q<p$.}
\end{itemize}
2. The equivalences of 1, 5, 6 and 7 of Theorem~\ref{thm:5.4} follow also 
from recent work of N.~Randrianantoanina, which establishes these also 
for semi-finite von-Neumann algebras $\N$ and $1\le p<\infty$, $p\ne2$ 
(\cite{Ra1} and \cite{Ra2}). 
\end{rems}

\begin{proof} 
We show $1\implies 2\implies 4\implies 6\implies 7\implies 1$, 
$4\implies 3\implies 2$, and $4\implies 5\implies 3$ in case $p>1$. 
Of course $1\implies 2$ and $7\implies 1$ are trivial. 
So is $4\implies 3$, in virtue of Lemma~\ref{lem:2.3}.

$2\implies 4$.
Fix $\delta >0$. 
Choosing an ``almost isometric'' copy of $\ell_n^p$ in $X$ by 
Krivine's theorem, we shall show that for $n$ large enough, one of 
the natural basis elements $f_i$ of this copy is such that $\tilde\om_p
(f_i,\delta)$ is almost equal to 1.

Define $\lambda$ by 
\begin{equation}\label{eq:5.1} 
\lambda = \sup \{\tilde\om_p (x,\delta) :x\in X,\ \|x\|\le1\}\ .
\end{equation}
Let $C>1$, and using Krivine's theorem, choose $f_1,\ldots,f_n\in \Ba (X)$ 
with $(f_1,\ldots,f_n)$ $C$-equivalent to the $\ell_n^p$ basis. 
In particular, we have that 
\begin{equation}\label{eq:5.2} 
\Big\|\sum_{i=1}^n \pm f_i\Big\|_p 
\ge \frac1C n^{1/p} \text{ for all choices of }
\pm\ .
\end{equation}
Again by the final assertion of Lemma~\ref{lem:2.3}, we may choose for 
each $i$ a $\psi_i \in\N$ so that 
\begin{equation}\label{eq:5.3} 
\|\psi_i \|_\infty \le \delta^{-1/p}\ \text{ and }\ 
\|f_i - \psi_i\|\le\tilde\om_p (f_i,\delta) \le\lambda\ .
\end{equation}
Thus letting $\beta$ be as in the proof of Lemma~\ref{lem:3.1}, 
again we have 
\begin{equation}\label{eq:5.4}
\begin{split}
\frac1C n^{1/p} 
& \le \|\sum f_i\otimes r_i\|_{L^p(\beta)}\ \text{ by \eqref{eq:5.2}}\\
& \le \|\sum \psi_i\otimes r_i\|_{L^2(\beta)} + 
\|\sum (f_i-\psi_i)\otimes r_i\|_{L^p(\beta)}\\
&\le \delta^{-1/p}\sqrt{n} + \lambda n^{1/p} 
\end{split}
\end{equation}
by \eqref{eq:5.3} {\it and\/} the fact that $L^p(\beta)$ is type $p$ 
{\it with constant one\/}.

Thus 
\begin{equation}\label{eq:5.5} 
\frac1C - \frac1{\delta^{1/p} n^{\frac1p-\frac12}} \le\lambda\ .
\end{equation}
Since $C>1$ and $n$ are arbitrary, we obtain that $\lambda=1$, 
proving $2\implies 4$. 

$4\implies 6$.
We first note  that assuming 4, then given $1>\ep>0$, we may choose 
$f\in X$ with $\|f\|_p =1$ and $P\in \P(\N)$ with $\T(P)<\ep$ so that 
\begin{equation}\label{eq:5.6}
\|fP\|_p > 1-\ep\ \text{ and }\ \|f(I-P)\|_p <\ep\ .
\end{equation}
Indeed, choose $f$ in $X$ of norm one so that $\tilde\om_p(f,\ep)>1-\ep$. 
Then choose $P$ a spectral projection for $|f|$ with 
$\|fP\|_p > (1-\ep^p)^{1/p}$. 
But then since $P$ commutes with $|f|$, 
\begin{equation}\label{eq:5.7} 
\|fP\|_p^p = \T (|f|^p P)\ \text{ and }\ \|f(I-P)\|_p^p = \T(|f|^p (I-P))\ ,
\end{equation}
whence 
\begin{align}
1&\ge\T(|f|^p P) +\T (|f|^p (I-P)) \ge (1-\ep) + \|f(I-P)\|_p^p\label{eq:5.8}\\
1&\ge \T(|f|^pP) + \T(|f|^p (I-P)) \label{eq:5.9}\\
&\ge 1-\ep^p + \|f(I-P)\|_p^p\ ,\notag
\end{align}
so $\|f(I-P)\|_p <\ep$ as desired. 
Now since $|f|$ and $|f^*|$ are unitarily equivalent in $\N$, we 
also obtain the existence of a $Q\in P(\N)$ with $\T(Q)<\ep$ so that 
\begin{equation}\label{eq:5.10}
\|Qf\|_p >1-\ep\ \text{ and }\ \|f(I-Q)\|_p <\ep\ .
\end{equation}
Then let $R= P\vee Q$. 
We have 
\begin{equation}\label{eq:5.11}
\T (R)<2\ep\ \text{ and }\ \|f-RfR\| <2\ep\ .
\end{equation}
Indeed, the first estimate is trivial; but 
\begin{equation*}
f-RfR  = f(I-R)+ (I-R)fR 
= f(I-P) (I-R) + (I-R)(I-Q) fR
\end{equation*}
and so \eqref{eq:5.11} follows from \eqref{eq:5.6} and \eqref{eq:5.10}. 

Now  using that for $\ep>0$, 
$f$ of norm~1 in $X$ and $R$ may be chosen satisfying \eqref{eq:5.11} 
we choose inductively $f_1,f_2,\ldots$ in $X$ of norm one, 
$1>\delta_1 >\delta_2 >\cdots >0$, and $Q_1,Q_2,\ldots$ in $\P(\N)$ so 
that for all $j$, 
\begin{gather}
\|f_j - Q_j f_j Q_j\|_p < \frac1{2^j}\ \text{ and }\ \T(Q_j)\le 
\frac{\delta_j}{2^j} 
\label{eq:5.12}\\
\om_p (f_j,\delta_{j+1}) <\frac1{2^j}  \ .
\label{eq:5.13}
\end{gather}
To see this is possible, just choose $\delta_1=1/2$, then choose $f_1$ and 
$Q_1$ thanks to \eqref{eq:5.11}. 
Suppose $f_1,\ldots,f_n$, and $\delta_n$ chosen. 
By uniform integrability of $\{|f_n|^p\}$, choose 
$\delta_{n+1}<\delta_n$ so that $\om_p (f_n,\delta_{n+1})< 1/2^{n+1}$.  
Then choose $f_{n+1}$ and $Q_{n+1}$ satisfying \eqref{eq:5.12} for $j=n+1$.

Now define projections $P_j$ and $\tilde Q_j$ by \eqref{eq:3.19}.  
The $P_j$'s are orthogonal and 
by the argument for the last part of Proposition~\ref{prop:2.2}, fixing $j$, 
we have 
\begin{equation}\label{eq:5.14} 
\begin{split}
\T(\tilde Q_j) \le \sum_{k>j} \T(Q_k) 
& \le \delta_{j+1} \sum_{k>j} \frac1{2^k}\ \text{ by \eqref{eq:5.12}}\\
& < \delta_{j+1}\ .
\end{split}
\end{equation}
Hence 
\begin{equation*}
\|\tilde Q_j f_j \|_p \le \om_p (f_j^*,\delta_{j+1}) 
= \om_p (f_j,\delta_{j+1}) <\frac1{2^j}
\end{equation*}
(by \eqref{eq:5.13}) and also 
\begin{equation*}
\|f_j\tilde Q_j\|_p \le \om_p (f_j,\delta_{j+1}) <\frac1{2^j}\ .
\end{equation*}
Hence 
\begin{equation}\label{eq:5.15} 
\|\tilde Q_j f_j Q_j \|_p < \frac1{2^j}\ \text{ and }\ 
\|Q_j f_j \tilde Q_j\|_p < \frac1{2^j}\ .
\end{equation}
Hence finally we have by \eqref{eq:5.12} and \eqref{eq:5.15}, 
\begin{equation}\label{eq:5.16} 
\|f_j -P_jf_jP_j\| \le \frac3{2^j}\ \text{ for all }\ j\ .
\end{equation}
Thus $(f_j)$ is almost disjointly supported, proving that 6 holds. 

$6\implies 7$ 
is a standard perturbation argument in Banach space theory. 
Assuming 6 holds, we may choose a normalized disjointly supported 
sequence $(g_n)$ in $L^p(\N)$ and a sequence $(f_n)$ in $X$ so that 
\begin{equation}\label{eq:5.17} 
\sum \|g_n-f_n\|_p <\infty\ .
\end{equation}
But now $(g_n)$ is 1-equivalent to the $\ell^p$-basis, and a simple 
perturbation argument gives that given $\ep>0$, there is an $N$ so 
that $(f_n)_{n\ge N}$ is $(1+\ep)$-equivalent to the $\ell^p$ basis. 
(Thus $(f_n)$ is ``almost isometrically equivalent'' to the 
$\ell^p$ basis.) 

$3\implies 2$. 
We have that if $p=1$, $X$ contains a subspace isomorphic to $\ell^1$ 
by Lemma~\ref{lem:3.10}, so assume $p>1$. 
We may choose a sequence $(f_n)$ of norm-1 elements of $X$, 
$\delta_1 >\delta_2>\cdots$ with $\delta_n\to0$ and $\eta >0$ so that 
\begin{equation}\label{eq:5.18} 
\om_p (f_n,\delta_n) >\eta\ \text{ for all }\ n\ .
\end{equation}

By passing to a subsequence, we may assume without loss of generality 
that $(f_n)$ is weakly convergent, with weak limit $f$, say. 
But 
\begin{equation}\label{eq:5.19}
\om_p (f_n-f,\delta_n) \ge \om_p (f_n,\delta_n) - \om_p (f,\delta_n)
\end{equation}
and hence
\begin{equation}\label{eq:5.20} 
\varliminf_{n\to\infty} \om_p (f_n-f,\delta_n) \ge\eta\ .
\end{equation}
That is, we have now obtained a weakly null sequence $(g_n)$ in $X$ 
so that 
\begin{equation}\label{eq:5.21}
(|g_n|^p)\ \text{ is not uniformly integrable.}
\end{equation}
By Corollary~\ref{cor:3.4}, after passing to a subsequence of $(g_n)$, 
we may assume 
\begin{equation}\label{eq:5.22} 
(g_n\otimes r_n)\text{ is $C$-equivalent to the usual 
$\ell^p$-basis in $L^p(\beta)$ for some $C$.}
\end{equation}

Now Lemma~\ref{lem:3.7} yields that for all $n$, there exist 
$m_1<m_2<\cdots < m_n$ so that $g_{m_1},\ldots,g_{m_n}$ is 
4-unconditional, and hence 
\begin{equation}\label{eq:5.23} 
(g_{m_i})_{i=1}^n \text{ is $4C$-equivalent to the $\ell_n^p$-basis.}
\end{equation}

This proves that 2 holds. 
Now assume $p>1$. 

$4\implies 5$. 
Let $\ep>0$ and choose $f\in X$ with $\|f\|_p =1$ and 
$P\in \P (\N)$ with $\T(P) <\ep$ so that \eqref{eq:5.6} holds. 
Then of course 
\begin{equation}\label{eq:5.24} 
\|f(I-P)\|_1 <\ep\ .
\end{equation}
Now letting $\frac1p + \frac1q =1$, 
\begin{equation}\label{eq:5.25} 
\|fP\|_1 \le \|f\|_p \|P\|_q \le \ep^{1/q}\ \text{ by H\"older's 
inequality.}
\end{equation}
Thus 
\begin{equation}\label{eq:5.26} 
\|f\|_1 < \ep + \ep^{1/q}\ .
\end{equation}
Since $\|f\|_p =1$ and $\ep>0$ is arbitrary, 5 holds.

$5\implies 3$.
Suppose 5 holds, yet 3 were false. 
Choose $0<\delta$ so that 
\begin{equation}\label{eq:5.27} 
\tilde\om_p (f,\delta) \le \frac12\ \text{ for all }\ f\in\Ba (X)\ .
\end{equation}
Let $f\in X$, $\|f\|_p =1$. 
By the last statement of Lemma~\ref{lem:2.3}, choose $P$ a spectral 
projection for $|f|$ so that $fP\in\N$ with 
\begin{equation}\label{eq:5.28} 
\|f(I-P)\|_p \le \frac12\ \text{ and }\ 
\|fP\|_\infty \le \delta^{-1/p}\ .
\end{equation} 
Then 
\begin{equation}\label{eq:5.29}
\begin{split}
\frac1{2^p} \le \|fP\|_p^p 
& = \T(|f|^pP)\ \text{ (since $P\leftrightarrow |f|$)}\\
& = \T(|f| \, |f|^{p-1}P)\\
& \le \|f\|_1 \delta^{1-\frac1p}\ .
\end{split}
\end{equation}
That is, 
\begin{equation}\label{eq:5.30}
\|f\|_1 \ge 2^{-1/p} \delta^{\frac1p-1} \defeq C\ .
\end{equation}
\eqref{eq:5.30} yields that $\|g\|_p \le C\|g\|_1$ for all $g\in  X$; 
i.e., 5 does not hold, a contradiction. 
This completes the proof of the theorem.
\end{proof}

The final result of this section deals with the Banach-Saks property. 

\begin{DEF}\label{def:5.5}
Let $X$ be a Banach space, and $1<p<\infty$. 

{\rm (a)} 
Let $(x_n)$ be a weakly null sequence in $X$. 
$(x_n)$ is called 
\begin{itemize}
\item[(i)] a Banach-Saks sequence if 
\begin{equation}\label{eq:5.31} 
\lim_{n\to\infty} n^{-1}\Big\|\sum_{j=1}^n y_j\Big\|=0
\text{ for all subsequences } (y_j)\text{ of } (x_j)\ .
\end{equation}
\item[(ii)] a $p$-Banach-Saks sequence if 
\begin{equation}\label{eq:5.32} 
\text{ there is a $C<\infty$ so that }
\varlimsup_{n\to\infty} n^{-1/p} \Big\|\sum_{j=1}^n y_j\Big\|\le C
\text{ for all subsequences $(y_j)$ of $(x_j)$.}
\end{equation}
\item[(iii)] a strong $p$-Banach-Saks sequence if 
\begin{equation}\label{eq:5.33}
\lim_{n\to\infty} n^{-1/p} \Big\|\sum_{j=1}^n y_j\Big\|=0
\text{ for all subsequences $(y_j)$ of $(x_j)$.}
\end{equation}
\end{itemize}

{\rm (b)}
$X$ is said to have the Banach-Saks property (resp. the $p$-Banach-Saks 
property) (resp. the strong $p$-Banach-Saks property) 
if every weakly null sequence in $X$ has a Banach-Saks 
(resp. $p$-Banach-Saks) (resp. strong $p$-Banach-Saks) subsequence. 
\end{DEF}

The classical paper of Banach-Saks \cite{BS} yields that commutative $L^p$ 
spaces have the $p$-Banach-Saks property, for $1<p\le2$; the fact that 
$L^1$-spaces have the Banach-Saks property was proved later by 
Szlenk \cite{Sz}. 
Our last result yields in particular a generalization to the spaces $L^p(\N)$, 
$\N$ finite. 
Most of its assertions follow very quickly from our previous results. 

\begin{prop}\label{prop:5.6} 
Let $\N$ be a finite von-Neumann algebra and $1<p<2$. 
\begin{itemize}
\item[1.] $L^1(\N)$ has the Banach-Saks property and $L^p(\N)$ has 
the $p$-Banach-Saks property.
\item[2.] A weakly null sequence $(f_n)$ in $L^p(\N)$ has a strong 
$p$-Banach-Saks subsequence if $(|f_n|^p)$ is uniformly integrable. 
If $(|f_n|^p)$ is not uniformly integrable, $(f_n)$ has a subsequence 
$(f'_n)$ so that for some $c>0$ and all subsequences $(y_j)$ of 
$(f'_j)$, 
\begin{equation}\label{eq:5.34}
\varliminf n^{-1/p} \Big\|\sum_{j=1}^n y_j\Big\| \ge c\ .
\end{equation}
\item[3.] A closed linear subspace $X$ of $L^p(\N)$ has the strong
$p$-Banach-Saks property if and only if $X$ has no subspace 
isomorphic to $\ell^p$.
\end{itemize}
\end{prop}

\begin{proof} 
Corollary~\ref{cor:3.8} together with Proposition~\ref{prop:2.2} 
yields that $L^1(\N)$ has the Banach-Saks property. 
It also yields the first assertion in 2. 
Suppose $(|f_n|^p)$ is not uniformly integrable and assume (without 
loss of generality) that $\|f_n\|_p \le 1$ for all $n$. 
Applying Corollary~\ref{cor:3.4} and Lemma~\ref{lem:3.7}, we may 
choose a subsequence $(f'_n)$ of $(f_n)$ so that for some $C\ge1$, 
\begin{equation}\label{eq:5.35}
(f'_n\otimes r_n)\text{ is $C$-equivalent to the usual $\ell^p$-basis.}
\end{equation}
and
\begin{equation}\label{eq:5.36} 
(f'_{n_1},\ldots,f'_{n_{2^k}})\text{ is 4-unconditional for all } 
k\le n_1<n_2<\cdots < n_{2^k}\ .
\end{equation}
Suppose $(y_j)$ is a subsequence of $(f'_j)$. 
Then it follows that for all $k$, 
\begin{equation}\label{eq:5.37}
(y_{k+1},\ldots,y_{k+2^k})\text{ is $(4C)$-equivalent to the 
$\ell_{2^k}^p$-basis.}
\end{equation}
Let $n$ be a ``large'' integer and choose $k$ with 
\begin{equation}\label{eq:5.38}
2^{k-1} \le n< 2^k\ .
\end{equation}
Then 
\begin{equation}\label{eq:5.39}
\Big\|\sum_{j=k+1}^n y_j\Big\| \ge \frac{(n-k)^{1/p}}{4C}\ \text{ by 
\eqref{eq:5.37}}\ .
\end{equation}
Thus 
\begin{equation}\label{eq:5.40}
\Big\| \sum_{j=1}^n y_j \Big\|_p 
\ge \frac{(n-k)^{1/p}}{4C} - k 
\ge \frac{(n-\log_2 n-1)^{1/p}}{4C} - \log_2 n -1\ .
\end{equation}
Hence 
\begin{equation}\label{eq:5.41}
\varliminf_{n\to\infty} n^{-1/p} \Big\|\sum_{j=1}^n y_j\Big\|_p \ge 
\frac1{4C}\ .
\end{equation}
This completes the proof of assertion 2 of the Proposition. 
But we also have that 
\begin{equation}\label{eq:5.42} 
\Big\|\sum_{j=k+1}^n y_j\Big\|_p \le 4C(n-k)^{1/p}\ 
\text{ by \eqref{eq:5.37}}, 
\end{equation} 
and so 
\begin{equation}\label{eq:5.43} 
\Big\|\sum_{j=1}^n y_j\Big\|_p \le 4C (n-\log_2 n)^{1/p} 
+ \log_2 n +1\ ,
\end{equation}
thus 
\begin{equation}\label{eq:5.44} 
\varlimsup_{n\to\infty} n^{-1/p} \Big\|\sum_{j=1}^n y_j\Big\|_p 
\le 4C\ .
\end{equation}

This proves that $L^p(\N)$ has the $p$-Banach-Saks property, for any 
weakly null sequence $(f_n)$ in $L^p(\N)$ either has $(|f_n|^p)$ 
uniformly integrable (and hence a strong $p$-Banach-Saks subsequence), 
or a subsequence $(f'_n)$ as above. 

The final assertion of the Proposition follows immediately from 
Theorem~\ref{thm:5.4} and assertion 2.
\end{proof}

\begin{rem}
Of course Hilbert space has the 2-Banach Saks property. 
Actually, it can be shown that $L^p(\N)$ 
has the 2-Banach Saks property for $p>2$ and $\N$ finite, and 
this is best possible (in general). 
Indeed, if $(f_j)$ is a weakly null sequence in $L^p(\N)$, then if 
$\|f_j\|_p\to 0$, $(f_j)$ trivially has a $p$-Banach Saks subsequence; 
the same is true if $(f_j)$ has a subsequence  equivalent to 
the $\ell^p$-basis  (and of course a $p$-Banach Saks sequence 
is a 2-Banach Saks sequence).
Otherwise, combining arguments in [S1]
Theorem 2.4 with the arguments in Proposition 5.6, we see that there
exists a subsequence $(f'_j)$ of $(f_j)$ such that its all subsequences
$(y_n)$ are 2-Banach Saks.
\end{rem}

We conclude this section with a brief discussion of the following open 

\begin{Problem}
Let $1<p<2$ and $(f_n)$ be a seminormalized weakly null sequence in 
$L^p(\N)$ ($\N$ a finite von~Neumann algebra) such that 
$(|f_n|^p)$ is not uniformly integrable. 
Does $(f_n)$ have a subsequence equivalent to the usual $\ell^p$ basis?
\end{Problem}

As pointed out previously, the answer is affirmative if $(f_n)$ has an 
unconditional subsequence. 
Actually, it can be proved that if $(f_n)$ 
{\it satisfies the hypotheses of this Problem, it 
has a subsequence $(f'_n)$ which dominates the $\ell^p$-basis
and moreover has spreading model equivalent to the $\ell^p$-basis\/}.
(The last assertion follows immediately from our proof of 
Proposition~\ref{prop:5.6}.) 
It may then be shown that the above Problem is equivalent to the 
following one (in which the hypothesis concerning $(|f_n|^p)$ no longer 
enters). 

\begin{Problemprime}
Let $(f_n)$ be a seminormalized basic sequence in $L^p(\N)$, $p$ 
and $\N$ as above. 
Does $(f_n)$ have a subsequence $(f'_n)$ which is dominated by 
the $\ell^p$-basis?
i.e., such that $\sum c_j f'_j$ converges in $L^p(\N)$ whenever 
$\sum |c_j|^p<\infty$?
\end{Problemprime}

\section{The Banach isomorphic classification of the spaces $L^p(\N)$ for 
$\N$ hyperfinite semi-finite}		
\setcounter{equation}{0}

We first fix some notation. 
Let $1\le p<\infty$. 
We let $S_p = (\bigoplus_{n=1}^\infty C_p^n)_p$ 
($=L^p(\oplus M_n)_\infty$). 
To avoid confusion, we denote by $L_p\otimes_pX$ the Bochner space 
$L_p(X,m)$, where $m$ is Lebesgue measure and $X$ is a Banach space. 
Thus e.g., $L_p\otimes_p C_p = L_p (C_p) = L^p (L^\infty (m) 
\bar\otimes B(\ell^2))$. 
$\R$ denotes the hyperfinite type~II factor, and $L^p(\R)\otimes_pC_p$ 
denotes $L^p (\R\bar\otimes B(\ell^2))$ (so $\R \bar\otimes B(\ell^2)$ 
is the hyperfinite type~II$_\infty$ factor). 

The main motivating result of this section is as follows. 

\begin{thm}\label{thm:6.1} 
Let $\N$ be a hyperfinite semi-finite infinite dimensional von-Neumann 
algebra, and let $1\le p<\infty$, $p\ne 2$. 
Then $L^p(\N)$ is (completely) isomorphic to precisely one of the 
following thirteen Banach spaces. 
\begin{gather*}
\ell_p\ ,\quad S_p\ ,\quad L_p\ ,\quad C_p\ ,\quad 
S_p\oplus L_p\ ,\quad 
C_p \oplus L_p\ ,\quad 
L_p\otimes_p S_p\ ,\quad
C_p\oplus (L_p\otimes_p S_p)\\
L^p (\R)\ ,\quad 
L_p \otimes_p C_p\ ,\quad 
C_p \oplus L^p (\R)\ ,\quad 
L^p (\R) \oplus (L_p\otimes_p C_p)\ ,\quad 
L^p (\R) \otimes_p C_p\ .
\end{gather*}
\end{thm}

Theorem~\ref{thm:6.1} is an immediate consequence of the following 
finer result concerning embeddings. 

\begin{thm}\label{thm:6.2} 
Let $1\le p<2$. 
If $\N$ is as in \ref{thm:6.1}, then $L^p(\N)$ is (completely) isomorphic 
to one of the thirteen spaces in the tree in Figure~1.	
If $X\ne Y$ are listed in the tree, then $X$ is Banach isomorphic to a 
subspace of $Y$ if and only if $X$ can be joined to $Y$ through a descending 
branch (in which case $X$ is completely isometric to a subspace of $Y$). 
\end{thm}

\begin{figure}
\begin{picture}(300,305)
\put(0,300){1.}
\put(160,300){\line(3,-2){58}}
\put(160,300){\line(-3,-2){115}}
\put(105,261){\line(3,-2){170}}
\put(43,225){\line(3,-2){175}}
\put(217,260){\line(-3,-2){115}}
\put(155,310){$\ell_p$}
\put(156,297){\LARGE$\bullet$}
\put(0,260){2.}
\put(82,262){$S_p$}
\put(101,259){\LARGE$\bullet$}
\put(213,258){\LARGE$\bullet$}
\put(222,260){$L_p$}
\put(0,220){3.}
\put(30,220){$C_p$}
\put(41,219){\LARGE$\bullet$}
\put(157.5,219){\LARGE$\bullet$}
\put(170,220){$S_p\oplus L_p$}
\put(0,180){4.}
\put(217,182){\line(-3,-2){115}}
\put(75,170){$C_p\oplus L_p$}
\put(100,181){\LARGE$\bullet$}
\put(215,181){\LARGE$\bullet$}
\put(222,185){$L_p\otimes_p S_p$}
\put(0,140){5.}
\put(275,145){\line(-3,-2){115}}
\put(159,143){\LARGE$\bullet$}
\put(170,142){$C_p\oplus L_p\otimes_p S_p$}
\put(273,143){\LARGE$\bullet$}
\put(284,145){$L^p(\R)$}
\put(0,100){6.}
\put(101,106){\line(3,-2){58}}
\put(70,90){$L_p\otimes_p C_p$}
\put(99,102){\LARGE$\bullet$}
\put(215,104){\LARGE$\bullet$}
\put(230,103){$C_p\oplus L^p(\R)$}
\put(0,60){7.}
\put(156,63){\LARGE$\bullet$}
\put(165,63){$L^p (\R)\oplus L_p\otimes_p C_p$}
\put(0,20){8.}
\put(159,64){\line(0,-2){40}}
\put(156,20){\LARGE$\bullet$}
\put(130,10){$L^p(\R)\otimes_p C_p$}
\end{picture}
\caption{}
\label{fig:tree}
\end{figure}

\begin{rem}
In the language of graph theory, Theorem~\ref{thm:6.2} asserts that the 
tree in Figure~1 is the {\it Hasse diagram\/} for the partially ordered set 
consisting of the equivalence classes of $L^p(\N)$ under Banach 
isomorphism (over $\N$ as in \ref{thm:6.1}), with the order relation: 
$[X]\le [Y]$ provided $X$ is isomorphic to a subspace of $Y$.
\end{rem}

Parts of Theorem~\ref{thm:6.2} require previously known results, some 
of which are very recent. 
It is established in \cite{S2} that the first nine spaces in the list 
in Theorem~\ref{thm:6.1} are isomorphically distinct when $p=1$, and 
exhaust the list of the possible Banach isomorphism types of $L^p(\N)$ 
for $\N$ type~I ($\N$ as in \ref{thm:6.1}), $p\ne 2$. 

Theorem~\ref{thm:6.2} yields the new result in the type~I case: 
$L_p\otimes_p C_p$ does not embed in $C_p\oplus (L_p\otimes_p S_p)$ 
for $1\le p<2$; 
(another new result in this case, that $C_p$ does not embed in 
$L_p\otimes_p S_p$, follows immediately from Corollary~1.2);
the other embedding results stated in \ref{thm:6.2} 
for the type~I case are given in \cite{S2}. 
We give here a new proof of the particular case that $L_p\otimes_p S_p$ 
does not embed in $L_p\oplus C_p$, using the Main Result of this paper.

We first proceed with the non-embedding results required for 
Theorem~\ref{thm:6.2}. 
The following theorem is crucial. 

\begin{thm}\label{thm:6.3} 
Let $\N$ be a finite von~Neumann algebra and $1\le p<2$. 
Then $L_p\otimes_p C_p$ is not isomorphic to a subspace of 
$C_p\oplus L^p(\N)$. 
\end{thm}

We now fix $1\le p<2$ for the remainder of this section. 

We first require

\begin{lem}\label{lem:6.4} 
Let $T:L_p \to C_p$ be a given bounded linear operator, and let $\ep>0$. 
Then there exists an $f\in L_p$ with $f$ $\{1,-1\}$-valued so that 
$\|Tf\|<\ep$. 
\end{lem}

\begin{sublemma}
The conclusion of \ref{lem:6.4} holds, replacing $C_p$ by $\ell^2$ 
in its hypotheses.
\end{sublemma}

\begin{proof} 
Fix $n$ a positive integer. 
Using the generalized parallelogram identity, 
\begin{equation}\label{eq:6.1}
\begin{split}
av_\pm  \Big\| T\sum_{j=1}^n \pm \chix_{[\frac{j-1}n,\frac{j}n)}\Big\|^2_2 
& = \sum_{j=1}^n \|T (\chix_{[\frac{j-1}n, \frac{j}n)})\|^2_2\\
\noalign{\vskip6pt}
& \le \|T\|^2 \sum_{j=1}^n \|\chix_{[\frac{j-1}n,\frac{j}n)} \|_p^2 \\
\noalign{\vskip6pt}
& = \|T\|^2 \frac{n}{n^{2/p}} 
= \|T\|^2 \frac1{n^{2/p-1}}\ .
\end{split}
\end{equation}
It follows that we may choose $\eta_j = \pm1$ for all $j$ so that 
\begin{equation}\label{eq:6.2} 
\Big\|T\biggl(\sum_{j=1}^n \eta_j\chix_{[\frac{j-1}n,\frac{j}n)}\biggr)\Big\|_2
\le \frac{\|T\|}{n^{\frac1p-\frac12}}\ .
\end{equation}
Now simply choose $n$ so that $\frac{\|T\|}{n^{\frac1p-\frac12}} <\ep$ 
and let $f = \sum_{j=1}^n \eta_j \chix_{[\frac{j-1}n,\frac{j}n)}$.
\end{proof}

\begin{proof}[Proof of Lemma~\ref{lem:6.4}] 
Let $(e_{ij})$ be the matrix units basis for $C_p$, and define for 
each $n$, 
\begin{equation}\label{eq:6.3} 
H_n = [e_{ij}: 1\le i\le n\text{ and }1\le j<\infty\text{ or } 
1\le i<\infty\text{ and } 1\le j\le n]\ .
\end{equation}
Let $P_n$ be the natural basis projection onto $H_n$; i.e., 
$P_n :C_p\to C_p$ is the projection with $P_n (e_{ij})=0$ 
if $e_{ij}\notin H_n$, 
$P_n (e_{ij}) = e_{ij}$ if $e_{ij} \in H_n$ 
(so $\|P_n\|\le 2$). 
Then $H_n$ is isomorphic to $\ell^2$, so by the sub-lemma we may choose 
$f_n$ in $L^p$ with $f_n$ $\{1,-1\}$-valued and 
\begin{equation}\label{eq:6.4}
\|P_nTf_n\| \le \frac1{2^n}\ .
\end{equation}
We claim that 
\begin{equation}\label{eq:6.5} 
\lim_{n\to\infty} \|Tf_n\| =0\ .
\end{equation}
Of course \eqref{eq:6.5} yields the conclusion of the Lemma. 
Suppose \eqref{eq:6.5} were false. 
It follows that $(f_n)$ has a subsequence $(f'_n)$ so that 
\begin{equation}\label{eq:6.6} 
(Tf'_n)\text{ is equivalent to the usual $\ell^p$-basis}
\end{equation}
and 
\begin{equation}\label{eq:6.7} 
(f'_n)\text{ converges weakly in } L^2\ .
\end{equation}
(\eqref{eq:6.6} follows because $(f'_n)$ may be chosen to be a small 
perturbation of a ``block-off-diagonal sequence'', by \ref{lem:6.4}). 

Of course $(f'_n)$ converges weakly in $L^p$ as well, hence $(Tf'_n)$ 
also converges weakly, a contradiction when $p=1$ since then $(Tf'_n)$ is 
equivalent to the $\ell^1$-basis. 

When $p>1$, letting $f$ be the weak limit of $(f_n)$, we have that 
$Tf=0$ since $Tf'_n\to0$ weakly. 
Moreover $\|f\|_\infty\le 2$, so letting $f''_n = f'_n-f$ for all $n$, 
$(f''_n)$ is a uniformly bounded weakly null sequence in $L^p$ with 
$(Tf''_n) = (Tf'_n)$ equivalent to the $\ell^p$-basis. 
Finally, since $(f''_n)$ is also semi-normalized in $L^p$, $(f''_n)$ has 
a subsequence $(g_n)$ equivalent to the usual $\ell^2$-basis. 
(Indeed, we may choose $(g_n)$ equivalent to the $\ell^2$-basis in 
$L^2$-norm, and unconditional. 
But then since $L^p$ has cotype~2, $(g_n)$ is equivalent to the 
$\ell^2$-basis in the $L^p$-norm). 
Still, $(Tg_n)$ is equivalent to the $\ell^p$-basis; this is impossible  
since $p<2$.
\end{proof}

We now apply our Main Result and Lemma~\ref{lem:6.4}, to give the 

\begin{proof}[Proof of Theorem~\ref{thm:6.3}]
Suppose to the contrary that $\N$ is a finite von~Neumann algebra and 
$T:L_p\otimes_p C_p\to C_p\oplus L^p(\N)$ is an isomorphic embedding. 
Of course we may assume that $\|T\|=1$; let $\ep  = \|T^{-1}\|^{-1}$. 
Thus we have 
\begin{equation}\label{eq:6.8}
\|Tf\| \ge \ep \|f\|\ \text{ for all }\ f\in L_p\otimes_p C_p\ .
\end{equation}
Let $P$ be the projection of $C_p \oplus L^p(\N)$ onto $C_p$ with kernel 
$L^p(\N)$, and set $Q= I-P$. 
Also, for each $i$ and $j$, let $Q_{ij}$ be the natural projection 
of $L_p\otimes_p C_p$ onto the space 
\begin{equation}\label{eq:6.9} 
E_{ij} \defeq \{f\otimes e_{ij} :f\in L_p\}\ .
\end{equation}
(As before, $e_{ij}$ denotes the $i,j^{\text{th}}$ matrix unit for $C_p$.
Visualizing $C_p$ as matrices of scalars and 
$L_p\otimes_p C_p$ as all matrices $(f_{ij})$ of functions 
in $L_p$ with 
\begin{equation*}
\|(f_{ij})\| =  \biggl( \int \|(f_{ij}(w))\|_{C_p}^p\,dw \biggr)^{1/p} 
<\infty\ ,
\end{equation*}
then $Q_{ij} ((f_{k\ell})) = f_{ij}\otimes e_{ij}$. 
$E_{ij}$ is just the space of matrices with all entries zero except 
in the $ij^{\text{th}}$ slot). 
Now fix $i$ and $j$. 
Of course $E_{ij}$ is isometric to $L_p$. 

Thus by Lemma~\ref{lem:6.4}, we may choose $f_{ij} \in L_p$ with 
$f_{ij}$ $\{1,-1\}$-valued so that 
\begin{equation}\label{eq:6.10} 
\|PT f_{ij} \otimes e_{ij}\| < \frac{\ep}{2^{i+j+2}}\ .
\end{equation}
Now letting $X = [f_{ij} \otimes e_{ij} :i,j=1,2,\ldots]$, then $X$ 
is a 1-$GC_p$ space, in the terminology of the Introduction. 
That is, every row and column of $(f_{ij}\otimes e_{ij})$ is 1-equivalent 
to the  $\ell^2$ basis, while every generalized  diagonal is 1-equivalent 
to the $\ell^p$ basis. 
Hence $X$ is not isomorphic to a subspace of $L^p(\N)$ by our Main 
Theorem (i.e. Corollary~\ref{cor:1.2}). 
However 
\begin{equation}\label{eq:6.11} 
QT|X\ \text{ is an isomorphic embedding.}
\end{equation}
Indeed, if $x= \sum c_{ij} (f_{ij}\otimes e_{ij})$ with only finitely 
many $c_{ij}$'s non zero, and $\|x\| =1$, then $|c_{ij}|\le1$ for all $i$ 
and $j$ (since the $Q_{ij}$'s are contractive and $\|f_{ij}\|=1$ 
for all $i$ and $j$), and so 
\begin{equation}\label{eq:6.12}
\begin{split} 
\|PTx\| & \le \max_{i,j} |c_{ij}| \sum_{i,j} \|T(f_{ij}\otimes e_{ij})\|\\
&\le \sum_{i=1}^\infty \sum_{j=1}^\infty \frac{\ep}{2^{i+j+2}} 
= \frac{\ep}2
\end{split}
\end{equation}
using  \eqref{eq:6.10} and our assumption that $T$ is a contraction.
Hence 
\begin{equation}\label{eq:6.13} 
\|QTx\| \ge \frac{\ep}2\ \text{ by \eqref{eq:6.8}.}
\end{equation}
This proves \eqref{eq:6.11}, and completes the proof by contradiction.
\end{proof}

Our localization result, Corollary~\ref{cor:1.4}, and the preceding 
proof, yield  an alternate proof of the following result, 
obtained in \cite{S2}.

\begin{prop}\label{prop:6.5} 
$L^p\otimes_p S_p$ is not isomorphic to a subspace of 
$C_p \oplus L_p$.
\end{prop}

\begin{proof} 
We have that $L^p \otimes_p S_p$ is (linearly isometric to) 
$(\bigoplus_{n=1}^\infty L_p \otimes_p C_p^n)_p$. 
Thus it suffices to prove that 
\begin{equation}\label{eq:6.14} 
\lim_{n\to\infty} \lambda_n = \infty 
\end{equation} 
where 
\begin{equation} \label{eq:6.15} 
\lambda_n = \inf \{d(L_p \otimes_p C_p^n,Y) : Y\text{ is a subspace of } 
C_p\oplus L_p\}
\end{equation}
and ``$d$'' denotes the Banach Mazur distance-coefficient (defined
just preceding Corollary~\ref{cor:1.4}). 

Now fix $n$, and let $T:L_p\otimes_p C_p^n \to Y\subset C_p\oplus L_p$ 
be an isomorphic embedding onto $Y$, with 
\begin{equation}\label{eq:6.16} 
\|T\| =1\ \text{ and }\ \|T^{-1}\| \le 2\lambda_n \ .
\end{equation}

Using the notation and reasoning in the proof of Theorem~\ref{thm:6.3}, and 
setting $\ep = 1/(2\lambda_n)$, we may choose for each $i$ and $j$ with 
$1\le i,j\le n$, a $\{1,-1\}$-valued $f_{ij}\in L^p$ 
satisfying \eqref{eq:6.10}. 
We thus obtain that $\|PT|X\| \le \ep/2$ by \eqref{eq:6.12}. 
Hence for all $x\in X$, 
\begin{equation}\label{eq:6.17} 
\|QT(x)\| \ge \left( \frac1{2\lambda_n} - \frac{\ep}2\right)\|x\|
= \frac1{4\lambda_n} \|x\|
\end{equation} 
using also \eqref{eq:6.16}. 
That is, setting $Z= QT(X)$, we have that 
\begin{equation}\label{eq:6.18}
d(X,Z) \le 4\lambda_n\ .
\end{equation}
Now $X$ is a 1-$GC_p^n$-space; thus 
\begin{equation}\label{eq:6.19}
4\lambda_n \ge \beta_{n,1}\ \text{ for all }\ n
\end{equation}
(in the notation of Corollary~\ref{cor:1.4}), so \eqref{eq:6.14} 
holds by Corollary~\ref{cor:1.4}.
\end{proof}

We also require the following rather deep result, due to M.~Junge \cite{J}. 

\begin{thm}\label{thm:6.6} 
$C_q$ is isomorphic to a subspace of $L^p(\R)$ for all $p<q<2$.
\end{thm}

Finally, we require the following (unpublished) result, due to 
G.~Pisier and Q.~Xu \cite{PX2}. 

\begin{lem}\label{lem:6.7} 
Let $X$ be a (closed linear) subspace of $L_p\otimes_p C_p$. 
Then either $X$ embeds in $L_p$ or $\ell^p$ embeds in $X$.
\end{lem}

For the sake of completeness, we sketch a proof. 
First, we give an important, quick consequence of these last two results. 

\begin{cor}\label{cor:6.8} 
$L^p(\R)$ is not isomorphic to a subspace of $L_p\otimes_p C_p$. 
\end{cor}

\begin{proof} 
By Theorem~\ref{thm:6.6}, it suffices to prove that $C_q$ does not embed 
in $L_p\otimes_p C_p$ if $p<q<2$. 
If $C_q$ did embed, then since it does not embed in $L_p$, it would 
have a subspace isomorphic to $\ell^p$, 
by Lemma~\ref{lem:6.7}. 
However it is a standard fact that every infinite-dimensional subspace of 
$C_p$ is either isomorphic to $\ell^2$ or contains a subspace isomorphic 
to $\ell^p$, a contradiction.\qed

We next sketch the proof of Lemma~\ref{lem:6.7} (which also yields the 
above mentioned standard fact). 

Let $(x_{ij})$ be a given matrix in a linear space $X$. 
Call a sequence $(f_k)$ in $X$ a {\it generalized block diagonal\/} 
of $(x_{ij})$ if there exist $i_1< i_2<\cdots$ and $j_1<j_2<\cdots$ so 
that for all $k$, 
\begin{equation}\label{eq:6.20} 
f_k \in [x_{ij} : i_k \le i<i_{k+1}\text{ and } j_k \le j< j_{k+1}]\ .
\end{equation}
Now if $(f_k)$ is a generalized block diagonal of the matrix 
$(e_{ij})$ consisting 
of non-zero terms, $e_{ij}$ the matrix units for $C_p$ (as above), then 
$(f_k/\|f_k\|)$ is isometrically equivalent to the $\ell^p$-basis. 
But then it follows immediately that if $(f_k)$ is a normalized generalized 
block diagonal of $(\bone \otimes e_{ij})$ (in $L^p\otimes_p C_p$) 
consisting of non-zero terms, $(f_k)$ is also isometrically equivalent 
to the $\ell^p$-basis. 
Indeed, for any scalars $c_1,c_2,\ldots$ with only finitely many non-zero 
terms, and any $0\le c_j\le1$, 
\begin{equation}\label{eq:6.21} 
\|\sum c_j f_j (w)\|_{C_p}^p = \sum |c_j|^p \ |f_j(w)|^p\ .
\end{equation}
Hence 
\begin{equation}\label{eq:6.22} 
\|\sum c_j f_j\|^p = \int \|\sum c_j f_j(w)\|_{C_p}^p\,dw = \sum |c_j|^p\ .
\end{equation} 
Now fix $n$, and let $H_n$ be the subspace of $C_p$ defined in the  proof 
of Lemma~\ref{lem:6.4} (specifically, in \eqref{eq:6.3}). 
Standard results yield that $L^p \otimes_p H_n$ embeds in $L^p$ 
(actually, $L^p \otimes_p H_n$ is isomorphic to $L^p$ if $p>1$),
and of course $I\otimes P_n$ is a projection onto $L^p\otimes_p H_n$ 
with $\|I\otimes P_n\|\le 2$ ($P_n$ as defined in the proof of 
\ref{lem:6.4}). 
Now let $X$ be as in Lemma~\ref{lem:6.7}, and suppose $X$ does not embed 
in $L_p$. 
Then for each $n$, we may choose an $x_n\in X$ with 
\begin{equation}
\|x_n\| =1 \ \text{ and }\ \|(I\otimes P_n)x_n\| <\frac1{2^n}\ .
\end{equation}
But it follows that for any $f\in L_p\otimes_pC_p$, 
\begin{equation}\label{eq:6.24} 
(I\otimes P_n) (f) \to f\ \text{ as }\ n\to\infty\ .
\end{equation}
A standard travelling hump argument now yields a normalized generalized 
block diagonal $(f_k)$ of $(\bone \otimes e_{ij})$ and a subsequence 
$(x'_j)$ of $(x_j)$ so that 
\begin{equation}\label{eq:6.25} 
\|x'_k - f_k\| < \frac1{2^k}\ \text{ for all }\ k\ .
\end{equation}
It follows immediately that $(x'_k)$ is equivalent to the $\ell^p$-basis.
\end{proof}

\begin{rem}
The last part of this proof also yields the fact (due to Y.~Friedman \cite{F})
that if $X$ is an 
infinite-dimensional subspace of $C_p$, then $X$ is isomorphic to $\ell^2$ 
or $\ell^p$ embeds in $X$. 
Indeed, assuming $X$ is not isomorphic to $\ell^2$, then since $H_n$ 
is isomorphic to $\ell^2$, we obtain for each $n$ and $x_n\in X$ with 
$\|x_n\|=1$ and $\|P_n x_n\| <\frac1{2^n}$. 
Again we then obtain a normalized block diagonal $(f_k)$ of $(e_{ij})$ 
and a subsequence $(x'_j)$ of $(x_j)$ satisfying \eqref{eq:6.25}, 
and then $(x'_k)$ is equivalent to the $\ell^p$ basis. 
\end{rem}

We now give the last and perhaps most delicate of the needed non-embedding 
results; its proof requires Theorem~\ref{thm:3.5}, 
the ``fine'' version of our Main Result.

\begin{thm}\label{thm:6.9} 
Let $\N$ be a finite von Neumann algebra. 
Then $L^p (\R) \otimes_p C_p$ is not isomorphic to a subspace of 
$L^p(\N) \oplus (L_p\otimes_p C_p)$.
\end{thm}

We first give some notation used in the proof. 
As always, $e_{ij}$'s denote the matrix units for $C_p$. 
Thus $L^p(\R)\otimes_p C_p = L^p (\R \bar\otimes B(\ell^2)) = $ 
the closed linear span of the elementary tensors $f\otimes e_{ij}$, 
$f\in L^p(\R)$, $i$ and $j$ arbitrary. 
We denote also the norm on $L^p(\R) \otimes_p C_p$ as $\|\cdot\|_p$. 
If $X$ is a closed linear subspace of $L^p(\R)$, 
\begin{equation}\label{eq:6.26} 
X\otimes_p C_p \defeq [x\otimes e_{ij} :x\in X,\ i,\ j\in \nat]
\end{equation}
(where the closed linear span above is taken in $L^p(\R)\otimes_p C_p$). 
Next, we need expressions for the norm on $L^p(\R)\otimes \Row$, 
$L^p(\R)\otimes$~Column. 
We easily see that given $x_1,\ldots,x_n$ in $L^p(\R)$, then for any $i$, 
\begin{equation}\label{eq:6.27} 
\Big\|\sum_{j=1}^n x_j \otimes e_{ij}\Big\|_p 
= \Big\|\biggl( \sum_{j=1}^n x_j x_j^*\biggr)^{1/2}\Big\|_p
\end{equation}
and 
\begin{equation}\label{eq:6.28}
\Big\|\sum_{j=1}^n x_j \otimes e_{ji}\Big\|_p 
= \Big\|\biggl( \sum_{j=1}^n x_j^* x_j\biggr)^{1/2} \Big\|_p\ .
\end{equation}

Evidently \eqref{eq:6.27} and \eqref{eq:6.28} show that if we consider 
a matrix of the form $(x_{ij}\otimes e_{ij})$ with $x_{ij}$ non-zero 
elements of $L^p(\R)$ for all $i$ and $j$, then all rows and columns 
of this matrix are 1-unconditional sequences.

The next result is really a ``localization'' of Lemma~\ref{lem:3.1} (and 
could be formulated instead for $L^p(\N)$, $\N$ any finite 
von~Neumann algebra). 

\begin{lem}\label{lem:6.10}
Let $X$ be a closed linear subspace of $L^p(\R)$ containing no subspace 
isomorphic to $\ell^p$. 
Then given $\ep>0$, there is an $N$ so that given any $n\ge N$ and 
$x_1,\ldots,x_n$ in $\Ba (X)$, 
\begin{equation}\label{eq:6.29} 
n^{-1/p} \Big\|\biggl( \sum_{i=1}^n x_i x_i^*\biggr)^{1/2}\Big\|_p 
\le \ep\ \text{ and }\ 
n^{-1/p} \Big\|\biggl( \sum_{i=1}^n x_i^* x_i\biggr)^{1/2}\Big\|_p 
\le \ep \ .
\end{equation}
\end{lem}

\begin{proof} 
Let $\T$ be the normal faithful tracial state in $\R$. 
By Theorem~\ref{thm:5.4}, $\{|x|^p :x\in \Ba(X)\}$ is uniformly integrable. 
Let $\eta>0$, to be decided later. 
Choose $\delta>0$ so that 
\begin{equation}\label{eq:6.30} 
\om (|x|^p,\delta) \le\eta^p\ \text{ for all }\ x\in \Ba (X)\ .
\end{equation}
Let $x_1,\ldots,x_n$ be elements of $\Ba (X)$. 
By the final statement of Lemma~\ref{lem:2.3}, we may choose for each $j$ a  
$P_j \in \P(\R)$ so that $x_j P_j \in \R$ with 
\begin{equation}\label{eq:6.31}
\|x_jP_j\|_\infty \le \delta^{-1/p} \ \text{ and }\ 
\|x_j (I-P_j)\|_p \le \eta\ .
\end{equation}
Then 
\begin{equation}\label{eq:6.32}
\begin{split}
\Big\| \biggl( \sum_{j=1}^n x_j x_j^*\biggr)^{1/2}\Big\|_p 
& = \Big\| \sum_{j=1}^n x_j \otimes e_{ij}\Big\|_p\ 
\text{ by \eqref{eq:6.27}}\\
& \le \Big\| \sum_{j=1}^n x_j P_j\otimes e_{1j}\Big\|_p 
+ \Big\| \sum_{j=1}^n x_j (I-P_j)\otimes e_{1j}\Big\|_p\ .
\end{split}
\end{equation}
Since $(x_j (I-P_j)\otimes e_{1j})_{j=1}^n$ is 1-unconditional and 
$L^p (\R) \otimes_p C_p$ is type $p$ with constant one, 
\begin{equation}\label{eq:6.33}
\begin{split}
\sum_{j=1}^n \|x_j (I-P_j)\otimes e_{1j}\|_p 
& \le \biggl( \sum_{j=1}^n \|x_j (I-P_j)\|_p^p\biggr)^{1/p}\\
& \le \eta n^{1/p}\ \text{ by \eqref{eq:6.31}}\ .
\end{split}
\end{equation}
Now 
\begin{equation}\label{eq:6.34}
\begin{split}
\Big\| \sum_{j=1}^n x_j P_j \otimes e_{1j}\Big\|_p 
& = \left[ \T\biggl( \sum_{j=1}^n x_j P_j x_j^* \biggr)^{p/2}\right]^{1/p}\\
& \le \left[\T \biggl( \sum_{j=1}^n x_j P_j x_j^*\biggr)\right]^{1/2}
\ \text{ (since $p<2$)}\\
& \le n^{1/2} \delta^{-1/p}\ \text{ by \eqref{eq:6.31}.}
\end{split}
\end{equation}
Thus \eqref{eq:6.32}--\eqref{eq:6.34} yield that 
\begin{equation}\label{eq:6.35} 
n^{-1/p} \Big\| \biggl( \sum_{j=1}^n x_j x_j^*\biggr)^{1/2}\Big\|_p 
\le \eta + \frac1{n^{\frac1p-\frac12}} \delta^{-1/p}\ .
\end{equation}
Evidently we now need only take $\eta\le \frac{\ep}2$; then choose $N$ so that 
$N^{-(\frac1p-\frac12)} \delta^{-1/p} \le \frac{\ep}2$; the identical 
argument for $(x_i^* x_i)_{i=1}^n$ now yields that \eqref{eq:6.29} 
holds for all $n\ge N$.
\end{proof}

We may now easily obtain our final needed preliminary result. 
(See the Remark following Theorem~\ref{thm:3.5} for the definition of: 
the rows or columns of a matrix contain  $\ell_n^p$-sequences.) 

\begin{cor}\label{cor:6.11} 
Let $X$ be a closed linear subspace of $L^p(\R)$ containing no subspace 
isomorphic to $\ell^p$, and let $(x_{ij})$ be a seminormalized matrix 
whose terms lie in $X$. 
Then the matrix $(x_{ij}\otimes e_{ij})$ in $X\otimes_p C_p$ has the 
following properties:
\begin{itemize}
\item[(i)] Neither the rows nor the columns contain $\ell_n^p$-sequences.
\item[(ii)] Every row and column is 1-unconditional.
\item[(iii)] Every generalized diagonal is equivalent to the usual 
$\ell^p$ basis.
\end{itemize}
\end{cor}

\begin{proof} 
(i) follows immediately from Lemma~\ref{lem:6.10} and \eqref{eq:6.27}, 
and the latter also immediately yields (ii). 
If  $(f_i)$ is a generalized diagonal of the matrix, then there exist 
projections $P_1,P_2,\ldots$, $Q_1,Q_2,\ldots$ in $\R\bar\otimes B(\ell^2)$ 
so that the $P_j$'s and the $Q_j$'s are pairwise orthogonal, with 
$f_j = P_j f_j Q_j$ for all $j$. 
(That is, $(f_j)$ is  ``right and left disjointly supported''.) 
It then follows that for any $n$ and scalars $c_1,\ldots,c_n$, 
\begin{equation}\label{eq:6.36} 
\Big\| \sum_{j=1}^n c_j f_j\Big\|_p 
= \biggl( \sum_{j=1}^n |c_j|^p \|f_j\|_p^p\biggr)^{1/p}\ ,
\end{equation}
which immediately yields (iii) since $(x_{ij}\otimes e_{ij})$ is 
semi-normalized.
\end{proof}

We are finally prepared for the 

\begin{proof}[Proof of Theorem~\ref{thm:6.9}]
Let $p<q<2$ and let $X$ be a subspace of $L^p(\R)$ so that $X$ is 
isomorphic to $C_q$ (using Junge's result, formulated as 
Theorem~\ref{thm:6.6} above).  
We claim that $X\otimes_p C_p$ is not isomorphic to a subspace of 
$L^p(\N) \oplus (L_p \otimes_p C_p)$ (which of course proves 
Theorem~\ref{thm:6.9}). 
Suppose to the contrary that $T:X\otimes_p C_p \to L^p (\N) 
\oplus (L_p\otimes_p C_p)$ is an isomorphic embedding. 
Assume without loss of generality that $\|T\|=1$. 
Let $\ep>0$ be chosen so that $\|Tf\| \ge \ep \|f\|$ for all 
$f\in X \otimes_p C_p$. 
Let $P$ denote the projection of $L^p(\N) \oplus (L_p \otimes_p C_p)$ 
onto $L^p(\N)$, with kernel $L_p \otimes_p C_p$; and set $Q= I-P$. 
Now fix $i$ and $j$. 
Then of course $X\otimes e_{ij}$ is isometric to $X$. 
Thus by Lemma~\ref{lem:6.7}, $QT|(X\otimes e_{ij})$ cannot be an 
isomorphic embedding (that is, $C_q$ does not embed in 
$L_p\otimes_p C_p$). 
Hence we may choose $x_{ij}\in X$ with  
\begin{equation}\label{eq:6.37} 
\|x_{ij}\| =1\ \text{ and }\ \|QT(x_{ij} \otimes e_{ij})\| 
< \frac{\ep}{2^{i+j+2}}\ .
\end{equation}

Now let $Y = [x_{ij}\otimes e_{ij} : i,j=1,2,\ldots]$. 
Since $\ell^p$ does not embed in $X$, the conclusion of 
Corollary~\ref{cor:6.11} holds for the matrix $(x_{ij}\otimes e_{ij})$.

It follows from \eqref{eq:6.37} that 
\begin{equation}\label{eq:6.38} 
\|QT|Y\| < \frac{\ep}2\ .
\end{equation}
Hence we obtain that 
\begin{equation}\label{eq:6.39} 
\|PT (y)\| \ge \frac{\ep}2 \|y\|\ \text{ for all }\ y\in Y\ .
\end{equation}
Thus $Y$ is isomorphic to a subspace $Z$ of $L^p(\N)$. 
Let $z_{ij} = PT (x_{ij} \otimes e_{ij})$ for all $i$ and $j$. 
Now since $PT|Y$ is an isomorphism, Corollary~\ref{cor:6.11} 
yields that there is a $u$ so that every row and column of $(z_{ij})$ 
is $u$-conditional, every generalized diagonal of $(z_{ij})$ is 
equivalent to the $\ell^p$-basis, yet neither the rows nor the columns 
of $(z_{ij})$ contain $\ell_n^p$-sequences. 
This is impossible by Theorem~\ref{thm:3.5}. 
\end{proof}

The following result is an immediate consequence of Theorem~\ref{thm:6.9} 
and known structural results for von-Neumann algebras. 

\begin{cor}\label{cor:6.12} 
Let $\N,\M$ be von Neumann algebras so that $\M$ has a direct summand 
of type II$_\infty$ or of type III. 
If $L^p(\M)$ is Banach isomorphic to a subspace of $L^p(\N)$, then 
also $\N$  has a direct summand of type II$_\infty$ or of type~III.
\end{cor}

\begin{proof} 
The hypotheses imply (via known results, cf.\ \cite{HS}) that 
$\R\bar\otimes B(\ell^2)$ is isomorphic to a von~Neumann subalgebra of 
$\M$, which is the range of a normal conditional expectation, whence 
$L^p(\R) \otimes_p C_p$ is completely isometric to a subspace of $L^p(\M)$. 
Since $L^p(\R)\otimes C_p$ is separable, we can assume without loss 
of generality that $\N$ acts on a separable Hilbert space. 
Then if $\N$ fails the conclusion, there exists a finite von~Neumann 
algebra $\tilde{\N}$ so that $\N$ is isomorphic to a subalgebra of 
$\tilde{\N} \oplus (L^\infty \bar\otimes B(\ell^2))$, and hence  
$L^p(\N)$ is completely isometric to a subspace of $L^p (\tilde{\N}) 
\oplus (L_p\otimes_p C_p)$. 
But then $L^p (\M)$ does not Banach embed in $L^p (\N)$, since 
$L^p (\R) \otimes_p C_p$ does not embed in $L^p(\tilde{\N})\oplus 
(L_p \otimes_p C_p)$ by Theorem~\ref{thm:6.9}.
\end{proof}

\begin{rem}
Of course Corollary~\ref{cor:6.8} (i.e., the results of Junge and 
Pisier-Xu cited above) also immediately yields that if $\M$ and $\N$ 
are von~Neumann algebras so that $\M$ has a type II$_1$ summand, 
and $L^p(\M)$ embeds in $L^p(\N)$, then $\N$ must have also have a 
summand of type~II or type~III. 
Combining these two results, we have that
{\it if $L^p(\M)$ is Banach isomorphic to a subspace of $L^p(\N)$ and 
$\M$ has no type~III summand, then $\N$ has a direct summand of type 
at least as large as these of the summands of $\N$.} 
It remains a most intriguing problem to see if one can eliminate the 
non-type~III summand hypothesis in this statement.
\end{rem}

We now complete the proof of Theorem~\ref{thm:6.2}. 
We shall formulate the ``positive'' results in the language of operator 
spaces; the reader unfamiliar with the relevant terms may just ignore 
the adjective ``complete'' in all the statements, for of course all 
positive operator space claims imply the pure Banach space ones. 
Given operator spaces $X$ and $Y$, let us say that $X$ 
{\it completely contractively factors through\/} $Y$ if $X$ is completely 
isometric to a subspace $X'$ of $Y$ such that there exists a completely 
contractive projection mapping $Y$ onto $X'$. 
Equivalently, there exist complete contractions $U:X\to Y$ and $V:Y\to X$ 
such that $V\circ U = I_X$, $I_X$ the identity operator on $X$, that is, 
\begin{equation}\label{eq:6.40}
\begin{picture}(100,70)
\put(30,55){$Y$}
\put(10,20){\vector(1,2){15}}
\put(7,32){$\scriptstyle U$}
\put(4,4){$X$}
\put(40,50){\vector(1,-2){15}}
\put(52,32){$\scriptstyle V$}
\put(22,8){\vector(1,0){25}}
\put(28,12){$\scriptstyle I_X$}
\put(54,4){$X$}
\end{picture}
\qquad\qquad \raise3ex\hbox{.}
\end{equation}
Now we easily see that 
\begin{equation}\label{eq:6.41}
(L^p (\R) \oplus L^p(\R) \oplus\cdots)_p \text{ completely contractively 
factors through } L^p (\R)\ .
\end{equation}
Indeed, simply let $P_1,P_2,\ldots$ be pairwise orthogonal non-zero 
projections in $\R$. 
As is well known, then $P_i \R P_i$ is isomorphic to $\R$ and hence 
$P_i L^p (\R) P_i$ is completely isometric to $L^p(\R)$ for all $i$; 
then the map  on $L^p(\R)$ defined by $f\to \sum P_i fP_i$ witnesses 
\eqref{eq:6.41}. 

Since $\R \bar\otimes \R$ is isomorphic to $\R$, 
\begin{equation}\label{eq:6.42} 
L^p(\R) \otimes_p L^p(\R) \defeq L^p(\R\bar\otimes \R) 
\text{ is completely isometric to } L^p(\R)\ .
\end{equation}

Using \eqref{eq:6.41} and \eqref{eq:6.42}, we may now easily see that if 
$Y$ is immediately below $X$ in the tree (and lying on a branch), then $X$ 
completely contractively factors through $Y$. 
Using the notation $X\overset{cc}{\hookrightarrow} Y$ to mean that 
$X$ completely contractively factors through $Y$, we see, e.g., that 
$L_p \ccarrow L^p(\R) \implies L_p\otimes_p C_n^p \ccarrow 
L^p(\R) \otimes_p C_p^n \ccarrow L^p (\R)\otimes_p L^p(\R)$, 
whence 
\begin{equation*}
L_p\otimes_p S_p = \biggl( \bigoplus_{n=1}^\infty (L_p \otimes_p C_p^n)
\biggr)_p \ccarrow \biggl( \bigoplus_{n=1}^\infty L_p \otimes_p L^p(\R)
\biggr)_p \ccarrow L^p(\R)\ ,
\end{equation*}
i.e., 
\begin{equation}\label{eq:6.43} 
L_p \otimes_p S_p \ccarrow L^p (\R)\ .
\end{equation}

Writing $X\approx Y$ to mean: $X$ is completely isometric to $Y$, 
we have 
\begin{equation}\label{eq:6.44} 
C_p \oplus (L_p\otimes_p S_p) \ccarrow C_p \oplus L_p \otimes_p C_p 
\ccarrow (L_p \otimes C_p) \otimes (L_p\otimes C_p) 
\approx L_p \otimes C_p 
\end{equation}
(where we use $\ell^p$-direct sums).

$X\ccarrow Y$ if $X$ is the level 7 space and $Y$ is the level 8 space, 
since the same argument for 
\eqref{eq:6.41} yields also 
\begin{equation}\label{eq:6.45} 
\Big( (L^p (\R) \otimes_p C_p ) \oplus (L^p (\R)\otimes_p C_p) 
\oplus \cdots \Big) \ccarrow L^p (\R) \otimes_p C_p\ .
\end{equation}
The reader may now easily check that the remaining ``positive'' assertions 
on the tree. 
For the far deeper negative assertions, let us use the notation: 
$X\not\hookrightarrow Y$ to mean that the Banach space $X$ is not isomorphic 
to a subspace of $Y$. 

Now suppose $X\ne Y$ are on the tree and $Y$ cannot be connected to $X$ 
by a descending branch; {\it we claim that\/} $X\not\hookrightarrow Y$. 

It suffices to prove this assertion by showing by induction on 
$j=2,3,\ldots$ that $X$ lies at level $j$ and 
\begin{align}
&\text{there is a $k\ge j$ so that $Y$ is at the $k^{\text{th}}$ level,
but if $Z$ is a higher}
\label{eq:6.46}\\
&\text{level than $k$, connected to $Y$, $Z\ne X$,
then $X$ is connected to $Z$} 
\nonumber\\
&\text{and moreover there is no $X'$ connected to $X$ but not to $Y$ with}
\nonumber\\
&\text{level $X'<j$}\nonumber
\end{align}
or 
\begin{align}
&\text{$Y$ is at the $(j-1)^{\text{st}}$ level, 
but if $Y$ is connected to $Z$ at level $k\ge j$}\label{eq:6.47}\\
&\text{with $Z\ne X$, then $X$ is connected to $Z$ and moreover 
if $Z$ is connected}
\nonumber\\
&\text{to $X$ with level $Z<j$, then $Z$ is connected to $Y$.}
\nonumber
\end{align}
\begin{itemize}
\item[$j=2$.] 
$S_p\not\hookrightarrow L_p$ is classical (and also follows from our 
Corollary~\ref{cor:1.4}). $L_p \not\hookrightarrow C_p$ since $\ell_q 
\hookrightarrow L_p$ if $p<q<2$ but $\ell_q\not\hookrightarrow C_p$.
\item[$j=3$.] 
$C_p\not\hookrightarrow L^p(\R)$, the main result of the paper.
\item[$j=4$.] 
$L_p\otimes_p S_p \not\hookrightarrow C_p \oplus L_p$ by 
Proposition~\ref{prop:6.5}. 
\item[$j=5$.] 
$L^p(\R) \not\hookrightarrow L_p\otimes_p C_p$ by Corollary~\ref{cor:6.8}.
\item[$j=6$.] 
$L_p\otimes_p C_p \not\hookrightarrow C_p\oplus L^p(\R)$ 
by Theorem~\ref{thm:6.3}.
\item[$j=7$.]
There is no $Y$ satisfying \eqref{eq:6.46} or \eqref{eq:6.47}.
\item[$j=8$.] 
Theorem~\ref{thm:6.9} gives the one required non-embedding result.
\end{itemize}

This completes the proof of the final statement of Theorem~\ref{thm:6.2}.
It remains to prove the first statement. 
This follows via the known type-decomposition and structure of 
hyperfinite von-Neumann algebras, and the following operator space version
of the Pe{\l}czy\'nski decomposition method (whose proof is exactly as 
Pe{\l}czy\'nski's proof for the Banach space case \cite{P}; 
see also p.54 of \cite{LT} and \cite{Ar}). 

\begin{lem}\label{lem:6.13}
Let $X$ and $Y$ be operator spaces so that 
\begin{itemize}
\item[(i)] each completely factors through the other 
\end{itemize}
and so that either 
\begin{itemize}
\item[(ii)] $X$ is completely isomorphic to 
$X\oplus X$ and $Y$ is completely isomorphic to $Y \oplus Y$ 
\end{itemize}
or
\begin{itemize}
\item[(ii$'$)] $X$ is completely isomorphic to $(X\oplus X\oplus \cdots)_q$ 
for some $q\in [1,\infty]$.
\end{itemize}
Then $X$ and $Y$ are completely isomorphic. 
\end{lem}

(We say that $X$ completely factors through $Y$ if $X$ is completely 
isomorphic to a completely complemented subspace of $Y$.) 

\begin{cor}\label{cor:6.14} 
If $(X\oplus X\oplus\cdots)_p$ completely factors through the operator 
space $X$, then $X$ is completely isomorphic to $(X\oplus X\oplus\cdots)_p$.
\end{cor}

\noindent {\it End of the proof of Theorem 6.2.}
$(X\oplus X\oplus \cdots)_p$ completely contractively factors 
through $X$ for all of the 13 spaces $X$ listed in Theorem~\ref{thm:6.2} 
(applying \eqref{eq:6.41}, \eqref{eq:6.45}, and the analogous results 
for $C_p$, $L_p$, and $L_p\otimes_p C_p$). 
Thus the conclusion of \ref{cor:6.14} applies.

Now let $\N$ be as in the statement of Theorem~\ref{thm:6.2}. 
If $\N$ is type~I, then by the results in \cite{S2} 
$L^p(\N)$ is completely isomorphic to one of 
the first nine spaces listed in Theorem~\ref{thm:6.1}, so assume that 
$\N$ is not type~I. 
Then we have that 
\begin{equation*}
\N = \N_{\text{I}} \oplus \N_{\text{II}_1} \oplus 
\N_{\text{II}_\infty}\ ,
\end{equation*}
where for each $i$, $\N_i = \{0\}$ or $\N_i$ is a hyperfinite 
von~Neumann algebra of type $i$, so that also 
$\N_{\text{II}_1} \oplus \N_{\text{II}_\infty}\ne 0$. 

Now suppose that $\N$ is finite. 
It then follows from work of A.~Connes \cite{C2} that 
\begin{equation}\label{eq:6.48new}
\N_{\text{I}} \oplus \N_{\text{II}_1} 
\text{ is isomorphic to a von-Neumann subalgebra of }
\R\ .
\end{equation}
Indeed, by disintegration and Proposition~6.5 of  \cite{C2}, any finite 
hyperfinite von~Neumann algebra (with separable predual) is a countable 
$\ell^\infty$-direct sum of von~Neumann algebras of the form 
$\A\bar\otimes\B$, where $\A$ is abelian and $\B$ is either $M_n$ for some 
$n<\infty $ or $\R$. 
But such an algebra $\A\bar\otimes\B$ can be realized as a sub-algebra 
of $\R$; since also $\R\bar\otimes\R$ is isomorphic to $\R$, and 
$(\R\oplus \R\oplus\cdots)_{\ell^\infty}$ is (isomorphic to) a 
von~Neumann subalgebra of $\R$, \eqref{eq:6.48new} holds. 
Since $\N_{\text{II}_1}\ne0$, we have by the above discussion that also 
\begin{equation}\label{eq:6.49new} 
\R\text{ is isomorphic to a von-Neumann subalgebra of }\N\ .
\end{equation} 

Thus, we have that if $\A$ or $\B$ equals $\N$ or $\R$, then 
\begin{align}
&\text{$\A$ is (isomorphic to) a subalgebra of $\B$, which is}
\label{eq:6.50}\\
&\text{the range of a normal conditional expectation.}
\nonumber
\end{align}

Now if \eqref{eq:6.49new} holds for any two von Neumann algebras $\A$ and 
$\B$, then $L^p(\A)$ completely contractively factors through $L^p(\B)$. 
Thus by Lemma~6.13 and Corollary~6.14 applied to $X= L^p(\R)$, we obtain 
that $L^p(\N)$ is isomorphic to $L^p(\R)$.

Now if $\N_{\text{II}_\infty} \ne 0$,  again using the deep results in 
\cite{C2}, 
$\N_{\text{II}_\infty}$ is (isomorphic to) 
$\M\bar\otimes B(\ell^2)$ where $\M$ is a finite hyperfinite von Neumann 
algebra, whence letting $\A$ and $\B$ equal 
$\N$ or $R\bar\otimes B(\ell^2)$, \eqref{eq:6.48new} holds, whence 
$L^p(\N)$ is completely isomorphic to $L^p (\R) \otimes_p C_p$ again by 
Lemma~\ref{lem:6.13} and Corollary~\ref{cor:6.14} applied to 
$L^p (\R) \otimes_p C_p$. 

Now assume $\N_{\text{II}_\infty} = \{0\}$, so $\N_{\text{II}_1} \ne\{0\}$, 
and suppose $\N$ is infinite; since $\N_{\text{II}_\infty} =\{0\}$, 
we must have that $\N_I$ is infinite. 
But then by the classification of the $L^p$ spaces of type~I algebras, 
we have that $L^p(\N_I)$ is completely isomorphic to either $C_p$, 
$L_p\otimes C_p$, $C_p \oplus L_p$, or $C_p \oplus (L_p\otimes_p S_p)$. 

But $C_p \oplus L_p \oplus L^p(\R)$ and 
$C_p\oplus (L_p\otimes_p S_p) \oplus L^p(\R)$ 
are both completely isomorphic to 
$C_p\oplus L^p(\R)$, by our analysis of the finite case. 
Hence $L^p(\N)$ is completely isomorphic either to $C_p \oplus L^p(\R)$ 
or to $(L_p\otimes_p C_p)\oplus L^p(\R)$, completing the 
entire proof.\qed

\section{$L^p(\N)$-isomorphism results for $\N$ type III hyperfinite 
or a free group von Neumann algebra}		
\setcounter{equation}{0}

We first formulate the results of this section for the case of preduals 
of von~Neumann algebras $\N$, i.e., $L^1(\N)$, and then show they hold 
also for the spaces $L^p(\N)$ for $1<p<\infty$, as in the preceding sections. 
The following result is an immediate consequence of Corollary~6.12. 
We prefer to give a quick proof just using Corollary~1.2. 

\begin{thm}\label{thm:7.1} 
Let $\N$ be a factor of type II$_1$ and let $\M$ be a factor of 
type II$_\infty$ or type III. 
Then the preduals $\N_*$ and $\M_*$ are not Banach space isomorphic. 
\end{thm}

\begin{proof} 
By the assumptions $\M$ is a properly infinite von~Neumann algebra, i.e.,
$\M \cong \M \bar\otimes  B(\ell^2)$ as von~Neumann algebras (where 
$\bar \otimes$ is the standard von~Neumann algebra tensor product). 
In particular $\M_*$ is isometrically isomorphic to $\M_* \otimes_\gamma C_1$ 
for some crossnorm $\gamma$ on the algebraic tensor product 
$\M_* \otimes C_1$, and therefore $C_1$ imbeds isometrically in $\M_*$. 
By Corollary~1.2, $C_1$ does not Banach space imbed in $\N_*$.
\end{proof} 

It would be interesting to know, whether a type II$_\infty$-factor and a 
type III-factor can be distinguished by the Banach space isomorphism classes 
of their preduals. 
(As noted in the Introduction, we do not know the answer for the special 
case of injective factors.)
In \cite{C1} Connes introduced a subclassification of factors of 
type III into factors of type III$_\lambda$, where $\lambda$ can take any 
value in the closed interval $[0,1]$. 
Theorem~7.2 below shows that the number $\lambda$ in this classification 
cannot be determined by the Banach space isomorphism class (or even 
operator space isomorphism class) of the predual. 
Recall from \cite{C2} and \cite{H}, that for each $\lambda\in (0,1]$, there 
is up to von~Neumann algebra isomorphism only one injective factor of 
type III$_\lambda$ acting on a separable Hilbert space. 
For $0<\lambda <1$ it is the Powers factor 
\[ 
R_\lambda = \bigotimes_{n=1}^\infty (M_2 (\complex),\varphi_\lambda)
\] 
where $\varphi_\lambda$ is the state on the $2\times 2$ complex matrices 
given by 
\[ 
\varphi_\lambda \left(\begin{matrix} x_{11}&x_{12}\\ x_{21}&x_{22}
\end{matrix}\right) = \frac{\lambda}{1+\lambda} x_{11} + 
\frac{1}{1+\lambda} x_{22} 
\] 
and for $\lambda =1$ it is the Araki-Woods factor $R_\infty$, which 
can be obtained (up to von~Neumann-isomorphism) as the tensor product 
of two Powers factors 
\[ R_\infty \cong R_{\lambda_1} \bar\otimes R_{\lambda_2}
\] 
provided $\frac{\log \lambda_1}{\log \lambda_2}\notin \que$.
On the hand there are uncountably many injective factors of type III$_0$ 
acting on a separable Hilbert space (cf.\ \cite{C1}, \cite{C2}). 
We will consider the predual of a von~Neumann algebra as an operator 
space with the standard dual operator space structure (cf. \cite{Bl}). 

\begin{thm}\label{thm:7.2} 
Let for $0<\lambda <1$, $R_\lambda$ denote the Powers factor of type 
III$_\lambda$ and let $R_\infty$ denote the Araki-Woods factor of 
type III$_1$.
\begin{itemize}
\item[(a)] For every $\lambda \in (0,1)$ the predual $(R_\lambda)_*$ 
is completely isomorphic to $(R_\infty)_*$.
\item[(b)] There is an uncountable family $(\N_i)_{i\in I}$ of mutually 
non-isomorphic (in the von~Neumann algebra sense) injective 
type III$_0$-factors on a separable Hilbert space for which $(\N_i)_*$ 
is completely isomorphic to $(R_\infty)_*$.
\end{itemize}
\end{thm} 

\begin{rem}
In \cite{ChrS}, Christensen and Sinclair proved that all injective infinite
dimensional factors acting on separable Hilbert space are completely 
isomorphic. 
This does not imply that their preduals are completely isomorphic. 
Indeed the unique injective type II$_1$-factor $\R$ and the unique injective 
type II$_\infty$-factor $\R \bar\otimes B(\ell^2)$ have non-isomorphic 
preduals by Theorem~7.1. 
Theorem~7.2 as well as the results in \cite{ChrS} are based on the completely 
bounded version of the Pe{\l}czy\'nski decomposition method stated 
as Lemma~6.13 above. 
\end{rem}

\begin{proof}[Proof of Theorem 7.2]
(a) Let $0<\lambda <1$ and put $\N = R_\lambda$,  $\M= R_\infty$. 
Since $\N$ is a properly infinite von~Neumann algebra, there exists two 
isometries $u_1,u_2\in \N$, such that $u_1u_1^*$ and $u_2 u_2^*$ are two 
orthogonal projections with sum~1. 
Define now 
\[ \Phi :\N \to \N \oplus \N \ \text{ by }\ \Phi (x) = (u_1^* x,u_2^*x)
\] 
and 
\[ \Psi :\N\oplus \N \to \N\ \text{ by }\ \Psi  (x,y) = (u,x+u_2y)
\] 
Then $\Phi \circ\Psi = \id_{\N\oplus \N} $ and $\Psi\circ\Phi =\id_{\N}$. 
Since $\Phi$ and $\Psi$ are normal (i.e., continuous) in the 
$\omega^*$-topologies on $\N$ and $\N\oplus \N$) and also are completely 
bounded maps, it follows that $\N_* \approx_{\cb} \N_* \otimes \N_*$. 
Similary we have $\M_*\approx_{\cb}\M_* \oplus \M_*$. 
Thus the pair $(\M_*,\N_*)$ satisfies (ii) in Lemma~6.13. 
We next check condition (i) in Lemma~6.13. 

Since $R_\infty \cong R_\lambda \bar\otimes R_\infty$ as von~Neumann algebras 
(cf.\ \cite[Sect.3.6]{C1}), we can without loss of generality assume 
that $\M = \N \bar\otimes \P$ where $\P \cong R_\infty$. 
Let $\varphi$ be a normal faithful state on $\P$ and define 
\[ \pi :\N\to \N\bar\otimes \P\ \text{ by }\ \pi (x) = x\otimes \bone\ ,
\]
and let $\rho :\N\bar\otimes \P \to \N$ be the left slice map given by 
$\varphi$, i.e., the unique normal linear map $\N \bar\otimes P\to\N$ 
for which 
\[ \rho (x\otimes y) = \varphi (y)x\ ,\qquad x\in \N \ ,\ y\in \P\ .
\]
Thus $\|\pi \|_{\cb}  = \|\rho\|_{\cb} =1$ and $\rho \circ\pi = \id_{\N}$. 
Hence $\id_{\N_*}$ has a completely bounded factorization through 
$\M_*$, i.e., $\N_*$ is $\cb$-isomorphic to a $\cb$-complemented subspace 
of $\M_*$. 
To prove the converse, we use that if $\varphi$  is a normal faithful 
state on the III$_1$-factor $\M=R_\infty$ and  $\alpha =\sigma_{t_0}^\varphi$
is the moduluar automorphism associated with $\varphi$ at $t_0 = 
-\frac{2\pi}{\log \lambda}$, then the crossed product $R_\infty \rtimes_\alpha 
\zed$ is a factor of type III$_\lambda$ (cf.\ \cite[proof of Lemma 2.9]{HW}).
Moreover injectivity of $R_\infty$ implies that the crossed product is 
injective (cf.\ \cite{C2}). 
Hence $R_\infty \rtimes_\alpha \zed \cong R_\lambda$ as von~Neumann 
algebras, so in this part of the proof we may assume that $\M\rtimes_\alpha
\zed=\N$. 
Further, after identifying $\M$ with its natural imbedding in the crossed 
product, we have that $\N$ is generated as a von~Neumann algebra by $\M$ 
and a certain unitary group $\{u^n \mid n\in\zed\}$ coming from the crossed 
product construction  (cf.\ \cite{C1}). 
Let $i:\M\hookrightarrow \M\rtimes_\alpha\zed$ be the imbedding and let 
$\ep :\M\rtimes_\alpha\zed\to i(\M)$ be the unique  normal faithful 
conditional expectation of $\M\rtimes_\alpha\zed$ onto $i(\M)$ for which 
$\ep (u^n)=0$, for $n\in\zed\smallsetminus \{0\}$ (see again \cite{C1}). 
Then $i$ and $\ep$ are normal maps and $i^{-1}\circ \ep\circ i=\id_{\M}$,  
so as above, we obtain that $\M_*$ is $\cb$-isomorphic to a 
$\cb$-complemented subspace of $\N_*$. 
Hence a) follows from Lemma~6.13. 

(b) Put again $\M = R_\infty$ and let $G\subseteq \real$ be a dense 
countable subgroup. 
Let $\varphi$ be a normal faithful state on $R_\infty$ and put 
$\N = R_\infty \rtimes_\alpha G$ where $\alpha :G\to \Aut (\M)$ is the 
restriction of the modular automorphism group $(\sigma_t^\varphi)_{t\in\real}$
to $G$. 
It follows from \cite{C1} (see the proof of \cite[Lemma 2.9]{HW}) 
that $\N_G$ is a factor of type III$_0$, 
which is also injective (by \cite{C2}). 
Moreover $T(\N_G) = G$, where $T$ is Connes $\pi$-invariant. 
Hence $G\ne G'$ implies, that $\N_G$ and $\N_{G'}$ are not 
von~Neumann-algebra isomorphic. 
It is easy to check, that there are uncountably many dense countable 
subgroups of $\real$. 
Put $\P = \N_G \bar\otimes R_\infty$. 
Since $R_\infty \bar\otimes R_\lambda\simeq R_\infty$  for $0<\lambda <1$, 
we have $\P\bar\otimes R_\lambda \cong \P$, $0<\lambda<1$, which by 
\cite[Theorem 3.6.1]{C1} implies that $\P$ is a factor of type III$_1$. 
Since $\P$ is also injective we have 
\[ \N_G \bar\otimes R_\infty \cong R_\infty = \M
\] 
as von Neumann algebras. 
As in the proof of (a), it now follows, that $\M_*$ is $\cb$-isomorphic 
to a $\cb$-complemented subspace of $(\N_G)_*$. 
Moreover, since $\M\rtimes_\alpha G$ is a crossed product with respect to 
a discrete group, there is again an embedding 
$i:\M\to\M\rtimes_\alpha G$ and a normal faithful conditional expectation 
$\ep : \M\rtimes_\alpha G\to i(\M)$, and the rest of the proof of (b) 
follows now exactly as in the proof of (a).
\end{proof}  

Let $L(F_n)$ denote the von Neumann algebra associated with the free 
group $F_n$ on $n$ generators. 
Then for $2\le n\le\infty$ $L(F_n)$ is a factor of type II$_1$. 
It is a long standing open problem to decide whether these 
II$_1$-factors are isomorphic as von~Neumann algebras. 
Due to work of Voiculescu, Dykema and Radulescu, it is known that either 
these factors are all isomorphic or $L(F_{n_1})\not\cong L(F_{n_2})$ whenever 
$2\le n_1,n_2\le\infty$ and $n_1\ne n_2$ (cf.\ \cite{VDN}). 
In \cite{Ar} Arias proved that the von~Neumann algebras $L(F_n)$, 
$2\le n\le\infty$ are isomorphic  as operator spaces. 
We show below, that also their preduals are  isomorphic as operator spaces. 
While Arias' proof uses mainly group theoretical considerations, the proof 
of Theorem~7.3 below relies on one rather deep result of Voiculescu, 
that $L(F_\infty) \cong M_k(L(F_\infty))$ as von~Neumann algebras for 
$k=2,3,\ldots$ (cf.\ \cite{Vo} or \cite{VDN}). 

\begin{thm}\label{thm:7.3} 
$L(F_n)_*$ is $\cb$-isomorphic to $L(F_\infty)_*$ for $n=2,3,\ldots$.
\end{thm}

\begin{proof} 
Let $n\in\nat$, $n\ge2$ and put $\N = L(F_n)$ and $\M= L(F_\infty)$. 
Since $F_n$ is isomorphic to a subgroup of $F_\infty$ and vice versa, 
$\N$ is von~Neumann-algebra isomorphic  to a subfactor $\N_1$ of $\M$ 
and $\M$ is von~Neumann-algebra isomorphic to a subfactor $\M_1$ of $\N$ 
(see \cite{Ar} for details). 
Moreover, let $\tau_{\M}$ and $\tau_{\N}$ be the unique normal faithful 
tracial states on $\M$ and $\N$ respectively. 
Then there is a unique normal faithful conditional expectation $\ep:\M
\xrightarrow{\text{onto}} \N_1$ preserving the trace $\tau_{\M}$ 
(resp.\ a unique normal faithful conditional expectation 
$\ep' :\N \xrightarrow{\text{onto}} \M$, preserving the trace $\tau_{\N}$). 
As in the proof of Theorem~7.2, this implies that $X= \M_*$ and $Y=\N_*$ 
satisfy condition (i) in Lemma~6.13. 
We next prove that (ii$'$) in Lemma~6.13 is satisfied with $q=1$. 
Since $\M= L(F_\infty)$ is a II$_1$-factor, we can choose a sequence of 
orthogonal projections $(p_i)_{i=1}^\infty$ in $\M$, such that 
$\tau (p_i)=2^{-i}$ and $\sum_{i=1}^\infty p_i=1$ (convergence in the 
strong operator topology). 
By Voiculescu's result quoted above, $L(F_\infty) \cong M_{2^i} (L(F_\infty))$
for $i=1,2,\ldots$ as von~Neumann-algebras, which implies that 
$p_i \M p_i \cong \M$ as von~Neumann-algebras. 

Indeed, Voiculescu's result yields that there are orthogonal equivalent 
projections $q_1,\ldots,q_{2^i}$ in $\M$ with $\sum_{j=1}^{2^i} q_j=\bone$ 
so that $q_1\M q_1\cong \M$. 
It follows (by uniqueness of $\tau_{\M}$) that $\tau (q_j) = \tau (q_{j'})$, 
for all $j$ and $j'$, and so $\tau(q_1) = 2^{-i}$. 
Since also $\tau_{\M} (P_i)=2^{-i}$ and $\M$ is a finite factor, 
$q_1$ and $p_i$ are equivalent, and hence $p_i\M p_i\cong q_1 \M q_1
\cong \M$ as desired. 

Put 
\[ Q = (\M \oplus \M \oplus \cdots)_{\ell^\infty} = \M \bar\otimes \ell^\infty
\ .\]
Then $Q$ is a von~Neumann algebra isomorphic to $Q_1 = \sum^{\oplus} p_i \M
p_i$, which is a von~Neumann subalgebra of $\M$. 
Moreover, there is a $\tau_{\M}$-preserving normal faithful conditional 
expectation $\ep'' : \M \xrightarrow{\text{onto}} Q_1$. 
Hence $Q_*$ is $\cb$-isomorphic to a $\cb$-complemented subspace of $\M_*$. 
Put as above $X= \M_*$. 
Then $Q_* = (X\oplus X\oplus\cdots)_{\ell^1}$ as operator spaces. 
Hence we have shown that  $(X\oplus X\oplus \cdots)_{\ell^1}$ completely 
factors through $X$, so $X$ and $(X\oplus X\oplus\cdots)_{\ell^1}$ are 
completely isomorphic by Corollary~6.14. 
This proves (ii$'$) iin Lemma~6.13 with $q=1$. 
Hence $X=\M_*$ and $Y = \N_*$ are completely isomorphic.
\end{proof}

In the rest of this section, we will show how Theorem~7.2 and Theorem~7.3 
can be generalized to the non-commutative $L^p$-spaces associated with the 
von~Neumann algebras in question. 
In \cite{Ko}, Kosaki proved, that the abstract $L^p$-spaces $L^p(\M)$, 
$1<p<\infty$ associated with a $\sigma$-finite ($=$ countably decomposable) 
von~Neumann algebra $\M$, can be obtained by the complex interpolation 
method applied to the pair $(\M,\M_*)$ with the imbedding $M\hookrightarrow 
\M_*$ given by the map $x\to x\varphi$, $x\in \M$, for a fixed normal 
faithful state $\varphi$ on $\M$. 
Assume next that $\N$ is a von~Neumann subalgebra of $\M$ and $\ep:\M\to\N$ 
is a normal faithful conditional expectation of $\M$ onto $\N$. 
By replacing $\varphi$ by $\varphi\circ\ep$, we can assume, that the  
state $\varphi$ used in Kosaki's imbedding is $\ep$-invariant. 
Next, the adjoint of $\ep$ defines an imbedding of $\N_*$ in $\M_*$ 
and $i^*$, the adjoint of the inclusion map $i:\N \to \M$ defines a 
$\cb$-contraction of $\M_*$ onto $\N_*$. 
Moreover, we have the following commuting diagram: 
\[
\begin{matrix}
\N&\xrightarrow{i}& \M &\xrightarrow{\ep} &  \N\\
\big\downarrow&&\big\downarrow&&\big\downarrow\\
\N_* &\xrightarrow{\ep^*}& \M_* &\xrightarrow{i^*}&\N_*
\end{matrix}
\]
where the vertical arrows are the Kosaki inclusions with 
respect to $\varphi_{1\N}$, $\varphi$ and $\varphi_{1\N}$ respectively. 
By the complex interpolation method we now get contractions 
$i_p :L^p(\N) \to L^p(\M)$ and $\ep_p : L^p (\M) \to L^p(\N)$, such that 
the following diagram commutes: 
\[
\begin{matrix}
\N & \xrightarrow{i} & \M &\xrightarrow{\ep} &\N&\\
\big\downarrow&&\big\downarrow&&\big\downarrow&\\
L^p(\N)&\xrightarrow{i_p}&L^p (\M)& \xrightarrow{\ep_p} &L^p(\N)&\\
\big\downarrow&&\big\downarrow&&\big\downarrow&\\
\N_* &\xrightarrow{\ep^*} &\M_* & \xrightarrow{i^*}&\N&\qquad .
\end{matrix}
\] 
Further, if we consider $L^p(\N)$ and $L^p(\M)$ as operator spaces with the 
operator spaces structure introduce by Pisier in \cite{Pi1}, we get that 
$i_p$ and $\ep_p$ are complete contractions. 
Hence we have proved: 

\begin{lem}\label{lem:7.4} 
Let $\M$ be a $\sigma$-finite von Neumann algebra, and $\N\subseteq \M$ 
a sub von~Neumann algebra, which is the range of a normal faithful 
conditional expectation $\ep: \M\to \N$. 
Then for every $1<p<\infty$, $L^p(\N)$ is $\cb$-isometrically 
isomorphic to a $\cb$-contractively complemented subspace of $L^p(\M)$.
\end{lem}

Lemma 7.4 implies that the proofs of Theorem~7.2 and Theorem~7.3 can be 
repeated almost word for word to cover the $L^p$-case. 
Note that the argument for $\N_* \oplus \N_*\approx \N_*$ and 
$\M_* \oplus \M_* \approx \M_*$ in the beginning of Theorem~7.2 also works 
for the $L^p$-spaces, when $L^p(\N)$ (resp. $L^p(\M)$) are equipped with 
the natural left $\M$-module structure (resp. left $\N$-module structure).  
Hence we get:

\begin{thm}\label{lem:7.5} 
Let $R_\lambda$, $0<\lambda <1$ and $R_\infty$ be as in Theorem~7.2 and 
let $1\le p<\infty$. 
Then 
\begin{itemize}
\item[(a)] $L^p(R_\lambda) \approx_{\cb} L^p(R_\infty)$. 
\item[(b)] There is an uncountable family of mutually non-isomorphic (in the 
von~Neumann algebra sense) injective type III$_0$-factors on a separable 
Hilbert space, for which $L^p(N_i)\approx_{\cb} L^p(R_\infty)$ 
for all $i\in I$. 
\item[(c)] For every $n\in\nat$, $n\ge2$, $L^p(L(F_n))\approx_{\cb} 
L^p(L(F_\infty))$.
\end{itemize}
\end{thm}

\end{document}